\documentclass{amsart}
\usepackage{epsfig}
\usepackage{graphicx}
\usepackage{color}
\usepackage{amssymb}
\usepackage{esint}
\usepackage{tikz-cd}
\usetikzlibrary{matrix,arrows.meta}

\newcommand\R{\mathbb R}
\newcommand\N{\mathbb N}
\newcommand\C{\mathbb C}

\newcommand\loc{{\text{\upshape loc}}}

\newcommand\const{\text{\upshape Const.}}

\newcommand\Pot{{\mathcal{P}}}
\newcommand\Potc{{\mathcal{P}^c}}
\newcommand\Lagr{{\mathcal{L}}}
\newcommand\Eul{{\mathcal{E}}}
\newcommand\Eulc{{\mathcal{E}^c}}
\newcommand\PcL{{_{\Potc\!\Lagr}}}
\newcommand\EcL{{_{\Eulc\!\Lagr}}}
\newcommand\LPc{{_{\Lagr \Potc}}}
\newcommand\LEc{{_{\Lagr \Eulc}}}
\newcommand\PE{{_{\Pot \Eul}}}
\newcommand\PcEc{{_{\Potc\!\Eulc}}}
\newcommand\EP{{_{\Eul \Pot}}}
\newcommand\EcPc{{_{\Eulc\!\Potc}}}
\newcommand\vvect{\mathbf{v}}
\newcommand\rvect{\mathbf{r}}
\newcommand\svect{\mathbf{s}}
\newcommand\tvect{\mathbf{t}}
\newcommand\phivect{\boldsymbol{\phi}}
\DeclareMathOperator\Div{div}
\DeclareMathOperator\curl{curl}
\DeclareMathOperator\DivGamma{div_{_\Gamma}}
\DeclareMathOperator\Tr{Tr}
\DeclareMathOperator\Rg{Rg}

\numberwithin{equation}{section}
\newtheorem{thm}{Theorem}
\numberwithin{thm}{section}
\newtheorem{cor}[thm]{Corollary}
\newtheorem{lem}[thm]{Lemma}
\newtheorem{prop}[thm]{Proposition}
\theoremstyle{definition}
\newtheorem{definition}[thm]{Definition}
\theoremstyle{remark}
\newtheorem{rem}[thm]{Remark}

\begin{document}
\title[Three evolution problems modelling the interaction...]
{Three evolution problems modelling the interaction between acoustic waves and non--locally reacting surfaces}
\author{Enzo Vitillaro}
\address[E.~Vitillaro]
       {Dipartimento di Matematica e Informatica, Universit\`a di Perugia\\
       Via Vanvitelli,1 06123 Perugia ITALY}
\email{enzo.vitillaro@unipg.it}

\subjclass{35L05, 35L10, 35L20, 35L51, 76N30}

\keywords{Wave equation, hyperbolic systems of second order, acoustic boundary conditions, semigroup theory, acoustic waves, compressible fluids}


\thanks{Supported by the European Union, Next Generation EU}

\begin{abstract} The paper deals with three evolution problems arising in the physical modelling of acoustic phenomena of small amplitude in a fluid, bounded by a surface of extended reaction. The first one is the widely studied wave equation with acoustic boundary conditions, which derivation  from the physical model is not fully mathematically satisfactory. The other two models studied in the paper, in the Lagrangian and Eulerian settings, are physically transparent. In the paper the first model is derived from the other two in a rigorous way, also for solutions merely belonging to the natural energy spaces.
\end{abstract}

\maketitle

\section{Introduction and main results} \label{intro}
\subsection{Presentation of the problem and literature overview}

We deal with three evolution problems arising in the physical modelling of acoustic phenomena of small amplitude in a fluid, bounded by a surface of extended reaction. The problems are posed in a bounded and simply connected domain $\Omega$ of $\R^3$ with boundary $\Gamma=\partial\Omega$ of class
$C^{r,1}$ (see \cite{grisvard}) with  $r=2,3,\ldots,\infty$,
where $\Gamma=\Gamma_0\cup\Gamma_1$, $\overline{\Gamma_0}\cap \overline{\Gamma_1}=\emptyset$, $\Gamma_1$ being nonempty. These properties of $\Omega$, $\Gamma_0$ and $\Gamma_1$ will be assumed throughout the paper without further reference.

The first problem we shall consider, to which in the sequel we shall refer as to  \emph{the Eulerian model}, is the  boundary-value problem
$$(\mathcal{E})\qquad
\begin{cases}
p_t+B\Div \vvect =0 \qquad &\text{in
$\R\times\Omega$,}\\
\rho_0\vvect_t=-\nabla p\qquad &\text{in
$\R\times\Omega$,}\\
\curl\vvect=0\qquad &\text{in
$\R\times\Omega$,}\\
\mu v_{tt}- \DivGamma (\sigma \nabla_\Gamma v)+\delta v_t+\kappa v+p =0\qquad
&\text{on
$\R\times \Gamma_1$,}\\
\vvect\cdot{\boldsymbol{\nu}} =0 \quad \text{on $\R\times \Gamma_0$,}\qquad\vvect\cdot{\boldsymbol{\nu}} =-v_t\qquad
&\text{on
$\R\times \Gamma_1$,}
\end{cases}
$$
where $p=p(t,x)$, $\vvect=\vvect(t,x)$, $v=v(t,y)$, $t\in\R$, $x\in\Omega$, $y\in\Gamma_1$,
$p$ and $v$ taking complex \begin{footnote}{All unknowns and function spaces in the paper will be complex, as common in Acoustics. The reduction to the real case is anyway trivial.}\end{footnote} values while $\vvect$ takes values in $\C^3$.
The reader should remark the difference between the notations $\vvect$ and $v$.
In problem ($\mathcal{E}$) we respectively denote by $\nabla$, $\Div$ and $\curl$
the gradient, divergence and rotor operator, while  $\DivGamma$ and $\nabla_\Gamma$ stand for the
Riemannian version on $\Gamma$ of the first two operators. Moreover by
${\boldsymbol{\nu}}$ we denote the outward normal to $\Omega$, $B$ and $\rho_0$ are fixed positive constants and
$\mu,\sigma,\delta$ and  $\kappa$ are given real  functions on $\Gamma_1$
satisfying the assumption
\begin{itemize}
\item[(A)] $\mu,\sigma,\delta,\kappa\in W^{r-1,\infty}(\Gamma_1)$ with $\min_{\Gamma_1}\mu=\mu_0>0$, $\min_{\Gamma_1}\sigma=\sigma_0>0$,
\end{itemize}
where $W^{\infty,\infty}(\Gamma_1)$ stands for $C^\infty(\Gamma_1)$ when $r=\infty$.
Also assumption (A) will be kept without further reference throughout the paper.

In the paper we shall also deal with \emph{the constrained Eulerian model}, in short  ($\mathcal{E}^c$), constituted by adding to  ($\mathcal{E})$ the  integral condition
\begin{equation}\label{1}
  \int_\Omega p(t,\cdot)=B\int_{\Gamma_1}v(t),\qquad\text{for all $t\in\R$.}
\end{equation}

The complete physical derivation of problems $(\mathcal{E})$ and  $(\mathcal{E}^c)$ was given  in \cite[Chapter~7]{mugnvit}. Here we just recall that $p$ and  $\vvect$ respectively stand for the acoustic excess pressure and the velocity $\mathbf{v}$ field of a fluid contained in $\Omega$, while $v$ stands for the normal displacement inside $\Omega$ of the moving part $\Gamma_1$ of the boundary. In \cite{mugnvit} a mathematical study of $(\Eul)$ and $(\Eulc)$ was not performed. Moreover,
as to the author's knowledge, these problems do not possesses a mathematical literature, although they are strongly related to the third problem we shall study.

The second problem we shall consider, to which in the sequel we shall refer as to \emph{the Lagrangian model}, is the boundary-value problem
$$(\mathcal{L})\qquad
\begin{cases}
\rho_0\rvect_{tt}-B\nabla \Div \rvect=0\qquad &\text{in
$\R\times\Omega$,}\\
\curl\rvect=0\qquad &\text{in
$\R\times\Omega$,}\\
\mu v_{tt}- \DivGamma (\sigma \nabla_\Gamma v)+\delta v_t+\kappa v-B\Div \rvect =0\qquad
&\text{on
$\R\times \Gamma_1$,}\\
\rvect\cdot{\boldsymbol{\nu}} =0 \quad \text{on $\R\times \Gamma_0$,}\qquad\rvect\cdot{\boldsymbol{\nu}} =-v\qquad
&\text{on
$\R\times \Gamma_1$,}
\end{cases}
$$
where we keep the notation and assumptions of the previous problem and $\rvect=\rvect(t,x)$,  $t\in\R$, $x\in\Omega$, takes values in $\C^3$.
Since $\curl\rvect=0$ one could equivalently write, in a distributional sense, equation $(\mathcal{L})_1$ as  $\rho_0\rvect_{tt}-B\Delta \rvect=0$,  where $\Delta\rvect$ stands for the componentwise Laplacian of $\rvect$.
Also the complete physical derivation of problem ($\mathcal{L}$) was given in \cite[Chapter~7]{mugnvit},
 although it was written in the alternative form
$$(\mathcal{L}')\qquad
\begin{cases}
p+B\Div \rvect =0 \qquad &\text{in
$\R\times\Omega$,}\\
\rho_0\rvect_{tt}=-\nabla p\qquad &\text{in
$\R\times\Omega$,}\\
\curl\rvect=0\qquad &\text{in
$\R\times\Omega$,}\\
\mu v_{tt}- \DivGamma (\sigma \nabla_\Gamma v)+\delta v_t+\kappa v+p =0\qquad
&\text{on
$\R\times \Gamma_1$,}\\
\rvect\cdot{\boldsymbol{\nu}} =0 \quad \text{on $\R\times \Gamma_0$,}\qquad\rvect\cdot{\boldsymbol{\nu}} =-v\qquad
&\text{on
$\R\times \Gamma_1$,}
\end{cases}
$$
more similar to problem $(\mathcal{E})$. Since in problem ($\mathcal{L}'$) on can use the first equation to eliminate $p$, the two formulations are completely equivalent. In them $p$ and $v$ have the same physical meaning as in problem $(\mathcal{E})$, altough in a Lagrangian framework, while $\rvect$ stands  for the the displacement of fluid  particles from a reference  position.
In the paper we shall deal only with problem ($\mathcal{L}$).
Also problem ($\mathcal{L}$) does not possess, as to the author's knowledge, a mathematical literature but it is strongly related, as we shall see, with the third problem we shall consider, to which in the sequel we shall refer as to  \emph{the potential model},
is the boundary-value problem
$$(\mathcal{P})\qquad
\begin{cases} \rho_0u_{tt}-B\Delta u=0 \qquad &\text{in
$\R\times\Omega$,}\\
\mu v_{tt}- \DivGamma (\sigma \nabla_\Gamma v)+\delta v_t+\kappa v+\rho u_t =0\qquad
&\text{on
$\R\times \Gamma_1$,}\\
\partial_{\boldsymbol{\nu}} u=0 \quad \text{on $\R\times \Gamma_0$,}\qquad v_t =\partial_{\boldsymbol{\nu}} u\qquad
&\text{on
$\R\times \Gamma_1$,}
\end{cases}
$$
where $u=u(t,x)$,  $t\in\R$, $x\in\Omega$,  $\Delta$ denotes the Laplacian operator, and we keep the notation and assumption introduced for $(\mathcal{E})$
and $(\mathcal{L})$.
In the paper we shall also deal with the
\emph{the constrained potential model}, in short  ($\mathcal{P}^c$), constituted by adding to  ($\mathcal{P})$ the  integral condition
\begin{equation}\label{2}
  \rho_0\int_\Omega u_t(t)=B\int_{\Gamma_1}v(t),\qquad\text{for all $t\in\R$.}
\end{equation}
Unlike ($\mathcal{E}$) and ($\mathcal{L}$), problem ($\mathcal{P}$) possesses a wide mathematical literature, in which the wave equation ($\mathcal{P}$)$_1$ is equivalently written as $u_{tt}-c^2\Delta u=0$, with $c^2=B/\rho_0$, and it is known as wave equation with acoustic boundary conditions.
This type of  boundary conditions have been introduced by Beale and Rosencrans, for bounded or external domains, in \cite{beale, beale2,BealeRosencrans} when $\sigma\equiv 0$ and $\Gamma_0=\emptyset$. In this case the boundary $\Gamma=\Gamma_1$ is called, using the terminology of  \cite[pp.~256]{morseingard}, {\em locally reacting}, since each point of it reacts like an harmonic oscillator.

The same model, when also $\delta=0$, has been proposed in \cite{WeitzKeller} when $\Omega$ is a strip in $\R^2$ or $\R^3$, and in \cite{Peters}
when $\Omega$ is the half-space in $\R^3$, to describe acoustic wave propagation in an ice-covered ocean. We refer the interested reader to~\cite{Belinsky} for a historical overview of these and related problems in Mathematical Physics.

After their introduction acoustic boundary conditions for locally reacting surfaces have been the subject of several papers. See for example \cite{Alcantara, CFL2004,CavFrota, FL2006, FG2000, GGG,  graber, Hao, JGPhD, JGSB2012, KJR2016, KT2008, LLX2018,  Maatoug2017,  mugnolo}.
When one dismisses the simplifying assumption that neighboring points do not interact,  such surfaces, again using the terminology of \cite[p.266]{morseingard}), are called of {\em extended reaction}.
We shall call those which react like a membrane {\em non-locally reacting}, since other types of reactions can be considered.

The simplest case in which $\sigma\equiv\sigma_0$, so   the operator $\DivGamma(\sigma \nabla_\Gamma)$ reduces, up to $\sigma_0$, to the Laplace--Beltrami operator $\Delta_\Gamma$ was briefly considered in \cite[\S 6]{beale} and then studied in \cite{FMV2011,  FMV2014, VF2013bis,VF2016, VF2017} and in the recent paper \cite{becklin2019global}.
In all of them the authors assume that $\Gamma_0\not=\emptyset$ and that the homogeneous Neumann boundary condition on it is replaced by the (mathematically more attracting) homogeneous Dirichlet boundary condition. See also \cite{AMP,ArElst} for somehow related problems.

As to the physical derivation of problem $(\Pot)$, two rough physical derivations of it from the Eulerian and the Lagrangian models were given in \cite[Chapter 7]{mugnvit}.  Although being essentially correct, from the mathematical point of view these derivations were not so clear.

The first aim of the present paper is then to rigorously study the relations among problems  $(\Eul)$, $(\Eul^c)$, $(\Lagr)$, $(\Pot)$ and $(\Pot^c)$. To achieve such a result is mandatory, from a mathematical point of view, to recall a  well--posedness result for $(\Pot)$ and to prove analogous results for the other four problems listed above. This one is the second aim of the paper.  As a byproduct of the analysis we shall also get   optimal regularity results for problems $(\Eul)$, $(\Eul^c)$ and  $(\Lagr)$.
  To illustrate these results in the clearest possible way we start from well--posedness.

In the already quoted paper \cite{mugnvit} the authors widely studied, under more general assumptions and with $c^2=B/\rho_0$, the initial--value problem associated to $(\mathcal{P})$, i.e. the initial--and boundary value problem
$$(\mathcal{P}_0)\qquad
\begin{cases} \rho_0u_{tt}-B\Delta u=0 \qquad &\text{in
$\R\times\Omega$,}\\
\mu v_{tt}- \DivGamma (\sigma \nabla_\Gamma v)+\delta v_t+\kappa v+\rho u_t =0\qquad
&\text{on
$\R\times \Gamma_1$,}\\
\partial_{\boldsymbol{\nu}} u=0 \quad \text{on $\R\times \Gamma_0$,}\qquad v_t =\partial_{\boldsymbol{\nu}} u\qquad
&\text{on
$\R\times \Gamma_1$,}\\
u(0,x)=u_0(x),\quad u_t(0,x)=u_1(x) &
 \text{in $\Omega$,}\\
v(0,x)=v_0(x),\quad v_t(0,x)=v_1(x) &
 \text{on $\Gamma_1$.}
\end{cases}
$$
The starting point of the present paper, see \S~\ref{Section 3} below, will be to recall  several results from \cite{mugnvit} with some simplifications essentially due to the higher regularity of $\Gamma$ assumed here.
In this section we only need to recall the main well--posedness result in the quoted paper.
 To state it we introduce, for  $n\in\N$, $n\le r$, the Hilbert space
\begin{equation}\label{1.3}
{\mathcal H}_\Pot^n=H^n(\Omega)\times H^n(\Gamma_1)\times H^{n-1}(\Omega)\times H^{n-1}(\Gamma_1),
\end{equation}
using the standard notation $H^0=L^2$ which we shall keep throughout the paper.
We also introduce, for $n\in\N\cup\{\infty\}$, $n\le r$, the Fr\'echet  space
\begin{equation}\label{1.4}
X_\Pot^n=\bigcap_{i=0}^n C^i(\R;H^{n-i}(\Omega)\times H^{n-i}(\Gamma_1)).
\end{equation}
The phase space for problem $(\Pot)$ is then $\mathcal{H}_\Pot^1$, and we shall consider weak solutions of $(\Pot)$ and $(\Pot^c)$,  that is couples $(u,v)\in X_\Pot^1$ which satisfy $(\Pot)$ is a suitable distributional sense, see Definition~\ref{Definition 3.1}  below. Moreover we shall mean initial conditions in problem $(\Pot_0)$  in the sense of the space $X_\Pot^1$.
\begin{thm}[\bf Well--posedness for $(\Pot_0)$]\label{Theorem 1.1} For all data  $U_{0\Pot}=(u_0,v_0,u_1,v_1)\in \mathcal{H}_\Pot^1$ problem $(\Pot_0)$ has a unique weak solution $(u,v)\in X_\Pot^1$, continuously depending on  $U_{0\Pot}$ in the topologies of the respective spaces.
Moreover $(u,v)\in X_\Pot^2$ if and only if $U_{0\Pot}\in \mathcal{H}_\Pot^2$, $\partial_{\boldsymbol{\nu}}u_0=0$ on $\Gamma_0$ and $\partial_{\boldsymbol{\nu}}u_0=v_1$ on $\Gamma_1$, such data being dense in $\mathcal{H}_\Pot^1$. In this case $u$ and $v$ satisfy
$(\Pot)_1$ a.e. in $\R\times\Omega$, $(\Pot)_2$--$(\Pot)_3$ a.e. on $\R\times\Gamma_1$ and $(u,v)$ continuously depends  on data in the topologies of $\mathcal{H}_\Pot^2$ and $X_\Pot^2$.
Finally, for all  $U_{0\Pot}\in \mathcal{H}_\Pot^1$ and $s,t\in\R$, $(u,v)$ satisfies the energy identity
\begin{footnote}{Here and in the sequel $|\cdot|_\Gamma$ denotes the norm naturally associated to the  metric naturally induced by the Euclidean metric on the tangent bundle, see \S~\ref{Section 2.2}  below.}\end{footnote}
\begin{equation}\label{1.5}
\int_\Omega \!\!\left(\frac {\rho_0}2|\nabla u|^2\!+ \!\frac {\rho_0^2}{2B} |u_t|^2\right)
+\frac 12 \int_{\Gamma_1}\!\!\!\! \left(\sigma |\nabla_\Gamma v|_\Gamma^2+\mu |v_t|^2 +\kappa|v|^2\right)\Bigg|_s^t\! =\!-\!\int_s^t\!\!\int_{\Gamma_1}\!\!\delta|v_t|^2.
\end{equation}
\end{thm}

\subsection{Main results I: well--posedness.}\label{Section 1.2} We start from problem $(\Pot^c)$. The initial value problem associated to it is obtained by adding the integral condition \eqref{2} to $(\Pot_0)$. In the sequel it will be denoted by $(\Pot^c_0)$.

Setting, for  $n\in\N$, $n\le r$, the closed subspaces (respectively) of ${\mathcal H}_\Pot^n$ and $X_\Pot^n$,
\begin{equation}\label{1.6}
\begin{aligned}
{\mathcal H}_\Potc^n=&\left\{(u,v,w,z)\in {\mathcal H}_\Pot^n: \rho_0\int_\Omega w=B\int_{\Gamma_1}v\right\},\\
 X_\Potc^n=&\{(u,v)\in X_\Pot^n: \text{\eqref{2} holds}\},
 \end{aligned}
\end{equation}
weak solutions of $(\Pot^c)$ will be weak solutions of $(\Pot)$ belonging to $X_\Potc^1$, and the phase space associated to $(\Pot^c)$ is simply ${\mathcal H}_\Potc^1$.

As an easy application of Theorem~\ref{Theorem 1.1}, in \S~\ref{Section 3.4}  we shall prove the following result.
\begin{cor}[\bf Well--posedness for $(\Pot_0^c)$]\label{Corollary 1.1} The statement of Theorem~\ref{Theorem 1.1} continues to hold when one replaces $(\Pot_0)$, ${\mathcal H}_\Pot^n$ and $X_\Pot^n$  with $(\Pot_0^c)$, ${\mathcal H}_\Potc^n$ and $X_\Potc^n$.
\end{cor}
We now turn to the Lagrangian model. The initial--value problem associated to it is the problem
$$(\mathcal{L}_0)\qquad
\begin{cases}
\rho_0\rvect_{tt}-B\nabla \Div \rvect=0\qquad &\text{in
$\R\times\Omega$,}\\
\curl\rvect=0\qquad &\text{in
$\R\times\Omega$,}\\
\mu v_{tt}- \DivGamma (\sigma \nabla_\Gamma v)+\delta v_t+\kappa v-B\Div \rvect =0\qquad
&\text{on
$\R\times \Gamma_1$,}\\
\rvect\cdot{\boldsymbol{\nu}} =0 \quad\text{on $\R\times \Gamma_0$,}\qquad
\rvect\cdot{\boldsymbol{\nu}} =-v\qquad
&\text{on
$\R\times \Gamma_1$,}\\
\rvect(0,x)=\rvect_0(x),\quad \rvect_t(0,x)=\rvect_1(x) &
 \text{in $\Omega$,}\\
v(0,x)=v_0(x),\quad v_t(0,x)=v_1(x) &
 \text{on $\Gamma_1$.}
\end{cases}
$$
Denoting for simplicity $H^n(\Omega)^3=H^n(\Omega;\C^3)$ we introduce,
for $n\in\N_0$, its closed subspace
\begin{equation}\label{1.7}
 H^n_{\curl 0}(\Omega)=\{\rvect\in H^n(\Omega)^3: \curl \rvect =0\}
\end{equation}
where, at least when $n=0$, $\curl\rvect$ is taken in the sense of distributions.
We also set, for $n\in\N_0$, $n\le r$, the Hilbert space
$H^n_\Lagr=H^n_{\curl 0}(\Omega)\times H^n(\Gamma_1)$
and, when $n\in\N$, $n\le r$, its closed subspace
\begin{equation}\label{1.9}
 \mathbb{H}^n_\Lagr=\{(\rvect,v)\in H^n_\Lagr: \rvect\cdot\boldsymbol{\nu}=-v
 \quad\text{on $\Gamma_1$},\quad
 \rvect\cdot\boldsymbol{\nu}=0\quad\text{on $\Gamma_0$}\},
\end{equation}
where $\rvect$  is taken in the trace sense on $\Gamma$. For the same $n$'s we also
set the Hilbert space
\begin{equation}\label{1.10}
 \mathcal{H}^n_\Lagr=\{(\rvect,v,\svect,z)\!\in\! H^n_{\curl 0}(\Omega)\!\times\! H^n(\Gamma_1)\!\times
 \!H^{n-1}_{\curl 0}(\Omega)\!\times\! H^{n-1}(\Gamma_1): (\rvect,v)\in\mathbb{H}^n_\Lagr\},
\end{equation}
endowed with the norm inherited from the product. Moreover, for $n\in\N\cup\{\infty\}$, $n\le r$,
we set the Fr\'echet  space
\begin{equation}\label{1.11}
X_\Lagr^n=\bigcap_{i=0}^{n-1} C^i(\R;\mathbb{H}_\Lagr^{n-i})\cap C^n(\R; H_\Lagr^0).
\end{equation}
The phase space for $(\Lagr)$ is $\mathcal{H}_\Lagr^1$, and in \S~\ref{Section 4.4}
we shall define weak solutions for $(\Lagr)$  as couples $(\rvect,v)\in X_\Lagr^1$
satisfying $(\Lagr)$ is a suitable distributional sense, see Definition~\ref{Definition 4.1}, while we shall mean initial conditions in $(\Lagr_0)$  in the sense of the space $X_\Lagr^1$. The definition of weak solutions, together with their uniqueness,
see Theorem~\ref{Theorem 4.3}  below, is an essential outcome of the paper. In \S~\ref{Section 4.5}  we shall
prove the following result.
\begin{thm}[\bf Well--posedness for $(\Lagr_0)$]\label{Theorem 1.2}
 For all data  $U_{0\Lagr}=(\rvect_0,v_0,\rvect_1,v_1)\in \mathcal{H}_\Lagr^1$
 problem $(\Lagr_0)$ has a unique weak solution $(\rvect,v)\in X_\Lagr^1$,
 continuously depending on  $U_{0\Lagr}$ in the topologies of the respective spaces.

Moreover $(\rvect,v)\in X_\Lagr^2$ if and only if $U_{0\Lagr}\in \mathcal{H}_\Lagr^2$,
and $(\rvect_1,v_1)\in \mathbb{H}_\Lagr^1$, such data being dense in $\mathcal{H}_\Lagr^1$.
In this case $\rvect$ and $v$ satisfy $(\Lagr)_1$ -- $(\Lagr)_2$
a.e. in $\R\times\Omega$, $(\Lagr)_3$--$(\Lagr)_5$ a.e. on $\R\times\Gamma$ and $(\rvect,v)$
continuously depends  on data in the topologies of $\mathcal{H}_\Lagr^2$ and $X_\Lagr^2$.

Finally, for all  $U_{0\Lagr}\in \mathcal{H}_\Lagr^1$ and $s,t\in\R$,
$(\rvect,v)$ satisfies the energy identity
\begin{equation}\label{1.12}
   \int_\Omega \!\!\!\left(\tfrac {\rho_0}2|\rvect_t|^2+ \tfrac B2 |\Div\rvect|^2\right)
+\tfrac 12 \int_{\Gamma_1}\!\!\!\left(\sigma |\nabla_\Gamma v|_\Gamma^2+\mu |v_t|^2
+\kappa|v|^2\right)\Bigg|_s^t\!\! =-\int_s^t\!\!\int_{\Gamma_1}\!\!\delta|v_t|^2.
\end{equation}
\end{thm}
Our second main result concerns the initial--value problems associated to $(\Eul)$ and $(\Eulc)$, the first one being the problem
$$(\mathcal{E}_0)\qquad
\begin{cases}
p_t+B\Div \vvect =0 \qquad &\text{in
$\R\times\Omega$,}\\
\rho_0\vvect_t=-\nabla p\qquad &\text{in
$\R\times\Omega$,}\\
\curl\vvect=0\qquad &\text{in
$\R\times\Omega$,}\\
\mu v_{tt}- \DivGamma (\sigma \nabla_\Gamma v)+\delta v_t+\kappa v+p =0\qquad
&\text{on
$\R\times \Gamma_1$,}\\
\vvect\cdot{\boldsymbol{\nu}} =0 \quad\text{on $\R\times \Gamma_0$,}\qquad \vvect\cdot{\boldsymbol{\nu}} =-v_t\qquad
&\text{on
$\R\times \Gamma_1$,}\\
p(0,x)=p_0(x),\quad \vvect(0,x)=\vvect_0(x) &
 \text{in $\Omega$,}\\
v(0,x)=v_0(x),\quad v_t(0,x)=v_1(x) &
 \text{on $\Gamma_1$,}
\end{cases}
$$
while the second one is simply the problem obtained by adding the integral condition \eqref{1} to $(\Eul_0)$, to which in the sequel we shall refer as problem $(\Eul^c_0)$. To deal with them we set, for $n\in\N$, $n\le r$, the Hilbert spaces
\begin{equation}\label{1.13}
\begin{aligned}
  &\mathcal{H}_\Eul^n=H^{n-1}(\Omega)\!\times\! H^{n-1}_{\curl 0}(\Omega)\!\times\! H^n(\Gamma_1)\!\times \!H^{n-1}(\Gamma_1),\\
   &\mathcal{H}_\Eulc^n=\left\{(p,\vvect,v,z)\in \mathcal{H}_\Eul^n: \int_\Omega p=B\int_{\Gamma_1}v\right\},
   \end{aligned}
\end{equation}
and, for $n\in\N\cup\{\infty\}$, $n\le r$, the Fr\'echet spaces
$$X_\Eul^n\!=\!\bigcap_{i=0}^{n-1}\! C^i(\R;H^{n-1-i}(\Omega))\!\times\! \bigcap_{i=0}^{n-1} \!C^i(\R;H^{n-1-i}_{\curl 0}(\Omega))\!\times\! \bigcap_{i=0}^n \!C^i(\R;H^{n-i}(\Gamma_1))
$$
and $X_\Eulc^n=\{(p,\vvect,v)\in X_\Eul^n: \text{\eqref{1} holds}\}.$
The phase space for $(\Eul)$ and  $(\Eulc)$ are $\mathcal{H}_\Eul^1$ and $\mathcal{H}_\Eulc^1$, and in \S~\ref{Section 5.4}
we shall define weak solutions of $(\Eul)$ (respectively of $(\Eulc)$) as triples $(p,\vvect,v)$ in $X_\Eul^1$ (in $X_\Eulc^1$) satisfying $(\Eul)$ is a suitable distributional sense, see Definition~\ref{Definition 5.1}, while we shall mean initial conditions in $(\Eul_0)$ and $(\Eul^c_0)$ in the sense of the space $X_\Eul^1$. The definition of weak solutions, together with their uniqueness,
see Theorem~\ref{Theorem 5.3}  below, is an essential outcome of the paper. In \S~\ref{Section 5.5}  we shall
then prove the following result.
\begin{thm}[\bf Well--posedness for $(\Eul_0)$ and $(\Eul^c_0)$]\label{Theorem 1.3}
 For all data  $U_{0\Eul}\in\mathcal{H}_\Eul^1$, $U_{0\Eul}=(p_0,\vvect_0,v_0,v_1)$,
 problem $(\Eul_0)$ has a unique weak solution $(p,\vvect,v)\in X_\Eul^1$,
 continuously depending on  $U_{0\Eul}$ in the topologies of the respective spaces.
Moreover $(p,\vvect,v)\in X_\Eul^2$ if and only if $U_{0\Eul}\in \mathcal{H}_\Eul^2$,
and $(\vvect_0,v_1)\in \mathbb{H}_\Lagr^1$, such data being dense in $\mathcal{H}_\Eul^1$.
In this case $p$, $\vvect$ and $v$ satisfy $(\Eul)_1$ -- $(\Eul)_3$
a.e. in $\R\times\Omega$, $(\Eul)_4$--$(\Eul)_6$ a.e. on $\R\times\Gamma$ and $(p,\vvect,v)$
continuously depends  on data in the topologies of $\mathcal{H}_\Eul^2$ and $X_\Eul^2$.
For all  $U_{0\Eul}\in \mathcal{H}_\Eul^1$ and $s,t\in\R$,
$(p,\vvect,v)$ satisfies the energy identity
\begin{equation}\label{1.15}
  \int_\Omega \!\!\left(\tfrac {\rho_0}2|\vvect|^2+ \tfrac 1{2B} |p|^2\right)
+\tfrac 12 \int_{\Gamma_1}\!\!\left(\sigma |\nabla_\Gamma v|_\Gamma^2+\mu |v_t|^2
+\kappa|v|^2\right)\Bigg|_s^t =-\int_s^t\int_{\Gamma_1}\delta|v_t|^2.
\end{equation}
Finally all previous assertions continue to hold when replacing $(\Eul_0)$, $\mathcal{E}_\Eul^n$ and $X_\Eul^n$, $n=1,2$, with $(\Eul^c_0)$, $\mathcal{E}_\Eulc^n$ and $X_\Eulc^n$.
\end{thm}

\subsection{Main results II: relations.} Theorems~\ref{Theorem 1.1}--\ref{Theorem 1.3} and Corollary~\ref{Corollary 1.1} look like variants of a single result. This fact depends on the relations between all problems treated. To illustrate them  it is useful to disregard initial conditions and introduce, for  $n\in\N\cup\{\infty\}$, $n\le r$,  the Fr\'echet spaces
\begin{equation}\label{1.18}
\begin{aligned}
&\mathcal{S}^n_\Pot =\{(u,v)\in X_\Pot^n: (u,v)\text{ is a weak solution of }(\Pot)\},\\
& \mathcal{S}^n_\Eul =\{(p,\vvect,v)\in X_\Eul^n: (p,\vvect,v)\text{ is a weak solution of }(\Eul)\},\\
&\mathcal{S}^n_\Lagr =\{(\rvect,v)\in X_\Lagr^n: (\rvect,v)\text{ is a weak solution of }(\Lagr)\},\\
&\mathcal{S}^n_\Potc=\mathcal{S}^n_\Pot\cap X_\Potc^1, \qquad \mathcal{S}^n_\Eulc=\mathcal{S}^n_\Eul\cap X_\Eulc^1,
\end{aligned}
\end{equation}
endowed with the topologies of the ambient spaces, the last two spaces clearly denoting  spaces of weak solutions of $(\Potc)$ and $(\Eulc)$. To relate solutions of $(\Pot)$ and  $(\Potc)$ to solutions of the other models we have to consider $u$ up to functions which are constant in $\R\times\Omega$. Denoting by  $\C_{\R\times\Omega}$ the space they constitute, we thus  have to consider solutions $(u,v)$ up to $\C_{X_\Pot}=\C_{\R\times\Omega}\times\{0\}$, the elements of which are trivial solutions of $(\Potc)$ and hence of $(\Pot)$. We then introduce the standard quotient Fr\'echet spaces (see \cite[Chapter I, p.31]{Schaefer})
$\dot{\mathcal{S}}^n_\Pot:=\mathcal{S}^n_\Pot/\C_{X_\Pot}$ and  $\dot{\mathcal{S}}^n_\Potc:=\mathcal{S}^n_\Potc/\C_{X_\Pot}$,
so trivially $\dot{\mathcal{S}}^n_\Potc\subset \dot{\mathcal{S}}^n_\Pot$,
which elements are of the form $(u,v)+\C_{X_\Pot}$ for $(u,v)$ in the respective class.

Denoting by $\mathcal{L}(X;Y)$ the space of continuous linear operators between two Fr\'echet spaces $X$ and $Y$, and $\mathcal{L}(X)=\mathcal{L}(X;X)$, we can state our third main result.
\begin{thm}[\bf The relation between $(\Potc)$ and  $(\Lagr)$]\label{Theorem 1.4}
For any weak solution $(u,v)$ of $(\Potc)$ and any $t\in\R$ there is unique $\rvect(t)\in H^1(\Omega)^3$ solving the problem
\begin{equation}\label{1.21}
  \begin{cases}
  -B\Div \rvect(t)=\rho_0u_t(t),\quad &\text{in $\Omega$,}\\
  \quad\curl\rvect(t)=0\quad &\text{in $\Omega$,}\\
\quad \rvect(t)\cdot\boldsymbol{\nu}=-v(t)\quad &\text{on $\Gamma_1$,}\\
\quad \rvect(t)\cdot\boldsymbol{\nu}=0\quad &\text{on $\Gamma_0$,}
  \end{cases}
\end{equation}
and the couple $(\rvect,v)\in X_\Lagr^1$ is a weak solution of $(\Lagr)$ which also satisfies
\begin{equation}\label{1.22}
 \rvect_t(t)=-\nabla u(t)\qquad\text{for all $t\in\R$.}
\end{equation}
Moreover the map $(u,v)\mapsto (\rvect,v)$ defines a surjective operator $\Psi_\PcL\in\mathcal{L}(\mathcal{S}^1_\Potc;\mathcal{S}^1_\Lagr)$ such that $ \text{Ker } \Psi_\PcL= \C_{X_\Pot},$
so $\Psi_\PcL$ subordinates a bijective isomorphism $\dot{\Psi}_\PcL\in\mathcal{L}(\dot{\mathcal{S}}^1_\Potc;\mathcal{S}^1_\Lagr)$ given by $\dot{\Psi}_\PcL [(u,v)+\C_{X_\Pot}]=\Psi_\PcL (u,v)$ for all $(u,v)\in\mathcal{S}_\Potc^1$.

Its inverse is the operator $\dot{\Psi}_\LPc\in\mathcal{L}(\mathcal{S}^1_\Lagr;\dot{\mathcal{S}}^1_\Potc)$,
$\dot{\Psi}_\LPc (\rvect,v)=(u,v)+\C_{X_\Pot}$ for all $(\rvect,v)\in\mathcal{S}_\Lagr^1$,
defined by setting $u$, up to a space--time constant, by
\begin{equation}\label{1.26}
  u(t)=u(0)-\tfrac B{\rho_0}\int_0^t \Div \rvect(\tau)\,d\tau,\quad\text{and}\quad -\nabla u(0)=\rvect_t(0),
\end{equation}
so $u$ also satisfies \eqref{1.22}.
Finally, when  $(u,v)$ and $(\rvect,v)$ are the weak solutions of $(\Pot^c_0)$ and $(\Lagr_0)$ in Corollary~\ref{Corollary 1.1} and Theorem~\ref{Theorem 1.2}, we have
\begin{equation}\label{1.27}
  (\rvect,v)=\Psi_\PcL(u,v)\qquad\Longleftrightarrow\qquad  -\nabla u_0=\rvect_1\,\,\text{and}\,\, -B\Div\rvect_0=\rho_0u_1.
\end{equation}
\end{thm}
\begin{rem}\label{Remark 1.1}
Since, in the equivalent version $(\Lagr')$ of the Lagrangian model one has $p=-B\Div\rvect$, equations  \eqref{1.21}--\eqref{1.22} show that the relation established in Theorem~\ref{Theorem 1.4} fulfills the equations $-\nabla u=\rvect_t$ and $p=\rho_0 u_t$, hence giving a  derivation of problem $(\Potc)$ mathematically clearer than the one given in \cite{mugnvit}. Moreover \eqref{1.26} makes explicit the construction of the velocity potential $u$.
\end{rem}
We can now give our fourth main result.
\begin{thm}[\bf The relations between $(\Pot)$ and $(\Eul)$, $(\Potc)$ and $(\Eulc)$ ]\label{Theorem 1.5}
For any weak solution $(u,v)$ of $(\Pot)$ (of  $(\Potc)$) the triple
\begin{equation}\label{1.28}
(p,\vvect,v)=(\rho_0 u_t,-\nabla u,v)
\end{equation}
is a weak solution of $(\Eul)$ (respectively of  $(\Eulc)$). Moreover the map $(u,v)\mapsto (p,\vvect,v)$ defines a surjective operator $\Psi_\PE\mathcal{L}(\mathcal{S}^1_\Pot;\mathcal{S}^1_\Eul)$, which  restricts to a surjective operator $\Psi_\PcEc\mathcal{L}(\mathcal{S}^1_\Potc;\mathcal{S}^1_\Eulc)$, and
$\text{Ker } \Psi_\PE = \text{Ker } \Psi_\PcEc =\C_{X_\Pot}$,
so $\Psi_\PE$ and $\Psi_\PcEc$  subordinate bijective isomorphisms
\begin{equation}\label{1.30}
\begin{alignedat}2
  &\dot{\Psi}_\PE\in\mathcal{L}(\dot{\mathcal{S}}^1_\Pot;\mathcal{S}^1_\Eul),\quad
&&\dot{\Psi}_\PE [(u,v)+\C_{X_\Pot}]=\Psi_\PE (u,v)\,\, \text{for all $(u,v)\in\mathcal{S}_\Pot^1$,}\\
& \dot{\Psi}_\PcEc\in\mathcal{L}(\dot{\mathcal{S}}^1_\Potc;\mathcal{S}^1_\Eulc),\quad
&&\dot{\Psi}_\PcEc =(\dot{\Psi}_\PcEc)_{|\dot{\mathcal{S}}_\Potc^1}.
\end{alignedat}
\end{equation}
The inverse of $\dot{\Psi}_\PE$ is the operator
\begin{equation}\label{1.31}
 \dot{\Psi}_\EP\in \mathcal{L}(\mathcal{S}^1_\Eul;\dot{\mathcal{S}}^1_\Pot),\quad
\dot{\Psi}_\EP (p,\vvect,v)=(u,v)+\C_{X_\Pot}\,\, \text{for all $(p,\vvect,v)\in\mathcal{S}_\Eul^1$,}
\end{equation}
defined by setting $u$, up to a space--time constant, by
\begin{equation}\label{1.32}
  u(t)=u(0)-\tfrac 1{\rho_0}\int_0^t p(\tau)\,d\tau,\quad\text{and}\quad -\nabla u(0)=\vvect(0),
\end{equation}
so $u$ also satisfies
\begin{equation}\label{1.33}
 -\nabla u(t)=\vvect(t)\qquad\text{for all $t\in\R$.}
\end{equation}
The inverse of $\dot{\Psi}_\PcEc$ is the operator
\begin{equation}\label{1.34}
 \dot{\Psi}_\EcPc\in\mathcal{L}(\dot{\mathcal{S}}^1_\Eulc;\mathcal{S}^1_\Potc),\quad
\dot{\Psi}_\EcPc =(\dot{\Psi}_\EP)_{|\dot{\mathcal{S}}_\Eulc^1}.
\end{equation}
Finally, when  $(u,v)$ and $(p,\vvect,v)$ are the weak solutions of $(\Pot_0)$ and $(\Eul_0)$,  or of $(\Pot^c_0)$ and $(\Eul^c_0)$ in Theorem~\ref{Theorem 1.1}, Corollary~\ref{Corollary 1.1} and Theorem~\ref{Theorem 1.3}, we have
\begin{equation}\label{1.35}
  (p,\vvect,v)=\Psi_\PE(u,v)\qquad\Longleftrightarrow\qquad  -\nabla u_0=\vvect_0\,\,\text{and}\,\, p_0=\rho_0u_1.
\end{equation}
\end{thm}
\begin{rem}\label{Remark 1.2}
Equations  \eqref{1.32} -- \eqref{1.33} show that the relation established in Theorem~\ref{Theorem 1.5} fulfills the equations $-\nabla u=\vvect$ and $\rho_0 u_t=p$, giving a  mathematical derivation of problems $(\Pot)$ and  $(\Potc)$ which is clearer than the one given in \cite{mugnvit}. They  also make explicit the construction of the velocity potential $u$.
\end{rem}
By simply combining Theorems~\ref{Theorem 1.4} and \ref{Theorem 1.5}, see \S~\ref{Section 5.6} below, one gets the relation between the Lagrangian model and the constrained Eulerian one.
\begin{cor}[\bf The relation between $(\Eulc)$ and  $(\Lagr)$]\label{Corollary 1.2}
For any weak solution $(p,\vvect,v)$ of $(\Eulc)$ and any $t\in\R$ there is unique $\rvect(t)\in H^1(\Omega)^3$ solving the problem
\begin{equation}\label{1.36}
  \begin{cases}
  -B\Div \rvect(t)=p(t),\quad &\text{in $\Omega$,}\\
  \quad \curl\rvect(t)=0\quad &\text{in $\Omega$,}\\
\quad \rvect(t)\cdot\boldsymbol{\nu}=-v(t)\quad &\text{on $\Gamma_1$,}\\
\quad \rvect(t)\cdot\boldsymbol{\nu}=0\quad &\text{on $\Gamma_0$,}
  \end{cases}
\end{equation}
and the couple $(\rvect,v)\in X_\Lagr^1$ is a weak solution of $(\Lagr)$ which also satisfies
\begin{equation}\label{1.37}
 \rvect_t(t)=\vvect(t)\qquad\text{for all $t\in\R$.}
\end{equation}
Hence $\rvect$ can also be expressed as
\begin{equation}\label{5.52}
 \rvect(t)=\rvect(0)+\int_0^t \vvect(\tau)\,d\tau,\quad\text{where $\rvect(0)$ solves problem \eqref{1.36} for $t=0$.}
\end{equation}
Moreover the map $(p,\vvect,v)\mapsto (\rvect,v)$ defines a bijective isomorphism
\begin{equation}\label{1.38}
\Psi_\EcL=\dot{\Psi}_\PcL\cdot \dot{\Psi}_\EcPc\in\mathcal{L}(\mathcal{S}^1_\Eulc;\mathcal{S}^1_\Lagr)
\end{equation}
which inverse $\Psi_\LEc=\dot{\Psi}_\PcEc\cdot \dot{\Psi}_\LPc\in\mathcal{L}(\mathcal{S}^1_\Lagr;\mathcal{S}^1_\Eulc)$ is given by
\begin{equation}\label{1.39}
\Psi_\LEc(\rvect,v)=(p,\vvect,v)=(-B\Div\rvect,\rvect_t,v).
\end{equation}
Finally, when $(\rvect,v)$ and $(p,\vvect,v)$ respectively denote the solutions of the Cauchy problems $(\Lagr_0)$ and $(\Eul_0^c)$ given by Theorems~\ref{Theorem 1.2} and \ref{Theorem 1.3} we have
\begin{equation}\label{1.40}
  (\rvect,v)=\Psi_\EcL(p,\vvect,v)\qquad\Longleftrightarrow\qquad  -B\Div \rvect_0=p_0\,\,\text{and}\,\, \rvect_1=\vvect_0.
\end{equation}
\end{cor}
The relations between the isomorphisms obtained in Theorems~\ref{Theorem 1.4}--\ref{Theorem 1.5} and in Corollary~\ref{Corollary 1.2} are illustrated by the  commutative diagrams
\tikzcdset{row sep/normal=2.5em}
\tikzcdset{column sep/normal=8em}
\begin{equation}\label{1.41}
\begin{tikzcd}
[column sep=large]
\dot{\mathcal{S}}^1_\Pot\arrow[r,shift left,"\dot{\Psi}_\PE"]\arrow[r,leftarrow,shift right,"\dot{\Psi}_\EP"']&\mathcal{S}^1_\Eul
\end{tikzcd}
\qquad\qquad
\begin{tikzcd}
\dot{\mathcal{S}}^1_\Potc\arrow[r,shift left,"\dot{\Psi}_\PcEc"]\arrow[r,leftarrow,shift right,"\dot{\Psi}_\EcPc"']
\arrow[d,shift left,"\dot{\Psi}_\LPc"]\arrow[d,leftarrow,shift right,"\dot{\Psi}_\PcL"']&\mathcal{S}^1_\Eulc\\
\mathcal{S}^1_\Lagr\arrow[ru,shift left,"\Psi_\LPc"]\arrow[ru,leftarrow,shift right,"\Psi_\PcL"']&
\end{tikzcd}
\end{equation}
which also shows that a Lagrangian counterpart of problems $(\Pot)$,  $(\Lagr)$  is missing.
These relations may probably be best understood at an abstract level by considering the groups of linear operators associated to the various problems considered, see Theorems~\ref{Theorem 4.1}, \ref{Theorem 5.1} and \ref{Theorem 5.5} below.

\noindent{\bf Organization of the paper.} The organization of the paper is simple: in \S~\ref{Section 2} we give all preliminaries needed in the paper while  \S~\ref{Section 3},  \S~\ref{Section 4} and \S~\ref{Section 5} are respectively devoted to the potential, the Lagrangian and the Eulerian models.

\section{Notation and preliminaries}\label{Section 2}
\subsection{Notation}\label{Section 2.1}
We shall denote $\widetilde{\N}=\N\cup\{\infty\}$. Borrowing a convention in Physics, for vectors in $\C^3$ and vector--valued functions we shall use boldface. For $\mathbf{x}=(x_1,x_2,x_3), \mathbf{y}=(y_1,y_2,y_3)\in \C^3$ we shall denote $\mathbf{x}\cdot \mathbf{y}=\sum_{i=1}^3 x_iy_i$ and by $\overline{\mathbf{x}}$ the vector conjugated to $\mathbf{x}$.

We shall use the standard notation  for functions spaces on $\Omega$, referring to \cite{adamsfournier}. As already done for formula \eqref{1.7}, where  $H^n(\Omega)^3=H^n(\Omega;\C^3)$, to simplify the notation we shall systematically identify the $\C^3$--valued versions of all spaces above with the Cartesian cubes of the corresponding scalar spaces. Moreover $\|\cdot\|_p$, $1\le p\le \infty$, will denote the norm in $L^p(\Omega)$ and in $L^p(\Omega)^3$, since no confusion will arise.

Moreover, for any Fr\'{e}chet space $X$  we shall denote by $X'$ its dual, by $\langle\cdot,\cdot\rangle_X$ the duality product and we shall use the standard notation for Bochner--Lebesgue and Sobolev spaces of $X$--valued functions when $X$ is a Banach space.
\subsection{Function spaces and operators on $\Gamma$}\label{Section 2.2}The assumption made on $\Omega$, $\Gamma_0$ and $\Gamma_1$ assures that $\Gamma$ inherits from $\R^3$ the structure of a Riemannian surface of class $C^r$,  so in the sequel we shall use some notation of geometric nature, quite common when $r=\infty$, see \cite{taylor}, which can be easily extended to $r<\infty$, see for example \cite{mugnvit} or \cite{Dresda1,Dresda2}.

Moreover, since $\overline{\Gamma_0}\cap \overline{\Gamma_1}=\emptyset$, both $\Gamma_0$ and $\Gamma_1$ are relatively open on $\Gamma$. Hence all geometrical concepts apply to them as well and, to avoid repetition, in the sequel we shall denote by $\Gamma'$ any relatively open subset of $\Gamma$.

We shall denote by $(\cdot,\cdot)_\Gamma$ the Riemannian metric inherited from $\R^3$ and uniquely extended to an Hermitian metric on the complexified  tangent bundle $T(\Gamma')$, and also the associated bundles metric on the complexified cotangent bundle $T^*(\Gamma')$. By $|\cdot|_\Gamma^2= (\cdot,\cdot)_\Gamma$  we shall denote the associated bundle norms.

The standard surface elements $\omega$ associated to $(\cdot,\cdot)_\Gamma$ is then the density of the Lebesgue surface measure on $\Gamma$, coinciding with the restriction to $\Gamma$ of the Hausdorff measure $\mathcal{H}^2$, i.e. $\omega=d\mathcal{H}^2$. In the sequel $\Gamma$ will be equipped without further comments  with this measure, with the corresponding  notions of a.e. equivalence, integrals and Lebesgue spaces $L^p(\Gamma')$, $1\le p\le \infty$. For simplicity we shall denote $\|\cdot\|_{p,\Gamma'}=\|\cdot\|_{L^p(\Gamma')}$. Moreover the notation $d\mathcal{H}^2$ will be dropped from boundary integrals, and a.e. equivalence on $\R\times\Gamma$ will be referred to the Hausdorff measure $\mathcal{H}^3$ in $\R^4$.

Sobolev spaces on $\Gamma'$ are treated in many textbooks when $r=\infty$, see for example
\cite{hebey,lionsmagenes1}. The case $r<\infty$  is treated in  \cite{grisvard} and, when $\Gamma$ is possibly non--compact, in \cite{mugnvit}. Here we shall refer, for simplicity, to \cite{grisvard}, and we shall use the standard notation.
Moreover, mainly to simplify the notation, since $\overline{\Gamma_0}\cap \overline{\Gamma_1}=\emptyset$, by identifying the elements of $W^{s,q}(\Gamma_i)$, $i=0,1$,  with their trivial extension to $\Gamma$, we have the splitting
\begin{equation}\label{2.8}
W^{s,q}(\Gamma)=W^{s,q}(\Gamma_0)\oplus W^{s,q}(\Gamma_1)\qquad\text{for $s\in\R$, $-r\le s\le r$, $1\le q<\infty$,}
\end{equation}

We refer to \cite{mugnvit} for details on the  Riemannian gradient operator $\nabla_\Gamma$ and on the
Riemannian divergence operator $\Div_\Gamma$. Here we just recall that one gets the operator
\begin{equation}\label{2.9}
 \Div_\Gamma(\sigma\nabla_\Gamma)\in\mathcal{L}\left(H^m(\Gamma_1),H^{m-2}(\Gamma_1)\right)\qquad\text{for $m\in\N$, $2\le m\le r$,}
\end{equation}
and that, $\Gamma_1$ being compact, one gets that
\begin{equation}\label{2.10}
\int_{\Gamma_1} -\Div_\Gamma(\sigma\nabla_\Gamma u)v=\int_{\Gamma_1}\sigma (\nabla_\Gamma u,\nabla_\Gamma \overline{v})_\Gamma
\end{equation}
for all  $u\in H^2(\Gamma_1)$, $v\in H^1(\Gamma_1)$.
Finally, in the sequel we shall use the well--known Trace Theorem, i.e. the  existence of the trace operator
$\Tr\in \mathcal{L}\left(H^m(\Omega),H^{m-1/2}(\Gamma)\right)$, for $m\in\N$, $m\le r$.
We shall denote, as usual, $\Tr u=u_{|\Gamma}$. By $u_{|\Gamma_0}$ and $u_{|\Gamma_1}$ we shall denote the restrictions of $u_{|\Gamma}$ to $\Gamma_0$ and $\Gamma_1$ and, when clear, we shall omit trace related subscripts.
\subsection{Quotient spaces and projections}\label{Section 2.3}Denoting by $\C_\Omega$ the space of constant functions in $\Omega$, since $\Omega$ is bounded we can introduce the quotient Hilbert and Fr\'{e}chet spaces
\begin{footnote}{Unfortunately there is no a standard notation for them in the literature.}
\end{footnote}
$\dot{H}^n(\Omega)=H^n(\Omega)/\C_\Omega$, $n\in\N$, and $\dot{C}^\infty(\overline{\Omega})=C^\infty(\overline{\Omega})/\C_\Omega$,
endowed with the standard quotient norm and pseudo-norms (see \cite{Schaefer}).
We denote by $\sim$ the equivalence relation defined on $H^1(\Omega)$ by $u\sim v$ if and only if $u-v\in\C_\Omega$ and by
$u +\C_\Omega $ the equivalence class of $u$ with respect to $\sim$. Consequently, denoting $\pi_0 u=u+\C_\Omega$, we have
$\pi_0 \in\mathcal{L}\left(H^n(\Omega); \dot{H}^n(\Omega)\right)\cap \mathcal{L}\left(C^\infty(\overline{\Omega}); \dot{C}^\infty(\overline{\Omega})\right)$ for all $n\in\N$.
We make the reader aware that in the last formula and (when useful) in the sequel we shall implicitly restrict linear operators.
\subsection{Basic notions of semigroup theory} We now recall, since in the sequel we shall repetitively use them,  some basic facts in semigroup theory. Given any unbounded linear operator $A:D(A)\subset H\to H$ on the Hilbert space $H$,
we shall consider {\em strong} (or { \em classical} ) and {\em generalized} (or {\em mild}) solutions of  the equation $U'+AU=0$\ in $H$, and  of the Cauchy problem
\begin{equation}
\label{2.15}
U'+AU=0\qquad\text{in $H$}, \qquad U(0)=U_0\in H,
\end{equation}
in the standard semigroup sense, see \cite{engelnagel}.
Moreover, when $-A:D(A)\subset H\to H$ is the generator of a strongly continous semigroup on the Hilbert space $H$, trivially adapting
\cite[Chapter II, pp. 124--125]{engelnagel} to semigroups of arbitrary growth bound, one can inductively set, for $n\in\N_0$, the Hilbert space
\begin{equation}\label{2.16}
D(A^0)=H,\qquad D(A^n)=\{U\in D(A^{n-1}): AU\in D(A^{n-1})\},
\end{equation}
endowed with the graph norm
and, for $n\in\N$, the part ${_n}A$
of $A$ in $D(A^n)$, that is the unbounded operator
\begin{equation}\label{2.18}
{_n}A: D({_n}A)\subset D(A^n)\to D(A^n), \quad D({_n}A)=D(A^{n+1}),\quad {_n}A= A_{|D(A^n)}.
\end{equation}
We also set, for $n\in\N$, the Fr\'{e}chet space
$Y^n=\bigcap_{i=0}^{n-1} C^i(\R; D(A^{n-1-i}))$.
By standard semigroup theory, one then gets the following result.
\begin{lem}\label{Lemma 2.1}Let $-A:D(A)\subset H\to H$ be the generator of a strongly continous group $\{T^1(t),t\in\R\}$ on the Hilbert space $H$. Then
\renewcommand{\labelenumi}{{\roman{enumi})}}
\begin{enumerate}
\item for all $U_0\in H$ problem \eqref{2.15} has a unique generalized solution $U\in Y^1$ given by $U(t)=T^1(t)[U_0]$ for all $t\in\R$, continously depending on $U_0$ in the topologies of the respective spaces;
\item for all $n\in\N$ the operator ${_n}A$ generates on $D(A^n)$ the strongly continous group $\{T^{n+1}(t),t\in\R\}$, given by
$T^{n+1}(t)=T^1(t)_{|D(A^n)}$ for all $t\in\R$;
\item the following three facts are equivalent: $U$ is a strong solution of \eqref{2.15},   $U_0\in D(A)$,  $U\in Y^2$;
\item for all $n\in\N$ we have $U\in Y^n$ if and only if $U_0\in D(A^{n-1})$ and, in this case, $U$ continously depends on $U_0$ in the topologies of the respective spaces.
\end{enumerate}
\end{lem}
\section{The potential models}\label{Section 3}
\subsection{Functional spaces} In addition to the main phase spaces ${\mathcal H}_\Pot^n$ and ${\mathcal H}_\Potc^n$ introduced in \S~\ref{intro} we shall also use (see \S~\ref{Section 2.3}) for $n\in\N$, $n\le r$, the Hilbert spaces
\begin{align*}
\dot{\mathcal H}_\Pot^n=& \dot{H}^n(\Omega)\times H^n(\Gamma_1)\times H^{n-1}(\Omega)\times H^{n-1}(\Gamma_1),\\
\dot{\mathcal H}_\Potc^n=& \left\{(u,v,w,z)\in \dot{\mathcal H}_\Pot^n: \rho_0\int_\Omega w=B\int_{\Gamma_1}v\right\},
\end{align*}
and
\begin{equation}\label{3.2}
H_\Potc^n=\left\{(v,w)\in H^n(\Gamma_1)\times H^{n-1}(\Omega): \rho_0\int_\Omega w=B\int_{\Gamma_1}v\right\}.
\end{equation}
When needed we shall use the trivial identification $\dot{\mathcal H}_\Potc^n=\dot{H}^n(\Omega)\times H_\Potc^n\times H^{n-1}(\Gamma_1)$.
We recall the compatibility conditions for problem $(\Pot_0)$ with $c^2=B/\rho_0$ (see \cite[Theorem 1.2.3]{mugnvit}) given for $n\in\N$, $2\le n\le r$ and data $U_0=(u_0,v_0, u_1,v_1)\in \mathcal{H}_\Pot^n $ by the equations
\begin{equation}\label{3.4}\left\{
\begin{aligned}
&\partial_{\boldsymbol{\nu}} \Delta^iu_0=0,\qquad \text{on $\Gamma_0$ \qquad for $i=0,\ldots,   \lfloor n/2 \rfloor -1$,}\\
&\partial_{\boldsymbol{\nu}} \Delta^iu_1=0,\qquad \text{on $\Gamma_0$ \qquad for $i=0,\ldots, \lfloor (n-1)/2\rfloor -1$, when $n\ge 3$,}\\
&\partial_{\boldsymbol{\nu}} u_0=v_1,\qquad\,\,\,\, \text{on $\Gamma_1$,}\\
&\mu \partial_{\boldsymbol{\nu}} u_1=\DivGamma (\sigma\nabla_\Gamma v_0)-\delta \partial_{\boldsymbol{\nu}} u_0-\kappa v_0-\rho_0 u_1,\qquad \text{on $\Gamma_1$, \quad when $n\ge 3$,}\\
&\begin{split}
B\mu \partial_{\boldsymbol{\nu}} \Delta^i\!u_0\!=\!\rho_0[\DivGamma (\sigma\nabla_\Gamma \partial_{\boldsymbol{\nu}}\Delta^{i-1}\!u_0)\!
-\!\delta \partial_{\boldsymbol{\nu}}\Delta^{i-1}\!u_1\!-\!\kappa\partial_{\boldsymbol{\nu}}\Delta^{i-1}\!u_0\!-\! B\Delta^i \!u_0]\\
\text{on $\Gamma_1$  \quad for $i=1,\ldots, \lfloor n/2 \rfloor -1$, \quad when $n\ge 4$,}\qquad\,\,\,&
\end{split}\\
&\begin{split}
B\mu \partial_{\boldsymbol{\nu}} \Delta^i\!u_1\!=-\!B\delta \partial_{\boldsymbol{\nu}}\Delta^iu_0\!\!+\!\!\rho_0[\DivGamma (\sigma\nabla_\Gamma \partial_{\boldsymbol{\nu}}\Delta^{i-1}\!u_1)\!
\!-\!\kappa\partial_{\boldsymbol{\nu}}\Delta^{i-1}\!u_1\!-\! B\Delta^i \!u_1]\\
\text{on $\Gamma_1$  \quad for $i=1,\ldots, \lfloor (n-1)/2 \rfloor -1$, \quad when $n\ge 5$.}\,\,\,&
\end{split}
\end{aligned}\right.
\end{equation}
In \eqref{3.4} and in the sequel  $\lfloor\cdot\rfloor$ stands for the integer part. Since the equations \eqref{3.4} make sense also for data $U_0\in \dot{\mathcal{H}}_\Pot^n $, i.e. for $u_0\in \dot{H}^n(\Omega)$, we can set for $n\in\N$, $2\le n\le r$,
\begin{equation}\label{3.5}
D^{n-1}_\Pot=\{U_0\in \mathcal{H}_\Pot^n: \text{\eqref{3.4} hold}\},\quad\text{and}\quad \dot{D}^{n-1}_\Pot=\{U_0\in \dot{\mathcal{H}}_\Pot^n: \text{\eqref{3.4} hold}\},
\end{equation}
which by \eqref{2.9} are closed subspaces of ${\mathcal H}_\Pot^n$ and $\dot{\mathcal H}_\Pot^n$, together with their subspaces
\begin{equation}\label{3.6}
D^{n-1}_\Potc=D^{n-1}_\Pot\cap \mathcal{H}_\Potc^1=D^{n-1}_\Pot\cap \mathcal{H}_\Potc^n,\qquad
\dot{D}^{n-1}_\Potc=\dot{D}^{n-1}_\Pot\cap \dot{\mathcal{H}}_\Potc^1=\dot{D}^{n-1}_\Pot\cap \dot{\mathcal{H}}_\Potc^n.
\end{equation}
All of them will be equipped with the norms inherited from ${\mathcal H}_\Pot^n$ and $\dot{\mathcal H}_\Pot^n$.
When $r=\infty$ we shall also use the product Fr\'{e}chet spaces
\begin{equation}\label{3.7}
\begin{aligned}
\mathcal{H}_\Pot^\infty & =   C^\infty(\overline{\Omega})\times C^\infty(\Gamma_1)\times C^\infty(\overline{\Omega})\times C^\infty(\Gamma_1),\\
\dot{\mathcal{H}}_\Pot^\infty & =  \dot{C}^\infty(\overline{\Omega})\times C^\infty(\Gamma_1)\times C^\infty(\overline{\Omega})\times C^\infty(\Gamma_1),
\end{aligned}
\end{equation}
and their closed subspaces
$$
\begin{alignedat}2
&D^\infty_\Pot=\{U_0\in \mathcal{H}_\Pot^\infty: \text{\eqref{3.4} hold  $\forall n\in\N$}\},\quad &&\dot{D}^\infty_\Pot=\{U_0\in \dot{\mathcal{H}}_\Pot^\infty : \text{\eqref{3.4} hold $\forall n\in\N$}\},\\
&\mathcal{H}_\Potc^\infty=\mathcal{H}_\Pot^\infty\cap \mathcal{H}^1_\Potc, \quad, \dot{\mathcal{H}}_\Potc^\infty=\dot{\mathcal{H}}_\Pot^\infty\cap \dot{\mathcal{H}}^1_\Potc, \quad &&
D_\Potc^\infty=D_\Pot^\infty\cap \mathcal{H}^1_\Potc, \quad \dot{D}_\Potc^\infty=\dot{D}_\Pot^\infty\cap \dot{\mathcal{H}}^1_\Potc.
\end{alignedat}
$$
Using Morrey's Theorem one easily gets that
\begin{equation}\label{3.9}
\begin{alignedat}4
& \mathcal{H}_\Pot^\infty =\bigcap_{n=1}^\infty  \mathcal{H}_\Pot^n,
&&\quad \dot{\mathcal{H}}_\Pot^\infty =\bigcap_{n=1}^\infty  \dot{\mathcal{H}}_\Pot^n,
&& \quad \mathcal{H}_\Potc^\infty =\bigcap_{n=1}^\infty  \mathcal{H}_\Potc^n,
&& \quad \dot{\mathcal{H}}_\Potc^\infty =\bigcap_{n=1}^\infty  \dot{\mathcal{H}}_\Potc^n,
\\
& D_\Pot^\infty =\bigcap_{n=1}^\infty  D_\Pot^n,
&&\quad \dot{D}_\Pot^\infty =\bigcap_{n=1}^\infty  \dot{D}_\Pot^n,
&& \quad D_\Potc^\infty =\bigcap_{n=1}^\infty  D_\Potc^n,
&& \quad \dot{D}_\Potc^\infty =\bigcap_{n=1}^\infty  \dot{D}_\Potc^n.
\end{alignedat}
\end{equation}
Introducing the one--dimensional subspace
$\C_{\mathcal{H}_\Pot}=\C_\Omega\times \{0\}\times \{0\}\times \{0\}$ of $\mathcal{H}_\Pot^1$, in the sequel we shall isometrically identify the quotient space $\mathcal{H}^1_\Pot/\C_{\mathcal{H}_\Pot}$ with $\dot{\mathcal{H}}^1_\Pot$. The associated  quotient (and hence surjective) map is
\begin{equation}\label{3.11}
Q\in \mathcal{L}(\mathcal{H}^1_\Pot; \dot{\mathcal{H}}^1_\Pot),\quad Q(u,v,w,z)=(\pi_0 u,v,w,z)\quad \forall (u,v,w,z)\in\mathcal{H}^1_\Pot.
\end{equation}
Since $\C_{\mathcal{H}_\Pot}\subset \mathcal{H}^n_\Potc\subset \mathcal{H}^n_\Pot$ for $n\in\widetilde{\N}$, $n\le r$ and
$\C_{\mathcal{H}_\Pot}\subset D^{m-1}_\Potc\subset D^{m-1}_\Pot$ for $m\in\widetilde{\N}$, $2\le m \le r$, and since one trivially has
\begin{equation}\label{3.12}
  Q\mathcal{H}^n_\Pot=\dot{\mathcal{H}}^n_\Pot, \quad Q\mathcal{H}^n_\Potc=\dot{\mathcal{H}}^n_\Potc,\quad
  QD^{m-1}_\Pot=\dot{D}^{m-1}_\Pot, \quad QD^{m-1}_\Potc=\dot{D}^{m-1}_\Potc,
\end{equation}
we can also make the isometrical identifications
\begin{equation}\label{3.13}
  \mathcal{H}^n_\Pot/\C_{\mathcal{H}_\Pot}\!\!=\!\!\dot{\mathcal{H}}^n_\Pot,\,\,
  \mathcal{H}^n_\Potc/\C_{\mathcal{H}_\Pot}\!\!=\!\!\dot{\mathcal{H}}^n_\Potc\, \,
  D^{m-1}_\Pot/\C_{\mathcal{H}_\Pot}\!\!=\!\!\dot{D}^{m-1}_\Pot\!\!\!,\,\,\,
  D^{m-1}_\Potc/\C_{\mathcal{H}_\Pot}\!\!=\!\!\dot{D}^{m-1}_\Potc,
\end{equation}
so getting the associated quotient (and hence surjective) maps
\begin{equation}\label{3.14}
Q\in \mathcal{L}(\mathcal{H}^n_\Pot; \dot{\mathcal{H}}^n_\Pot)\cap \mathcal{L}(\mathcal{H}^n_\Potc; \dot{\mathcal{H}}^n_\Potc)\cap\mathcal{L}(D^{m-1}_\Pot; \dot{D}^{m-1}_\Pot)\cap \mathcal{L}(D^{m-1}_\Potc; \dot{D}^{m-1}_\Potc),
\end{equation}
all of them being restrictions of $Q$ defined in \eqref{3.11}.

We also set, for $n\in\widetilde{\N}$, $n\le r$, in addition to the spaces $X^n_\Pot$ and $X^n_\Potc$ defined in \eqref{1.4} and \eqref{1.11}, the further Fr\'{e}chet spaces
\begin{equation}\label{3.15}
\begin{alignedat}2
&Y^n_\Pot= \bigcap_{i=0}^{n-1} C(\R; \mathcal{H}^{n-i}_\Pot), \qquad
&&Y^n_\Potc= \bigcap_{i=0}^{n-1} C(\R; \mathcal{H}^{n-i}_\Potc),\\
&\dot{Y}^n_\Pot= \bigcap_{i=0}^{n-1} C(\R; \dot{\mathcal{H}}^{n-i}_\Pot), \qquad
&&\dot{Y}^n_\Potc= \bigcap_{i=0}^{n-1} C(\R; \dot{\mathcal{H}}^{n-i}_\Potc).
\end{alignedat}
\end{equation}
We also remark that, by Morrey's Theorem, when $r=\infty$  one easily gets that
$$X^\infty_\Pot=C^\infty(\R\times\overline{\Omega})\times C^\infty(\R\times\Gamma_1).$$
The quotient map $Q$ in \eqref{3.11} clearly induces the "pointwise quotient" operator
\begin{equation}\label{3.16}
\mathcal{Q}\in \mathcal{L}(Y^1_\Pot; \dot{Y}^1_\Pot),\qquad (\mathcal{Q}U)t)=Q(U(t))\quad\text{for all $U\in Y^1_\Pot$ and $t\in\R$,}
\end{equation}
explicitly given, for  any $U=(u,v,w,z)\in Y^1_\Pot$ and $t\in\R$, by
\begin{equation}\label{3.17}
\mathcal{Q}U(t)=( \pi_0 u(t),v(t),w(t),z(t)).
\end{equation}
The choice of distinguishing $Q$ and $\mathcal{Q}$  is motivated by the fact that their kernels are different. Indeed, introducing the subspace $\C_{Y_\Pot}=\C_{\R\times\Omega}\times\{0\}\times\{0\}\times\{0\}$
of $Y^n_\Pot$ and $Y^n_\Potc$, for $n\in\widetilde{\N}$, $n\le \R$, by \eqref{3.17}  one easily gets that
$\ker \mathcal{Q}=C(\R; \C_{Y_\Pot})\simeq C(\R)$, while $\ker Q=\C_{\mathcal{H}_\Pot}\simeq \C$.
By \eqref{3.14} the operator $\mathcal{Q}$ trivially restricts, for $n\in\widetilde{\N}$, $n\le r$, to
\begin{equation}\label{3.19}
\mathcal{Q}\in \mathcal{L}(Y^n_\Pot; \dot{Y}^n_\Pot)\cap \mathcal{L}(Y^n_\Potc; \dot{Y}^n_\Potc).
\end{equation}
\subsection{The problem  $(\Pot_0)$ and the groups $\mathbf{\{T^n(t), t\in\R\}}$} \label{Section 3.2}

At first we make precise which types of solutions of $(\Pot)$ we shall consider in the sequel.
\begin{definition}\label{Definition 3.1}
We say that
\renewcommand{\labelenumi}{{\roman{enumi})}}
\begin{enumerate}
\item $(u,v)\in X^2_\Pot$ is a {\em strong solution} of $(\Pot)$ provided $(\Pot)_1$ holds a.e. in $\R\times\Omega$ and
$(\Pot)_2$--$(\Pot)_3$ hold a.e. on $\R\times\Gamma_1$, where $u_t$ and $\partial_{\boldsymbol{\nu}} u$ on $\R\times\Gamma_1$ are taken in the pointwise trace sense given in \S~\ref{Section 2.2};
\item $(u,v)\in X^1_\Pot$ is a {\em generalized solution} of $(\Pot)$ provided it is the limit in $X^1_\Pot$ of a sequence of strong solutions of it;
\item $(u,v)\in X^1_\Pot$ is a {\em weak solution} of $(\Pot)$ provided the distributional identities
\begin{equation}\label{3.20}
 \int_{-\infty}^\infty\left[-\rho_0\int_\Omega u_t\varphi_t+B\int_\Omega\nabla u\nabla \varphi-B\int_{\Gamma_1}v_t\varphi\right]=0,
\end{equation}
\begin{equation}\label{3.21}
 \int_{-\infty}^\infty\int_{\Gamma_1}\left[-\mu v_t\psi_t+\sigma(\nabla_\Gamma v,\nabla_\Gamma \overline{\psi})_\Gamma+\delta v_t\psi
 +\kappa v\psi-\rho_0 u\psi_t\right]=0,
 \end{equation}
hold for all $\varphi\in C^\infty_c(\R\times\R^3)$ and $\psi\in C^r_c(\R\times\Gamma_1)$.
\end{enumerate}
Moreover solutions of $(\Pot)$ ot the types i)--iii) above are said to be solutions of the same type of:
j) problem $(\Potc)$ when also \eqref{2} holds; jj)  problem $(\Pot_0)$ when also  $(\Pot_0)_4$--$(\Pot_0)_6$ hold in $X^1_\Pot$ and
jjj) problem $(\Pot_0^c)$ when both j) and jj) hold.
\end{definition}
\begin{rem}\label{Remark 3.1} Due to standard density properties on $\R\times\Gamma_1$, Definition~\ref{Definition 3.1}--iii)
is independent on the actual value of $r$, i.e. on the regularity of the boundary of $\Omega$.
\end{rem}
\begin{rem}\label{Remark 3.2}
The definition of strong (and consequently also of generalized) solution of $(\Pot)$ given above looks different from the one given in \cite[Chapter 4, \S 4.2]{mugnvit}. In the sequel we shall show that they are actually equivalent since $r\ge 2$.
Moreover the definition of weak solution of $(\Pot)$ in Definition~\ref{Definition 3.1} is actually stronger than \cite[Definition 4.2.1, \S 4.2]{mugnvit}. On the other hand, in view of Theorem~\ref{Theorem 3.1}  below, which asserts the existence of such a weak solution, the difference only concerns the more restrictive extent of the uniqueness property for it. Since this problem was exhaustively studied in the quoted paper, here we shall use Definition~\ref{Definition 3.1}--iii) for the sake of simplicity.
\end{rem}
The following result points out some trivial relations among the three types of solutions in Definition~\ref{Definition 3.1}.
\begin{lem}\label{Lemma 3.1}Let $(u,v)\in X^1_\Pot$ be a solution of $(\Pot)$ according to Definition~\ref{Definition 3.1}. Then strong $\Rightarrow$ generalized $\Rightarrow$ weak and, if $(u,v)\in X^2_\Pot$, weak $\Rightarrow$ strong.
\end{lem}
\begin{proof}Trivially strong $\Rightarrow$ generalized. Moreover, for any strong solution
$(u,v)$, integrating by parts in $\R\times\Omega$ and on $\R\times\Gamma_1$  one gets \eqref{3.20} and \eqref{3.21}, so  strong $\Rightarrow$ weak. Since the distributional identities \eqref{3.20}--\eqref{3.21} are stable with respect to the convergence in $X^1_\Pot$, by
Definition~\ref{Definition 3.1}--ii) we also get the implication generalized $\Rightarrow$ weak. To complete the proof let now $(u,v)\in X^2_\Pot$ be a weak solution. Backward performing in \eqref{3.20}--\eqref{3.21} the same operation done before to show that strong $\Rightarrow$ weak we get that equation $(\Pot)_2$ holds in the sense of distributions, and hence also a.e., in $\R\times\Omega$. Moreover we also get that
$\int_{-\infty}^\infty\int_\Gamma (\partial_{\boldsymbol{\nu}}u-\widetilde{v_t})\varphi=0
\qquad\text{for all $\varphi\in \mathcal{D}(\R\times\R^3)$,}$
where $\widetilde{v_t}$ is the trivial extension of $v_t$ to $\R\times\Gamma$. By density the last equation extends to test functions $\varphi\in C_c(\R\times\R^3)$. Since $\Gamma$ is compact and of class $C^r$ any $\varphi\in C_c(\R\times\Gamma_1)$ can be extended to $\varphi\in C_c(\R\times\R^3)$, so proving $(\Pot)_3$ a.e. on $\R\times\Gamma$. Moreover, by \eqref{3.21}, using \eqref{2.10}, we also get that $(\Pot)_2$ holds in $[C_c^r(\R\times\Gamma_1)]'$, and then a.e. on $\R\times\Gamma_1$, concluding the proof.
\end{proof}
\begin{rem}\label{Remark 3.3} A trivial consequence of Lemma~\ref{Lemma 3.1} is that solutions $(u,v)\in X^2_\Pot$ are equivalently strong, generalized or weak, and thus, for every such solution, we have $(u(t), v(t), u_t(t), v_t(t))\in D^1_\Pot$ for all $t\in\R$.
\end{rem}
Problem $(\Pot_0)$ was studied in \cite[Chapters 4 and 5]{mugnvit} by a semigroup approach. To recall it we introduce the unbounded operator $A_\Pot: D(A_\Pot)\subset \mathcal{H}^1_\Pot \to \mathcal{H}^1_\Pot$ given by
\begin{equation}\label{3.23}
 D(A_\Pot)=D^1_\Pot,\qquad A_\Pot\begin{pmatrix}u\\v\\w\\z\end{pmatrix} =
\begin{pmatrix}-w\\-z\\-(B/\rho_0)\Delta u\\
\frac 1\mu\left[-\DivGamma(\sigma\nabla_\Gamma v)+\delta z+\kappa v+\rho_0 w_{|\Gamma_1}\right]
\end{pmatrix},
\end{equation}
together with the abstract equation and Cauchy problem
\begin{gather}\label{3.24}U_\Pot'+A_\Pot U_\Pot=0\qquad\text{in $\mathcal{H}^1_\Pot$,}
\\\label{3.25}
U_\Pot'+A_\Pot U_\Pot=0\qquad\text{in $\mathcal{H}^1_\Pot$}, \qquad U_\Pot(0)=U_{0\Pot}\in \mathcal{H}^1_\Pot.
\end{gather}
We point out that $A_\Pot$ in \eqref{3.23} coincides with the operator $A$ defined in \cite[Chapter 4, (4.3)--(4.4)]{mugnvit} although it domain $D(A_\Pot)$ appears to be strictly contained in $D(A)$ and the operators $\Delta$, $\DivGamma(\sigma \nabla_\Gamma)$, $\partial_{\boldsymbol{\nu}}$ in \cite[Chapter 4, (4.3)--(4.4)]{mugnvit} are taken in a distributional sense. Indeed, while in \cite[Chapter 4]{mugnvit} also the case $r=1$ was considered, in the present paper we have $r\ge 2$. Hence we can use the characterization of $D(A)$ given in \cite[Chapter 5, Lemma 5.0.4]{mugnvit}, which gives $D(A)=D^1_\Pot$ as in \eqref{3.23} and (see \cite[Chapter 3, \S 3.3.2 and \S 3.35]{mugnvit}) the operators $\Delta$, $\DivGamma(\sigma \nabla_\Gamma)$ can be taken in the a.e. sense and $\partial_{\boldsymbol{\nu}}$ in the trace one, as done in \eqref{3.23}.

The relation between problems $(\Pot)$, $(\Pot_0)$ and their abstract versions \eqref{3.24}--\eqref{3.25} is given by the following result.
\begin{lem}\label{Lemma 3.2}
The couple $(u,v)\in X^2_\Pot$ is a strong solution of $(\Pot)$ if and only if $u$ and $v$ are the first two components of a strong solution $U_\Pot=(u,v,w,z)$ of \eqref{3.24}. Moreover in this case $U_\Pot\in Y^2_\Pot$.
The same relation occurs between generalized solutions $(u,v)\in X^1_\Pot$  of $(\Pot)$ and generalized solutions $U_\Pot\in Y^1_\Pot$ of \eqref{3.24}.
Finally, in previous statements, one can replace $(\Pot)$ with  $(\Pot_0)$ provided one also replaces \eqref{3.24} with \eqref{3.25}.
\end{lem}
\begin{proof}
For strong solutions the asserted relation essentially follows from  by \cite[Chapter 5, Lemma 5.0.4]{mugnvit}. Its generalization to strong solutions follows by a standard density argument, as in \cite[\S~4.2]{mugnvit}.
\end{proof}
Since generalized and strong solutions of $(\Pot)$ where defined in \cite[Chapter 4]{mugnvit} as couples $(u,v)$ constituted by the first two components of an homologous solution of \eqref{3.24}, Lemma~\ref{Lemma 3.2} also shows that in Definition~\ref{Definition 3.1} we just made more explicit the definition in the quoted paper.

We now recall the combination of \cite[Theorem 1.2.1 and Theorem 4.1.5]{mugnvit} in a form which is most adequate for our purpose.
Trivially Theorem~\ref{Theorem 1.1} is just a simplified form of the following statement.
\begin{thm}[\bf Well--posedness for $(\Pot_0)$ and \eqref{3.25}]\label{Theorem 3.1}
The operator $-A_\Pot$ is densely defined and it generates on $\mathcal{H}^1_\Pot$ a strongly continuous group $\{T^1(t), t\in\R\}$.
Consequently, for any $U_{0\Pot}=(u_0,v_0,u_1, v_1)\in \mathcal{H}^1_\Pot$, problem \eqref{3.25} has a unique generalized solution $U_\Pot\in Y^1_\Pot$ defined by $U_\Pot(t)=T^1(t)[U_{0\Pot}]$ for all $t\in\R$ and problem $(\Pot_0)$ has a unique generalized solution $(u,v)\in X^1_\Pot$. Moreover $U_\Pot=(u,v,u_t,v_t)$.
 Next, the solutions $U_\Pot$ and $(u,v)$ continuously depend on $U_{0\Pot}$ in the topology of the respective spaces, $(u,v)$ satisfies the energy identity \eqref{1.5} and is unique also among weak solutions of  $(\Pot_0)$.
Finally, for any $U_{0\Pot}\in \mathcal{H}^1_\Pot$, the following properties are equivalent:
i) $U_{0\Pot}\in D^1_\Pot$; ii) $U_\Pot$ is a strong solution; \quad iii) $(u,v)$ is a strong solution; \quad iv) $U_\Pot\in Y^2_\Pot$;
\quad v) $(u,v)\in X^2_\Pot$.
\end{thm}
\begin{proof} The first two sentences follow from \cite[Theorem 4.1.5]{mugnvit} and Lemma~\ref{Lemma 3.2}, the asserted continuous dependence follows by Lemma~\ref{Lemma 2.1}--i). The uniqueness of $(u,v)$ among weak solutions of $(\Pot_0)$ follows, since our definition of weak solutions is more restrictive that the one given in the quoted paper, follows by \cite[Theorem 1.2.1 or Lemma 4.2.5]{mugnvit}, while the energy identity is asserted in \cite[Theorem 1.2.1]{mugnvit}. To prove the final asserted equivalence we remark that by Lemma~\ref{Lemma 2.1}--iii) we get i) $\Leftrightarrow$ ii), while ii) $\Leftrightarrow$ iii) $\Leftrightarrow$ iv) follow from Lemma~\ref{Lemma 3.2}. Finally iv) $\Leftrightarrow$ v) is trivial.
\end{proof}

Theorem~\ref{Theorem 3.1} also shows that solutions of $(\Pot)$ are {\em equivalently} weak or generalized. Combining this remark with Remark~\ref{Remark 3.3} we thus recognize that {\em all types of solutions in Definition~\ref{Definition 3.1} coincide, strong solutions being defined only in the class $X^2_\Pot$}. Consequently in the sequel we  shall only deal with weak solutions of $(\Pot)$, i.e. with elements of the spaces $\mathcal{S}^n_\Pot$ defined in \eqref{1.18} for  $n\in\widetilde{\N}$, $n\le r$. For the same $n$'s we also set  the Fr\'echet spaces
\begin{equation}\label{3.26}
\mathcal{T}^n_\Pot=\{U_\Pot\in Y^n_\Pot: \quad \text{$U_\Pot$ is a generalized solution of \eqref{3.24}}\},
\end{equation}
endowed with the topology inherited from $Y^n_\Pot$. The connection between them and the spaces $\mathcal{S}^n_\Pot$ is given by the following result, which  trivially follows from Lemma~\ref{Lemma 3.2}, \eqref{1.4} and \eqref{3.15}.
\begin{lem}\label{Lemma 3.3}
The operator $\mathcal{I}\in\mathcal{L}(X^1_\Pot;Y^1_\Pot)$, defined by
$\mathcal{I}(u,v)=(u,v,u_t,v_t)$,
restricts for each $n\in\widetilde{\N}$, $n\le r$, to a bijective isomorphism $\mathcal{I}^n_\Pot\in\mathcal{L}(\mathcal{S}^n_\Pot;\mathcal{T}^n_\Pot)$ with inverse $\left(\mathcal{I}^n_\Pot\right)^{-1}\in\mathcal{L}(\mathcal{T}^n_\Pot;\mathcal{S}^n_\Pot)$ simply given by
$\left(\mathcal{I}^n_\Pot\right)^{-1}(u,v,w,z)=(u,v)$.
\end{lem}
To show that the regularity classes $\mathcal{S}^n_\Pot$ and $\mathcal{T}^n_\Pot$ introduced above contain physically significant solutions (also when $n\ge 3$) we now recall the regularity results \cite[Theorem~1.2.3 and Lemma 5.0.4]{mugnvit} combining them as follows.
\begin{thm}\label{Theorem 3.2}For $n\in\widetilde{\N}$, $2\le n\le r$, we have $D(A^{n-1}_\Pot)=D^{n-1}_\Pot$, the respective norms being equivalent. Consequently the operator $-\,\, {_n}A_\Pot$ given by \eqref{2.18} generates on $D^{n-1}_\Pot$ a strongly continuous group $\{T^n(t), t\in\R\}$ given by $T^n(t)=T^1(t)_{|D^{n-1}_\Pot}$ for all $t\in\R$.
Finally, for any $U_{0\Pot}\in\mathcal{H}^1_\Pot$, denoting by $(u,v)\in X^1_\Pot$ the weak solution given in Theorem~\ref{Theorem 3.1} and $U_\Pot=(u,v,u_t,v_t)$, for any $n\in\widetilde{\N}$, $2\le n\le r$, we have the equivalences
\begin{equation}\label{3.29}
  U_{0\Pot}\in D^{n-1}_\Pot \Leftrightarrow U_\Pot\in Y^n_\Pot \Leftrightarrow (u,v)\in X^n_\Pot.
\end{equation}
Moreover, when $U_{0\Pot}\in D^{n-1}_\Pot$, $U_\Pot$ and $(u,v)$ continuously depend on it in the topologies of the respective spaces.
\end{thm}
\subsection{The quotient groups $\mathbf{\{\mathcal{P}^n(t), t\in\R\}}$} \label{Section 3.3} The group $\{T^1(t),t\in\R\}$ introduced in Theorem~\ref{Theorem 3.1} and the induced groups  $\{T^n(t),t\in\R\}$ introduced in Theorem~\ref{Theorem 3.2} are not directly connected to the problems $(\Eul)$, $(\Eulc)$ and $(\Lagr)$. Hence we are now going to introduce their quotient groups with respect to the space $\C_{\mathcal{H}_\Pot}$ defined in \S~\ref{Section 3}. Indeed, since $(\Pot)$ possesses the trivial solutions $u(t,x)\equiv c\in\C$, $v(t,x)\equiv 0$, the space $\C_{\mathcal{H}_\Pot}$is invariant under the flow of all the groups $\{T^n(t),t\in\R\}$. The identifications \eqref{3.13} allow to make these quotient groups more concrete.

We preliminarily set the unbounded operator  $\dot{A}_\Pot: D(\dot{A}_\Pot)\subset \dot{\mathcal{H}}^1_\Pot \to \dot{\mathcal{H}}^1_\Pot$ defined by
\begin{gather}\label{3.30d}
 D(\dot{A}_\Pot)=\dot{D}^1_\Pot, \\
 \label{3.30}
 \dot{A}_\Pot\begin{pmatrix}\pi_0 u\\v\\w\\z\end{pmatrix} =
Q A_\Pot\begin{pmatrix}u\\v\\w\\z\end{pmatrix} =
\begin{pmatrix}-\pi_0 w\\-z\\-(B/\rho_0)\Delta u\\
\frac 1\mu\left[-\DivGamma(\sigma\nabla_\Gamma v)+\delta z+\kappa v+\rho_0 w_{|\Gamma_1}\right]
\end{pmatrix},
\end{gather}
for all $(u,v,w,z)\in D^1_\Pot$, together with the abstract equation and Cauchy problem
\begin{gather}\label{3.31}\dot{U}_\Pot'+\dot{A}_\Pot \dot{U}_\Pot=0\qquad\text{in $\dot{\mathcal{H}}^1_\Pot$,}
\\\label{3.32}
\dot{U}_\Pot'+\dot{A}_\Pot \dot{U}_\Pot=0\qquad\text{in $\dot{\mathcal{H}}^1_\Pot$}, \qquad \dot{U}_\Pot(0)=\dot{U}_{0\Pot}\in \dot{\mathcal{H}}^1_\Pot.
\end{gather}
The following result is a direct consequence of Theorems~\ref{Theorem 3.1}--\ref{Theorem 3.2}.
\begin{prop}[\bf Well--posedness and regularity for \eqref{3.32}]\label {Proposition 3.1}\phantom{A}

\renewcommand{\labelenumi}{{\Roman{enumi})}}
\begin{enumerate}
\item The operator $-\dot{A}_\Pot$ is densely defined and it generates on $\dot{\mathcal{H}}^1_\Pot=\mathcal{H}^1_\Pot/\C_{\mathcal{H}_\Pot}$ the quotient strongly continous group $\{\Pot^1(t), t\in\R\}$ defined by
    \begin{equation}\label{3.33}
      \Pot^1(t)[QU_{0\Pot}]=QT^1(t)[U_{0\Pot}]\qquad\text{for all $U_{0\Pot}\in \mathcal{H}^1_\Pot$ and $t\in\R$.}
    \end{equation}
    Consequently, for any $\dot{U}_{0\Pot}\in \dot{\mathcal{H}}^1_\Pot$, problem \eqref{3.32} has a unique generalized solution $\dot{U}_\Pot\in \dot{Y}_\Pot^1$ given by $\dot{U}_\Pot(t)=\Pot^1(t)[\dot{U}_{0\Pot}]$, $t\in\R$, continously depending on $\dot{U}_{0\Pot}$ in the topologies of the respective spaces.

    More explicitly, denoting by $U_\Pot\in Y^1_\Pot$ the unique generalized solution of \eqref{3.25} given by Theorem~\ref{Theorem 3.1} with any  initial datum $U_{0\Pot}\in Q^{-1}(\dot{U}_{0\Pot})$ and by $(u,v)$ the corresponding weak solution of $(\Pot_0)$, we have
    \begin{equation}\label{3.34}
     \dot{U}_\Pot(t)= Q U_\Pot(t)= (\pi_0 u(t), v(t), u_t(t), v_t(t))\qquad\text{for all $t\in\R$.}
    \end{equation}
    Moreover $\dot{U}_{0\Pot}\in \dot{D}^1_\Pot\Leftrightarrow$ $U_\Pot$ is a strong solution $\Leftrightarrow \dot{U}_\Pot\in \dot{Y}^2_\Pot$.
\item For $n\in\widetilde{\N}$, $2\le n\le r$, we have $D(\dot{A}^{n-1}_\Pot)=\dot{D}^{n-1}_\Pot$, the respective norms being equivalent. Consequently the operator $-\,\, {_n}\dot{A}_\Pot$ given by \eqref{2.18} generates on $\dot{D}^{n-1}_\Pot$ the strongly continuous group $\{\mathcal{P}^n(t), t\in\R\}$ given by $\Pot^n(t)=\Pot^1(t)_{|\dot{D}^{n-1}_\Pot}$ for all $t\in\R$.
    The group $\{\mathcal{P}^n(t), t\in\R\}$ can be equivalently defined as the quotient group on $\dot{D}^{n-1}_\Pot=D^{n-1}_\Pot$ given by
    \begin{equation}\label{3.35}
      \Pot^n(t)[QU_{0\Pot}]=QT^n(t)[U_{0\Pot}]\qquad\text{for all $U_{0\Pot}\in D^{n-1}_\Pot$ and $t\in\R$.}
    \end{equation}
\item For any $\dot{U}_{0\Pot}\in\mathcal{H}^1_\Pot$ and $n\in\widetilde{\N}$, $2\le n\le r$, we have  $\dot{U}_{0\Pot}\in \dot{D}^{n-1}_\Pot$ if and only if  $\dot{U}_\Pot\in \dot{Y}^n_\Pot$, and in this case $\dot{U}_\Pot$  continuously depends on $\dot{U}_{0\Pot}$ in the topologies of the respective spaces.
\end{enumerate}
\end{prop}
\begin{proof}To prove part I) we first remark that, since by Theorem~\ref{Theorem 3.1} $D(A_\Pot)=D^1_\Pot$ is dense in $\mathcal{H}^1_\Pot$ one immediately gets that $D(\dot{A}_\Pot)=\dot{D}^1_\Pot=D^1_\Pot/\C_{\mathcal{H}_\Pot}$ is dense in $\dot{\mathcal{H}}^1_\Pot=\mathcal{H}^1_\Pot/\C_{\mathcal{H}_\Pot}$, where we used \eqref{3.13}. Moreover a standard construction in semigroup theory shows, using the already remarked invariance of $\C_{\mathcal{H}_\Pot}$, that \eqref{3.33} defines on $\mathcal{H}^1_\Pot$ a (quotient) strongly continuous group. Its generator, using Theorem ~\ref{Theorem 3.1}, is given by the operator $-QA_\Pot$ with domain $QD(A_\Pot)$.
Hence, since $D(A_\Pot)=D^1_\Pot$, by \eqref{3.13} we have $QD(A_\Pot)=\dot{D}^1_\Pot=D(\dot{A}_\Pot)$. Comparing \eqref{3.23} and \eqref{3.30} one easily gets $QA_\Pot=\dot{A}_\Pot$. We then get the well--posedness of problem \eqref{3.32} asserted in the statement and, by Theorem~\ref{Theorem 3.1} and \eqref{3.11}, also \eqref{3.34}. By Lemma~\ref{Lemma 2.1}--iii) then $\dot{U}_{0\Pot}\in\dot{D}^1_\Pot$ if and only if $\dot{U}_\Pot$ is a strong solution of \eqref{3.32}. To complete the proof of part I) we first claim that the graph norm and the $\dot{\mathcal{H}}^2_\Pot$ -- norm are equivalent on $\dot{D}^1_\Pot$. By \eqref{3.30} one sees that $\dot{A}_\Pot\in\mathcal{L}(\dot{D}^1_\Pot; \dot{\mathcal{H}}^1_\Pot)$, so $\|\cdot\|_{D(\dot{A}_\Pot)}\le \const \|\cdot\|_{\dot{\mathcal{H}}^2_\Pot}$ on $\dot{D}^1_\Pot$ and, since $\dot{D}^1_\Pot$ is complete with respect to both norms, our claim follows by the Two Norms Theorem. Hence, if $\dot{U}_\Pot\in \dot{Y}^2_\Pot$, by \eqref{3.15} we have $\dot{U}_\Pot(0)\in \dot{H}^1_\Pot$, which by the definition of domain  yields $\dot{U}_{0\Pot}\in D(\dot{A}_\Pot)=\dot{D}^1_\Pot$. Conversely, if $\dot{U}_{0\Pot}\in D(\dot{A}_\Pot)=\dot{D}^1_\Pot$, by Lemma~\ref{Lemma 2.1}--iii) we get $\dot{U}_\Pot\in C(\R; D(\dot{A}_\Pot))\cap C^1(\R: \dot{\mathcal{H}}_\Pot)$, which by previous claim implies $\dot{U}_\Pot\in \dot{Y}^2_\Pot$, completing the proof of part I).

To prove part II) we first claim that $D(\dot{A}^{n-1}_\Pot)=\dot{D}^{n-1}_\Pot$ for $n\in\widetilde{\N}$, $2\le n\le r$. Since $D(\dot{A}_\Pot)=\dot{D}^1_\Pot$ by induction we suppose that $n\ge 3$ and $D(\dot{A}^{n-2}_\Pot)=\dot{D}^{n-2}_\Pot$. Using \eqref{2.16}, \eqref{3.13} and \eqref{3.30} we have
$D(\dot{A}^{n-1}_\Pot)=\{QU_{0\Pot}: U_{0\Pot}\in D^{n-2}_\Pot\quad\text{and}\quad QA_\Pot U_{0\Pot}\in \dot{D}^{n-2}_\Pot\}$.
Since $\C_{\mathcal{H}_\Pot}\subset D^{n-2}_\Pot$ in the previous formula we have $QA_\Pot U_{0\Pot}\in \dot{D}^{n-2}_\Pot$ if and only if $A_\Pot U_{0\Pot}\in D^{n-2}_\Pot$, so $D(\dot{A}^{n-1}_\Pot)=Q D^{n-1}_\Pot=\dot{D}^{n-1}_\Pot$, proving our claim. The asserted norm equivalence is then proved by the argument already used to get it in the case $n=2$ above once we recognize that
$\dot{A}_\Pot\in\mathcal{L}(\dot{D}^i_\Pot;\dot{D}^{i-1}_\Pot)$ for $i=1,\ldots, n-1$, where $\dot{D}^0_\Pot:=\dot{H}^1_\Pot$. But, by Theorem~\ref{Theorem 3.2}, $A_\Pot\in \mathcal{L}(D^i_\Pot;D^{i-1}_\Pot)$. Then, by \eqref{3.14}, $QA_\Pot\in\mathcal{L}(D^i_\Pot;\dot{D}^{i-1}_\Pot)$ and consequently, since $\dot{D}^i_\Pot=D^i_\Pot/\C_{\mathcal{H}_\Pot}$, we get
$\dot{A}_\Pot\in\mathcal{L}(\dot{D}^i_\Pot;\dot{D}^{i-1}_\Pot)$.
By applying Lemma~\ref{Lemma 2.1}--ii) we then get the strongly continuous group $\{\Pot^n(t), t\in\R\}$ on $\dot{D}^{n-1}_\Pot$, given by  $\Pot^n(t)=\Pot^1(t)_{|\dot{D}^{n-1}_\Pot}$ for all $t\in\R$, which has generator $- {_n}\dot{A}_\Pot$. Since $T^n(t)=T^1(t)_{|D^{n-1}_\Pot}$, \eqref{3.35} follows by \eqref{3.33}.

To prove part III) we first consider $n\in\N$ and we remark that when $\dot{U}_{0\Pot}\in \dot{D}^{n-1}_\Pot$, by Lemma~\ref{Lemma 2.1} we get
$\dot{U}_\Pot\in \dot{Y}^n_\Pot$ while, when $\dot{U}_\Pot\in \dot{Y}^n_\Pot$ we have $\dot{U}_\Pot'(0)\in \dot{H}^{n-1}_\Pot$,  so
$\dot{U}_{0\Pot}\in D({_{n-2}}\dot{A}_\Pot)=\dot{D}^{n-1}_\Pot)$. The case $n=\infty$ then follows by \eqref{3.9}.
\end{proof}
\begin{rem}\label{Remark 3.4} While the first two components $u$ and $v$ of a solution of \eqref{3.25} constitute a solution of $(\Pot)$, the first two components $\dot u$ and $v$ of a solution of \eqref{3.31} do not enjoy the same property, in general. We recall the example  given in \cite[\S 6.1]{mugnvit}. Let  $\kappa(x)\equiv k_0\in \R\setminus\{0\}$ and  the solutions $u(t,x)=u_1 t$, $v(t,x)\equiv -\rho_0u_1/k_0$. By \eqref{3.34} the strong solution $\dot{U}_\Pot$ of \eqref{3.32} with data $\dot{U}_{0\Pot}=(0,-\rho_0 u_1/k_0, u_1,0)$ is $\dot{U}_\Pot(t)=\dot{U}_{0\Pot}$ for all $t\in\R$. Trivially
$(0,-\rho_0 u_1/k_0)$ does not solve equation $(\Pot)_2$.
\end{rem}
We now set , for $n\in\widetilde{\N}$, $n\le r$, the Fr\'{e}chet space constituted by  the trajectories of the group $\{\Pot^n(t),t\in\R\}$, i.e.
\begin{equation}\label{3.33BIS}
\dot{\mathcal{T}}^n_\Pot=\{\dot{U}_\Pot\in \dot{Y}^n_\Pot: \dot{U}_\Pot\quad\text{is a generalized solution of \eqref{3.31}}\},
\end{equation}
endowed with the topology inherited from $\dot{Y}^n_\Pot$. The following result relates solutions of $(\Pot)$, \eqref{3.24} and \eqref{3.31} and it will be relevant in connection with problem $(\Eul)$.
\begin{prop}\label{Proposition 3.2}
The pointwise quotient operator $\mathcal{Q}$ defined in \eqref{3.16} restricts to a surjective operator $\mathcal{Q}^1_\Pot\in\mathcal{L}(\mathcal{T}^1_\Pot;\dot{\mathcal{T}}^1_\Pot)$, so the operator $\mathcal{J}^1_\Pot=\mathcal{Q}^1_\Pot\cdot \mathcal{I}^1_\Pot\in\mathcal{L}(\mathcal{S}^1_\Pot; \dot{\mathcal{T}}^1_\Pot)$
given by
\begin{equation}\label{3.34BIS}
  [\mathcal{J}^1_\Pot(u,v)](t)=(\pi_0 u(t), v(t), u_t(t), v_t(t))\qquad\text{for all $t\in\R$,}
\end{equation}
is surjective as well.
Moreover, for all $n\in \widetilde{\N}$, $n\le r$, the operators $\mathcal{Q}^1_\Pot$,  $\mathcal{J}^1_\Pot$
restrict to surjective operators $\mathcal{Q}^n_\Pot\in\mathcal{L}(\mathcal{T}^n_\Pot;\dot{\mathcal{T}}^n_\Pot)$,
$\mathcal{J}^n_\Pot\in\mathcal{L}(\mathcal{S}^n_\Pot; \dot{\mathcal{T}}^n_\Pot)$ and we have
\begin{equation}\label{3.35BIS}
  \ker\mathcal{Q}^n_\Pot=\C_{Y_\Pot}\qquad\text{and}\quad   \ker\mathcal{J}^n_\Pot=\C_{X_\Pot}.
\end{equation}
 Consequently $\mathcal{J}^n_\Pot$ subordinates a bijective isomorphism $\dot{\mathcal{J}}^n_\Pot\in\mathcal{L}(\dot{\mathcal{S}}^n_\Pot; \dot{\mathcal{T}}^n_\Pot)$.
\end{prop}
\begin{proof} Combining \eqref{3.16} with Proposition~\ref{Proposition 3.1} we get the surjective operator $\mathcal{Q}^1_\Pot\in\mathcal{L}(\mathcal{T}^1_\Pot;\dot{\mathcal{T}}^1_\Pot)$, so
by Lemma~\ref{Lemma 3.3} the operator $\mathcal{J}^1_\Pot=\mathcal{Q}^1_\Pot\cdot \mathcal{I}^1_\Pot\in\mathcal{L}(\mathcal{S}^1_\Pot; \dot{\mathcal{T}}^1_\Pot)$ is surjective and it satisfies \eqref{3.34BIS}. By \eqref{3.19} we get that $\mathcal{Q}^1_\Pot$ restricts to $\mathcal{Q}^n_\Pot\in\mathcal{L}(\mathcal{T}^n_\Pot;\dot{\mathcal{T}}^n_\Pot)$, so by Lemma~\ref{Lemma 3.3} the operator $\mathcal{J}^1_\Pot$ restricts to $\mathcal{J}^n_\Pot\in\mathcal{L}(\mathcal{S}^n_\Pot; \dot{\mathcal{T}}^n_\Pot)$.
To show that $\mathcal{Q}^n_\Pot$ is surjective also when $n\ge 2$ let us take $\dot{U}_\Pot\in\dot{\mathcal{T}}^n_\Pot$. By Proposition~\ref{Proposition 3.1}--III)  we then have $\dot{U}_\Pot(0)\in \dot{D}^{n-1}_\Pot$. By part I) of the same result
we have $\dot{U}_\Pot=\mathcal{Q}U_\Pot$, where $U_\Pot$ is the solution of the Cauchy problem \eqref{3.25} corresponding to any datum $U_{0\Pot}\in Q^{-1}(\dot{U}_\Pot(0))$ which by \eqref{3.12} is equivalent to require that $U_{0\Pot}\in D^{n-1}_\Pot$. Consequently, by Theorem~\ref{Theorem 3.2}, formula \eqref{3.29}, we have $U_\Pot\in Y^n_\Pot$ and then $\dot{U}_\Pot=\mathcal{Q}^n_\Pot U_\Pot$, proving that $\mathcal{Q}^n_\Pot$ is surjective. The surjectivity of $\mathcal{J}^n_\Pot$ then follows by Lemma~\ref{Lemma 3.3}.
To prove \eqref{3.35BIS} we first remark that, by Lemma~\ref{Lemma 3.3}, we just have to prove that $\ker\mathcal{J}^n_\Pot=\C_{X_\Pot}$. Moreover, since $\C_{X_\Pot}\subset\mathcal{S}^n_\Pot$, we just have to prove that
$\ker\mathcal{J}^1_\Pot=\C_{X_\Pot}$, this fact being  a trivial consequence of \eqref{3.34BIS}.
\end{proof}

Proposition~\ref{Proposition 3.2} shows that solutions of $(\Pot)$ correspond, up to an additive space--time constant, to trajectories of $\{T^1(t), t\in\R\}$, regularity classes being preserved.

\subsection{The problem  $(\Pot_0^c)$ and the groups $\mathbf{\{T_c^n(t), t\in\R\}}$, $\mathbf{\{\mathcal{P}_c^n(t), t\in\R\}}$ } \label{Section 3.4} In this subsection we are going to show that the results in \S~\ref{Section 3.2} -- \ref{Section 3.3} apply as well to the restricted potential problem $(\Potc)$. This fact depends on the following easy result, which (trivially) generalizes \cite[Lemma 6.2.1]{mugnvit}, with the same proof.
\begin{lem}\label{Lemma 3.4} For any $(u,v)\in\mathcal{S}^1_\Pot$ we have
\begin{equation}\label{3.36}
\rho_0\int_\Omega u_t(t)-B\int_{\Gamma_1}v(t)=\rho_0\int_\Omega u_t(s)-B\int_{\Gamma_1}v(s)\qquad\text{for all $s,t\in\R$.}
\end{equation}
Consequently the spaces $\mathcal{H}^1_\Potc$ and $\dot{\mathcal{H}}^1_\Potc$ are invariant with respect to the flows of the groups
$\{T^1(t), t\in\R\}$ and $\{\Pot^1(t), t\in\R\}$, as, for $n\in\N$, $2\le n\le r$, are the spaces $D^{n-1}_\Potc$ and $\dot{D}^{n-1}_\Potc$ with respect to the flows of the groups $\{T^n(t), t\in\R\}$ and $\{\Pot^n(t), t\in\R\}$. Moreover
\begin{equation}\label{3.37}
\Rg(A_\Pot)\subseteq \mathcal{H}^1_\Potc\qquad\text{and}\quad  \Rg(\dot{A}_\Pot)\subseteq \dot{\mathcal{H}}^1_\Potc.
\end{equation}
\end{lem}
To state the extensions of the results in \S~\ref{Section 3.2} -- \ref{Section 3.3} we now set, thanks to \eqref{3.37},
the unbounded operators $A_\Potc: D(A_\Potc)\subset \mathcal{H}^1_\Potc \to \mathcal{H}^1_\Potc$  and $\dot{A}_\Potc: D(\dot{A}_\Potc)\subset \dot{\mathcal{H}}^1_\Potc \to \dot{\mathcal{H}}^1_\Potc$, defined by
\begin{equation}\label{3.38}
D(A_\Potc)=D^1_\Potc, \quad A_\Potc=A_{\Pot|D^1_\Potc}, \qquad D(\dot{A}_\Potc)=\dot{D}^1_\Potc, \quad \dot{A}_\Potc=\dot{A}_{\Pot|\dot{D}^1_\Potc}.
\end{equation}
Moreover we set, for $n\in\widetilde{\N}$, $n\le r$, the Fr\'{e}chet spaces
\begin{equation}\label{3.39}
\mathcal{S}^n_\Potc=\mathcal{S}^n_\Pot\cap X^n_\Potc,\qquad \mathcal{T}^n_\Potc=\mathcal{T}^n_\Pot\cap Y^n_\Potc \qquad\text{and}\quad
\dot{\mathcal{T}}^n_\Potc=\dot{\mathcal{T}}^n_\Pot\cap \dot{Y}^n_\Potc,
\end{equation}
and we point out the following extension of Lemma~\ref{Lemma 3.3}.
\begin{cor}\label{Corollary 3.1}
 For each $n\in\widetilde{\N}$, $n\le r$, the operator  $\mathcal{I}^n_\Pot$ in Lemma~\ref{Lemma 3.3} further restricts
 to a bijective isomorphism $\mathcal{I}^n_\Potc\in\mathcal{L}(\mathcal{S}^n_\Potc;\mathcal{T}^n_\Potc)$, with inverse
  given as in Lemma~\ref{Lemma 3.3}.
\end{cor}
\begin{proof}
By \eqref{1.6}, \eqref{3.15} and \eqref{3.39} one trivially gets $\mathcal{I}^1_\Pot\cdot\mathcal{S}^1_\Potc=\mathcal{T}^1_\Potc$.
When $n>1$ the statement follows in the same way once one remarks that, differentiating in time, we have
$X^n_\Potc=\left\{(u,v)\in X^n_\Pot: \rho_0\int_\Omega \frac {\partial^i u}{\partial t^i}=B\int_{\Gamma_1}\frac {\partial^{i-1} v}{\partial t^i}\quad\text{for $i\in\N$, $i\le n$}\right\}$,
from which one easily gets that $\mathcal{I}^n_\Pot\cdot\mathcal{S}^n_\Potc=\mathcal{T}^n_\Potc$.
\end{proof}
Using Lemma~\ref{Lemma 3.4} and Corollary~\ref{Corollary 3.1} it is then straightforward to get the following consequences of
the results in \S~\ref{Section 3.2} -- \ref{Section 3.3} which will be used in sequel. Corollary~\ref{Corollary 1.1} is just a partial statement of it.
\begin{cor}[\bf Problem $\mathbf{(\Pot_0^c)}$ and the groups $\mathbf{\{T^n_c(t), t\in\R\}}$]\label{Corollary 3.2}The
statements of Theorem~\ref{Theorem 3.1} and \ref{Theorem 3.2} continue to hold when one respectively replaces (also in problem \eqref{3.25}):
 the spaces $\mathcal{H}^1_\Pot$, $D^{n-1}_\Pot$, $X^n_\Pot$, $Y^n_\Pot$ with the spaces $\mathcal{H}^1_\Potc$, $D^{n-1}_\Potc$, $X^n_\Potc$,  $Y^n_\Potc$;
the operators $A_\Pot$ and ${_n}A_\Pot$ with the operators $A_\Potc$ and ${_n}A_\Potc$;
problem $(\Pot_0)$ with problem $(\Pot_0^c)$;
the groups $\{T^n(t),t\in\R\}$ with the groups  $\{T^n_c(t),t\in\R\}$ defined (denoting $D^0_\Potc=\mathcal{H}^1_\Potc$) by
$T^n_c(t)=T^n(t)_{|D^{n-1}_\Potc}$ for $n\in\N$, $n\le r$.
\end{cor}
\begin{proof} Recalling \eqref{3.2} and \eqref{3.6}, the statement follows by  a standard construction in semigroup theory combined with  Theorem~\ref{Theorem 3.1},  \ref{Theorem 3.2} and  Corollary~\ref{Corollary 3.1}.
\end{proof}
Using Lemma~\ref{Lemma 3.4} one can also generalize Proposition~\ref{Proposition 3.1} as follows.
\begin{cor}[\bf Problem \eqref{3.32} in $\mathbf{\dot{\mathcal{H}}^1_\Potc}$ and the groups $\mathbf{\{\mathcal{P}_c^n(t), t\in\R\}}$]
\label{Corollary 3.3}
The statement of Proposition~\ref{Proposition 3.1} continues to hold when one respectively replaces, also in problem \eqref{3.32}:
the spaces $\dot{\mathcal{H}}^1_\Pot$, $\dot{D}^{n-1}_\Pot$, $\dot{Y}^n_\Pot$ with the spaces $\dot{\mathcal{H}}^1_\Potc$, $\dot{D}^{n-1}_\Potc$, $\dot{Y}^n_\Potc$; the operators $\dot{A}_\Pot$ and ${_n}\dot{A}_\Pot$ with the operators $\dot{A}_\Potc$ and ${_n}\dot{A}_\Potc$; the group $\{\Pot^n(t),t\in\R\}$ with the subspace group  $\{\Pot^n_c(t),t\in\R\}$ defined (denoting $\dot{D}^0_\Potc=\dot{\mathcal{H}}^1_\Potc$) by $\Pot^n_c(t)=\Pot^n(t)_{|\dot{D}^{n-1}_\Potc}$ for $n\in\N$, $n\le r$.
\end{cor}
\begin{proof}
We first claim that
\begin{equation}\label{3.40}
\dot{Y}^n_\Potc =\dot{Y}^1_\Potc\cap \dot{Y}^n_\Pot\qquad\text{for all $n\in\widetilde{\N}$, $n\le r$.}
\end{equation}
By \eqref{3.15} trivially $\dot{Y}^n_\Potc \subseteq\dot{Y}^1_\Potc\cap \dot{Y}^n_\Pot$ and, deriving in time as many times as  needed, one also gets the reverse inclusion. By Proposition~\ref{Proposition 3.1}, using Lemma~\ref{Lemma 3.4} and \eqref{3.40}, one gets the assertion.
\end{proof}
Also Proposition~\ref{Proposition 3.2} extends to problem $(\Potc)$ as follows.
\begin{cor}\label{Corollary 3.4}For all $n\in \widetilde{\N}$, $n\le r$, the operators $\mathcal{Q}^n_\Pot$ and  $\mathcal{J}^n_\Pot$
defined in Proposition~\ref{Proposition 3.2} further restrict to  surjective operators $\mathcal{Q}^n_\Potc\in\mathcal{L}(\mathcal{T}^n_\Potc;\dot{\mathcal{T}}^n_\Potc)$ and
$\mathcal{J}^n_\Potc\in\mathcal{L}(\mathcal{S}^n_\Potc; \dot{\mathcal{T}}^n_\Potc)$, still satisfying
$\ker\mathcal{Q}^n_\Potc=\C_{Y_\Pot}$ and $\ker\mathcal{J}^n_\Potc=\C_{X_\Pot}$.
Consequently the bijective isomorphism $\dot{\mathcal{J}}^n_\Pot$  in Proposition~\ref{Proposition 3.2} further restricts to a bijective isomorphism $\dot{\mathcal{J}}^n_\Potc\in\mathcal{L}(\dot{\mathcal{S}}^n_\Potc; \dot{\mathcal{T}}^n_\Potc)$.
\end{cor}
\begin{proof}By Proposition~\ref{Proposition 3.2} and \eqref{3.19} we have $\mathcal{Q}^n_\Potc\in\mathcal{L}(\mathcal{T}^n_\Potc;\dot{\mathcal{T}}^n_\Potc)$ and, since $\mathcal{Q}^n_\Pot$ is surjective and by \eqref{3.40} one has $\mathcal{Q}^n_\Pot\cdot\mathcal{T}^n_\Potc=\dot{\mathcal{T}}^n_\Potc$, also $\mathcal{Q}^n_\Potc$ is surjective. Moreover, since $\C_{Y_\Pot}\subset \mathcal{T}^n_\Potc$ and $\C_{X_\Pot}\subset \mathcal{S}^n_\Potc$, from \eqref{3.35BIS} one gets $\ker\mathcal{Q}^n_\Potc=\C_{Y_\Pot}$ and  $\ker\mathcal{J}^n_\Potc=\C_{X_\Pot}$. Consequently $\mathcal{J}^n_\Potc$ subordinates a bijective isomorphism $\dot{\mathcal{J}}^n_\Potc\in\mathcal{L}(\dot{\mathcal{S}}^n_\Potc; \dot{\mathcal{T}}^n_\Potc)$, trivially being the restriction of the operator $\dot{\mathcal{J}}^n_\Pot$ in Proposition~\ref{Proposition 3.2}.
\end{proof}
\section{The Lagrangian model}\label{Section 4}
\subsection{Preliminaries} \label{Section 4.1}
We start by introducing all functional spaces needed to deal with problem $(\Lagr)$, in addition to the ones already defined in \S~\ref{Section 1.2}. We start with the Hilbert space $L^2_\Lagr=L^2(\Omega)^3\times L^2(\Gamma_1)$ with the standard inner product.
Moreover, recalling the spaces $H^n_{\curl 0}(\Omega)$, $H^n_\Lagr$,  $\mathbb{H}^n_\Lagr$ and $\mathcal{H}^n_\Lagr$ defined in \eqref{1.7}--\eqref{1.10}, in the sequel we shall use, when useful, the trivial identifications
\begin{equation}\label{4.2}
  \mathcal{H}^n_\Lagr\,=\,\mathbb{H}^n_\Lagr\times H^{n-1}_{\curl 0}(\Omega)\times H^{n-1}(\Gamma_1)\,=\,\mathbb{H}^n_\Lagr\times H^{n-1}_\Lagr.
\end{equation}
To fix the notation we also anticipate the compatibility conditions for problem $(\Lagr_0)$, which are given, for $n\in\N$, $2\le n\le r$ and for data $U_{0\Lagr}=(\rvect_0,v_0,\rvect_1,v_1)\in \mathcal{H}^n_\Lagr$, by
\begin{equation}\label{4.3}\left\{
\begin{aligned}
&\rvect_1\cdot \boldsymbol{\nu}=-v_1,\qquad \text{on $\Gamma_1$,}\qquad \rvect_1\cdot \boldsymbol{\nu}=0,\quad \text{on $\Gamma_0$,} \\
&\Delta^i\rvect_0\cdot\boldsymbol{\nu} =0,\qquad \text{on $\Gamma_0$ \qquad for $i=1,\ldots, \lfloor (n-1)/2 \rfloor $, when $n\ge 3$,}\\
&\Delta^i\rvect_1\cdot\boldsymbol{\nu} =0,\qquad \text{on $\Gamma_0$ \qquad for $i=1,\ldots, \lfloor n/2\rfloor -1$, when $n\ge 4$,}\\
&\begin{split}
\frac{B\mu}{\rho_0} \Delta^i\!\rvect_0\!\cdot\!\boldsymbol{\nu}\!=\!\DivGamma (\sigma\nabla_\Gamma (\Delta^{i-1}\!\!\rvect_0\!\cdot\!\boldsymbol{\nu}))\!
-\!\delta \Delta^{i-1}\!\!\rvect_1\!\cdot\!\boldsymbol{\nu}\!-\!\kappa\Delta^{i-1}\!\!\rvect_0\!\cdot\!\boldsymbol{\nu}\!-\! B\Div\Delta^{i-1} \!\rvect_0\\
\text{\quad on $\Gamma_1$  \quad for $i=1,\ldots, \lfloor (n-1)/2\rfloor $, \quad when $n\ge 3$,}\qquad\,
\end{split}\\
&\begin{split}
\frac{B\mu}{\rho_0} \Delta^i\!\rvect_1\!\cdot\!\boldsymbol{\nu}\!=\!\DivGamma (\sigma\nabla_\Gamma (\Delta^{i-1}\!\!\rvect_1\!\cdot\!\boldsymbol{\nu}))\!
-\!\frac{B\delta}{\rho_0} \Delta^i\rvect_0\!\cdot\!\boldsymbol{\nu}\!-\!\kappa\Delta^{i-1}\!\!\rvect_1\!\cdot\!\boldsymbol{\nu}\!-\! B\Div\Delta^{i-1} \!\!\rvect_1\\
\text{\quad on $\Gamma_1$  \quad for $i=1,\ldots, \lfloor n/2\rfloor -1$, \quad when $n\ge 4$.}\qquad\quad
\end{split}\end{aligned}\right.
\end{equation}
By \eqref{2.9} and the Trace Theorem it is straightforward to check that
\begin{equation}\label{4.4}
  D^{n-1}_\Lagr:=\{U_{0\Lagr}\in \mathcal{H}^n_\Lagr: \quad\text{\eqref{4.3} hold}\}
\end{equation}
is a closed subspace  of $\mathcal{H}^n_\Lagr$ and then a Hilbert space. When $r=\infty$ we also set the Fr\'{e}chet space
$\mathcal{H}_\Lagr^\infty = \left[C^\infty(\overline{\Omega})^3\times C^\infty(\Gamma_1)\times C^\infty(\overline{\Omega})^3\times C^\infty(\Gamma_1)\right]\cap\mathcal{H}^1_\Lagr,
$
endowed with the topology of the product space in  square brackets, and its closed subspace
\begin{equation}\label{4.6}
 D_\Lagr^\infty = \{U_{0\Lagr}\in  \mathcal{H}_\Lagr^\infty: \quad\text{\eqref{4.3} hold for all $n\in\N$}\}.
\end{equation}
We remark that, by Morrey's Theorem, $D_\Lagr^\infty=\bigcap_{n=1}^\infty D^n_\Lagr$.
Recalling the space $X^n_\Lagr$ defined in \eqref{1.11} for $n\in\widetilde{\N}$, $n\le r$, we remark that, when $r=\infty$,
\begin{multline*}
X^\infty_\Lagr=\{(\rvect,v)\in C^\infty(\R\times\overline{\Omega})\times C^\infty(\R\times\Gamma_1):  \quad\curl \rvect=0\quad\text{in $\R\times\Omega$,}\\\rvect\cdot\boldsymbol{\nu}=-v\quad\text{on $\R\times\Gamma_1$,}\qquad \rvect\cdot\boldsymbol{\nu}=0\quad\text{on $\R\times\Gamma_0$}\}.
\end{multline*}
For the same $n$'s we also set the Fr\'{e}chet space
\begin{equation}\label{4.9}
Y^n_\Lagr=\bigcap_{i=0}^{n-1} C^i(\R;\mathcal{H}^{n-i}_\Lagr).
\end{equation}

At first we essentially recall a result in \cite{dautraylionsvol3}.
\begin{lem}\label{Lemma 4.1}
The operator $G\in\mathcal{L}(\dot{H}^1(\Omega); H^0_{\curl 0}(\Omega))$ defined by
\begin{equation}\label{4.10}
G\dot{u}=-\nabla \dot{u}  \qquad\text{for all $\dot{u}\in \dot{H}^1(\Omega)$}
\end{equation}
restricts for all $n\in\N$ to a bijective isomorphism between $\dot{H}^n(\Omega)$ and $H^{n-1}_{\curl 0}(\Omega)$.
\end{lem}
\begin{proof} $\Gamma$ being compact it has finitely many connected components. Moreover $\Omega$ is bounded, simply connected and of class $C^r$, $r\ge 2$. Hence it satisfies \cite[Chapter IX, \S 1, assumption (1.45) p. 217]{dautraylionsvol3}. In the authors' notation $N=0$. Hence one can  apply \cite[Chapter IX, \S 1, Proposition 2, p. 219]{dautraylionsvol3} with $\mathbb{H}_1=\{0\}$, this space being isomorphic to the first cohomology space. We then get that $H^0_{\curl 0}(\Omega)=\{\nabla\varphi, \, \varphi\in H^1(\Omega)\}$,
from which $H^{n-1}_{\curl 0}(\Omega)=\{\nabla\varphi, \, \varphi\in H^n(\Omega)\}$ trivially follows. Since $\ker \nabla=\C_\Omega$ in $H^1(\Omega)$, the proof is complete.
\end{proof}
The key result to connect problems $(\Potc)$ and $(\Lagr)$ is the following one.
\begin{lem}\label{Lemma 4.2}
For all $w\in L^2(\Omega)$ and $v\in H^{1/2}(\Gamma_1)$ satisfying the compatibility condition
\begin{equation}\label{4.11}
 \rho_0\int_\Omega w=B\int_{\Gamma_1}v,
\end{equation}
the problem
\begin{equation}\label{4.12}
  \begin{cases}
  -B\Div \rvect=\rho_0 w,\quad &\text{in $\Omega$,}\\
  \quad \curl\rvect=0\quad &\text{in $\Omega$,}\\
\quad \rvect\cdot\boldsymbol{\nu}=-v\quad &\text{on $\Gamma_1$,}\\
\quad \rvect\cdot\boldsymbol{\nu}=0\quad &\text{on $\Gamma_0$,}
  \end{cases}
\end{equation}
has a unique solution $\rvect\in H^1_{\curl 0}(\Omega)$. Consequently, recalling the spaces defined in \eqref{1.9} and \eqref{3.2},
the map $S(v,w)=(\rvect,v)$ for all $(v,w)\in H^1_\Potc$,
defines a bijective isomorphism $S\in\mathcal{L}(H^1_\Potc;\mathbb{H}^1_\Lagr)$,  with inverse $S^{-1}\in\mathcal{L}(\mathbb{H}^1_\Lagr;H^1_\Potc)$ given by
$S^{-1}(\rvect,v)=\left(v, - B\Div \rvect/\rho_0\right)$.
Moreover, for all $n\in\N$, $2\le n\le r$, and $q\in [2,\infty)$, there is a positive constant $c_{n,q}=c_{n,q}(\Omega)$ such that, for all $w\in W^{n-1,q}(\Omega)$ and $v\in W^{n-1/q,q}(\Gamma_1)$ satisfying \eqref{4.11} one has
\begin{equation}\label{4.15}
  \|\rvect\|_{[W^{n,q}(\Omega)]^3}\le c_{n,q}\left(\|w\|_{W^{n-1,q}(\Omega)}+\|v\|_{v\in W^{n-1/q,q}(\Gamma_1)}\right).
\end{equation}
Consequently $S\in\mathcal{L}(H^n_\Potc;\mathbb{H}^n_\Lagr)$ and $S^{-1}\in\mathcal{L}(\mathbb{H}^n_\Lagr;H^n_\Potc)$.
Finally, for all $w\in C^{r-1}(\overline{\Omega})$ and $v\in C^{r-1}(\Gamma_1)$ satisfying \eqref{4.11} the solution $\rvect$ of \eqref{4.12} can be extended to $\rvect\in [C^{r-1}(\R^3)]^3$.
\end{lem}
\begin{proof}
By Lemma~\ref{Lemma 4.1} solving problem \eqref{4.12} in the space $H^1_{\curl 0}(\Omega)$ it is equivalent to solving the inhomogeneous Neumann problem
\begin{equation}\label{4.16}
  \begin{cases}
  -B\Delta \varphi=\rho_0 w,\quad &\text{in $\Omega$,}\\
  \quad \partial_{\boldsymbol{\nu}}\varphi=-v\quad &\text{on $\Gamma_1$,}\\
\quad \partial_{\boldsymbol{\nu}}\varphi=0\quad &\text{on $\Gamma_0$,}
  \end{cases}
\end{equation}
with $\varphi\in \dot{H}^2(\Omega)$, and hence take $\rvect=\nabla\varphi$. Using the splitting \eqref{2.8} problem \eqref{4.16} can be rewritten as
$-B\Delta  \varphi=\rho_0 w$  in $\Omega$, $\partial_{\boldsymbol{\nu}}\varphi=-v$ on $\Gamma$,
with $v\in H^{1/2}(\Gamma)$, and  $v\in W^{n-1/q,q}(\Gamma_1)$  is rewritten as $v\in W^{n-1/q,q}(\Gamma)$.
By standard elliptic theory this   problem  has a unique solution $\varphi\in \dot{H}^2(\Omega)$ and, moreover, if also $w\in W^{n-1,q}(\Omega)$  and
$v\in W^{n-1/q,q}(\Gamma)$, one has $\varphi\in W^{n+1,q}(\Omega)/\C_\Omega$ and there is a positive constant
$c_{n,q}=c_{n,q}(\Omega)$ such that
$\|\varphi\|_{W^{n+1,q}(\Omega)/\C_\Omega}\le c_{n,q}\left(\|w\|_{W^{n-1,q}(\Omega)}+\|v\|_{v\in W^{n-1/q,q}(\Gamma_1)}\right)$.

We thus found the unique solution $\rvect =\nabla\varphi\in H^1_{\curl 0}(\Omega)$ of \eqref{4.12} and \eqref{4.15} holds. Taking $q=2$ in it we then get  that $S\in\mathcal{L}(H^n_\Potc;\mathbb{H}^n_\Lagr)$ for $n\in\N$, $n\le r$, with inverse being trivially given as in the statement, where the Divergence Theorem shows that $S^{-1}\in\mathcal{L}(\mathbb{H}^n_\Lagr;H^n_\Potc)$.
Let now assume that $w\in C^{r-1}(\overline{\Omega})$, $v\in C^{r-1}(\Gamma_1)$ and  \eqref{4.11} holds. When $r\in\N$, since $\Gamma_1$ is compact, by \eqref{4.15} we get $\rvect\in \bigcap_{q\in [2,\infty)} W^{r,q}(\Omega)^3$. By Stein Extension Theorem (see \cite[Chapter 5, Theorem 5.24, p. 154]{adamsfournier}) we can then extend $\rvect$ to $\rvect\in \bigcap_{q\in [2,\infty)} W^{r,q}(\R^3)^3$ and, by Morrey's Theorem, $\rvect\in C^{r-1}(\R^3)^3$. Since the extension operator in Stein Extension Theorem is independent of $r$, the case $r=\infty$ follows from the previous one.
\end{proof}
As a byproduct of Lemma~\ref{Lemma 4.2} we can now derive the following non--trivial density result, which will be useful in the sequel.
\begin{lem}\label{Lemma 4.3}
The space $\mathbb{Y}=\{(\rvect_{|\Omega},v)\in \mathbb{H}^1_\Lagr: \rvect\in C^{r-1}(\R^3)^3\}$ is dense in $\mathbb{H}^1_\Lagr$.
\end{lem}
\begin{proof} Using Lemma~\ref{Lemma 4.2} we can endow $\mathbb{H}^1_\Lagr$ with the inner product
$((\rvect,v), (\phivect,\psi))_{\mathbb{H}^1_\Lagr}:=\int_\Omega \Div\rvect\Div\overline{\phivect} + \int_{\Gamma_1} (\nabla_\Gamma v,\nabla_\Gamma \psi)_\Gamma+\int_{\Gamma_1} v\overline{\psi}$,
which induces a norm equivalent to the one inherited from $H^1(\Omega)^3\times H^1(\Gamma_1)$. To prove our assertion
we then take $(\rvect_0,v_0)\in\mathbb{H}^1_\Lagr$ such that
\begin{equation}\label{4.19}
 \int_\Omega \Div\rvect_0\Div\overline{\phivect} + \int_{\Gamma_1} (\nabla_\Gamma v_0,\nabla_\Gamma \psi)_\Gamma+\int_{\Gamma_1} v_0\overline{\psi}=0\quad\text{for all $(\phivect,\psi)\in \mathbb{Y}$,}
\end{equation}
claiming that $(\rvect_0,v_0)=0$.
We now take $\phivect_1\in \mathcal{D}(\Omega)^3$, so $\int_\Omega \Div \phivect_1=0$. By Lemma~\ref{Lemma 4.2} then there is $\phivect\in C^{r-1}(\R^3)^3$ such that $\Div\phivect=\Div\phivect_1$ and $\curl \phivect=0$ in $\Omega$, with $\phivect\cdot\boldsymbol{\nu}=0$ on $\Gamma$. Consequently $(\phivect_{|\Omega},0)\in \mathbb{Y}$ and, by \eqref{4.19}, we have
\begin{equation}\label{4.20}
 \int_\Omega \Div\rvect_0\Div\overline{\phivect_1} =\int_\Omega \Div\rvect_0\Div\overline{\phivect}= 0.
 \end{equation}
Since $\phivect_1\in \mathcal{D}(\Omega)^3$ is arbitrary, taking $\phivect_1=(\varphi,0,0), (0,\varphi,0), (0,0,\varphi)$ with $\varphi\in\mathcal{D}(\Omega)$ we then get that $\Div \rvect_0\in H^1(\Omega)$ and $\nabla \Div \rvect_0=0$ in $\Omega$ in the sense of distributions.  $\Omega$ being connected we thus have $\Div\rvect_0\in \C_\Omega$. Moreover, since $(\rvect_0,v_0)\in \mathbb{H}^1_\Lagr$, by \eqref{1.9} and the Divergence Theorem we have
\begin{equation}\label{4.21}
 \Div\rvect_0=-\tfrac 1{|\Omega|}\int_{\Gamma_1}v_0\qquad\text{in $\Omega$,}
\end{equation}
where $|\Omega|$ denotes the volume of $\Omega$. Plugging \eqref{4.21} into \eqref{4.19}, since $(\phivect,\psi)\in \mathbb{H}^1_\Lagr$, using the Divergence Theorem again we get
\begin{equation}\label{4.22}
 \frac 1{|\Omega|}\int_{\Gamma_1}v_0\int_{\Gamma_1}\overline{\psi} + \int_{\Gamma_1} (\nabla_\Gamma v_0,\nabla_\Gamma \psi)_\Gamma+\int_{\Gamma_1} v_0\overline{\psi}=0\quad\text{for all $(\phivect,\psi)\in \mathbb{Y}$.}
\end{equation}
We now claim that \eqref{4.22} actually holds for all $\psi\in H^1(\Gamma_1)$. To prove our claim we can assume, by density, that $\psi\in C^r(\Gamma_1)$. Applying Lemma~\ref{Lemma 4.2} again with $w=\tfrac B{\rho_0|\Omega|}\int_{\Gamma_1}\psi$ and $v=\psi$, so $w\in C^r(\overline{\Omega})$, $v\in C^r(\Gamma_1)$ and $\rho_0\int_\Omega w-B\int_{\Gamma_1}v=0$, we then get that there is $\phivect\in [C^{r-1}(\R^3)]^3$ such that $\curl \phivect=0$ in $\Omega$, $\phivect\cdot\boldsymbol{\nu}=-\psi$ on $\Gamma_1$ and $\phivect\cdot\boldsymbol{\nu}=0$ on $\Gamma_0$, so that $(\phivect,\psi)\in \mathbb{Y}$ and consequently \eqref{4.22} holds, proving our claim. Using it we can take $\psi=v_0$ in \eqref{4.22}, getting $v_0=0$ and, by \eqref{4.21}, $\Div \rvect_0=0$. Then $((\rvect_0,v_0),(\rvect_0,v_0))_{\mathbb{H}^1_\Lagr}=0$ and consequently $(\rvect_0,v_0)=0$.
\end{proof}
Our final preliminary result is a second non--trivial density result. To state it we introduce the Hilbert space
\begin{equation}\label{4.23}
  V_\Lagr=\{(\rvect,v)\in H^1(\Omega)^3\times H^1(\Gamma_1):\quad \rvect\cdot\boldsymbol{\nu}=-v\quad\text{on $\Gamma_1$,}
\quad \rvect\cdot\boldsymbol{\nu}=0\quad\text{on $\Gamma_0$}\}
\end{equation}
equipped with the norm inherited from the product. Trivially $V_\Lagr\hookrightarrow L^2_\Lagr$.
\begin{lem}\label{Lemma 4.4}
The embedding $V_\Lagr\hookrightarrow L^2_\Lagr$ is dense.
\end{lem}
\begin{proof}
Using the standard inner product of $L^2_\Lagr$, to prove our assertion we take $(\rvect_0,v_0)\in L^2_\Lagr$ such that
\begin{equation}\label{4.24}
\int_\Omega \rvect_0\cdot\overline{\phivect} + \int_{\Gamma_1} v_0\overline{\psi}=0\quad\text{for all $(\phivect,\psi)\in V_\Lagr$}
\end{equation}
and we claim that then $(\rvect_0,v_0)=0$. First taking in \eqref{4.24} test functions $(\phivect,\psi)=(\phivect_1,0)$, with $\phivect_1\in \mathcal{D}(\Omega)^3$, so trivially $(\phivect,\psi)\in V_\Lagr$, we get $\int_\Omega\rvect_0\cdot\phivect_1=0$ and consequently, being $\phivect_1$ arbitrary, $\rvect_0=0$. Hence \eqref{4.24} reads as
\begin{equation}\label{4.25}
\int_{\Gamma_1} v_0\overline{\psi}=0\quad\text{for all $(\phivect,\psi)\in V_\Lagr$}
\end{equation}
We now claim that \eqref{4.25} actually holds for all $\psi\in H^1(\Gamma_1)$, from which one immediately gets $v_0=0$ and concludes the proof. To prove our claim we apply Lemma~\ref{Lemma 4.2} once again with $w=\tfrac B{\rho_0|\Omega|}\int_{\Gamma_1}\psi$ and $v=\psi$. It follows that there is $\phivect\in H^1_{\curl 0}(\Omega)$ such that $(\phivect,\psi)\in V_\Lagr$, so \eqref{4.25} hold true and our claim is proved.
\end{proof}
\subsection{The main isomorphism}\label{Section 4.2} To deal with problem $(\Lagr)$ we first construct  the main isomorphism between its phase space $\mathcal{H}^1_\Lagr$ and the phase space $\mathcal{H}^1_\Potc$.
We set
\begin{equation}\label{4.26}
  F_\LPc (\rvect,v,\svect,z)=\left(\dot{u}, v, -B\Div \rvect/{\rho_0},z\right)\qquad\text{for all $(\rvect,v,\svect,z)\in\mathcal{H}^1_\Lagr$,}
\end{equation}
where $\dot{u}\in \dot{H}^1(\Omega)$ is the unique solution of the equation
\begin{equation}\label{4.27}
-\nabla\dot{u}=\svect
\end{equation}
given by Lemma~\ref{Lemma 4.1}. By \eqref{1.6}, \eqref{1.10} and the Divergence Theorem we have
$F_\LPc\in\mathcal{L}(\mathcal{H}^1_\Lagr; \dot{\mathcal{H}}^1_\Potc)$.
We also set
\begin{equation}\label{4.29}
 F_\PcL(\dot u, v,w,z)=(\rvect,v,-\nabla \dot u,z)\qquad\text{for all $(\dot u,v,w,z)\in\mathcal{H}^1_\Potc$,}
\end{equation}
where $\rvect\in H^1_{\curl 0}$ is the unique solution of problem \eqref{4.12} given by Lemma~\ref{Lemma 4.2}. Since
for any $\dot u\in \dot H^1(\Omega)$ we have $\curl\nabla\dot u=0$ in $\mathcal{D}'(\Omega)$, see \cite[Chapter 3, p. 202]{dautraylionsvol3}, one trivially gets
$F_\PcL\in \mathcal{L}(\dot{\mathcal{H}}^1_\Potc;\mathcal{H}^1_\Lagr)$.
By \eqref{3.15} and \eqref{4.9} the operators $F_\LPc$ and $F_\PcL$ trivially induce the operators
$\Phi_\LPc \in\mathcal{L}(Y^1_\Lagr; \dot{Y}^1_\Potc)$ and $\Phi_\PcL \in\mathcal{L}(\dot{Y}^1_\Potc; Y^1_\Lagr)$
given by
\begin{equation}\label{4.32}
 (\Phi_\LPc U_\Lagr)(t)=F_\LPc U_\Lagr (t),  \quad(\Phi_\PcL \dot{U}_\Pot)(t)=F_\PcL \dot{U}_\Pot (t)\quad \text{for all $t\in\R$.}
\end{equation}
The following result points out all properties of these operators needed in the sequel.
\begin{prop}\label{Proposition 4.1}
The operator $F_\LPc$ is a bijective isomorphism between $\mathcal{H}^1_\Lagr$ and $\dot{\mathcal{H}}^1_\Potc$ having
$F_\PcL$ as inverse. Consequently $\Phi_\LPc$ is a bijective isomorphism between $Y^1_\Lagr$ and $\dot{Y}^1_\Potc$ having
$\Phi_\PcL$ as inverse. Moreover, for all $n\in\N$, $2\le n\le r$, they restrict to
\begin{equation}\label{4.33}
\begin{aligned}
F_\LPc\in\mathcal{L}(\mathcal{H}^n_\Lagr; \dot{\mathcal{H}}^n_\Potc), \qquad & F_\PcL\in \mathcal{L}(\dot{\mathcal{H}}^n_\Potc;\mathcal{H}^n_\Lagr),\\
\Phi_\LPc \in\mathcal{L}(Y^n_\Lagr; \dot{Y}^n_\Potc),\qquad &
\Phi_\PcL \in\mathcal{L}(\dot{Y}^n_\Potc; Y^n_\Lagr),
\end{aligned}
\end{equation}
and, since
\begin{equation}\label{4.34}
  F_\LPc D^{n-1}_\Lagr=\dot{D}^{n-1}_\Potc,
\end{equation}
$F_\LPc$ and $F_\PcL$ also restrict to
$F_\LPc\in\mathcal{L}(D^{n-1}_\Lagr; \dot{D}^{n-1}_\Potc)$  and $F_\PcL\in \mathcal{L}(\dot{D}^{n-1}_\Potc;D^{n-1}_\Lagr)$.
\end{prop}
\begin{proof} To recognize that $F_\LPc$ and $F_\PcL$ are the inverse of each  other we simply recall that the identifications \eqref{4.2} and $\dot{\mathcal H}_\Potc^n=\dot{H}^n(\Omega)\times H_\Potc^n\times H^{n-1}(\Gamma_1)$   when $n=1$ read as
$\mathcal{H}^1_\Lagr=\mathbb{H}^1_\Lagr\times H^0_{\curl 0}(\Omega)\times L^2(\Gamma_1)$, and $\dot{\mathcal H}_\Potc^1=\dot{H}^1(\Omega)\times H_\Potc^1\times L^2(\Gamma_1)$.
Using them together with \eqref{4.26} and \eqref{4.29} we can represent the operators $F_\LPc$ and $F_\PcL$ in the matrix form
\begin{equation}\label{4.36}
 F_\LPc= \begin{pmatrix}0 & G^{-1}  & 0\\S^{-1} & 0 & 0 \\0 & 0 &\text{Id}\end{pmatrix}\qquad\text{and}\quad
 F_\PcL= \begin{pmatrix}0 & S  & 0\\G & 0 & 0 \\0 & 0 &\text{Id}\end{pmatrix},
\end{equation}
where $G$ and $S$ are the bijective isomorphism introduced in Lemmas~\ref{Lemma 4.1} and \ref{Lemma 4.2}, so trivially $F_\PcL=F_\LPc^{-1}$. The asserted properties of $\Phi_\LPc$ and $\Phi_\PcL$ then follow trivially.  Moreover, for $n\in\N$, $2\le n\le r$, using \eqref{4.36}, Lemmas~\ref{Lemma 4.1} -- \ref{Lemma 4.2} and recalling \eqref{3.15}, \eqref{4.9}, we get \eqref{4.33}.
Next,  from \eqref{4.33} and \eqref{4.34} it trivially follow that $F_\LPc\in\mathcal{L}(D^{n-1}_\Lagr; \dot{D}^{n-1}_\Potc)$ and $F_\PcL\in \mathcal{L}(\dot{D}^{n-1}_\Potc;D^{n-1}_\Lagr)$,  so we just have to prove \eqref{4.34} to complete the proof.

By \eqref{3.5}, \eqref{3.6}, \eqref{4.4} and \eqref{4.26} -- \eqref{4.27} to prove \eqref{4.34} reduces to prove that for all $(\rvect_0,v_0,\rvect_1,v_1)\in \mathcal{H}^n_\Lagr$ and $(u_0,v_0,u_1,v_1)\in \dot{\mathcal{H}}^n_\Potc$ such that
\begin{equation}\label{4.37}
\rho_0u_1=-B\Div \rvect_0\qquad\text{and}\quad \nabla u_0=-\rvect_1
\end{equation}
the compatibility conditions \eqref{3.4} and \eqref{4.3} are equivalent. To prove this fact
 we preliminarily remark that one easily gets  that
for all $\rvect\in H^0_{\curl 0}(\Omega)$ and $i\in\N$
$\curl \Delta ^i\rvect=0$  and $\Delta^i\rvect=(\nabla\Div)^i\rvect$ in $\mathcal{D}'(\Omega)^3$,
Then, starting from \eqref{4.37} and using it one gets that, in $\Omega$, one has
$\Delta^iu_0=-\Div\Delta^{i-1}\rvect_1$ for $i=1,\ldots,[n/2]$, $\nabla \Delta^iu_0=-\Delta^i\rvect_1$ for $i=0,\ldots,[(n-1)/2]$,
$\rho_0\Delta^iu_1=-B\Div\Delta^i\rvect_0$, for $i=0,\ldots,[(n-1)/2]$ and $\rho_0\nabla \Delta^iu_1=-B\Delta^{i+1}\rvect_0$ for $i=0,\ldots,[n/2]-1$. Passing these relations to traces one easily checks that \eqref{3.4} and \eqref{4.3} are equivalent.
\end{proof}
\subsection{Abstract analysis of problems $(\Lagr)$ and $(\Lagr_0)$}\label{Section 4.3}
To deal with problem $(\Lagr_0)$ in a semigroup setting we  reduce  it  to the first order problem
\begin{equation}\label{4.40}
\begin{cases}
\rvect_t=\svect\qquad &\text{in
$\R\times\Omega$,}\\
v_t =z\qquad
&\text{on
$\R\times \Gamma_1$,}\\
\rho_0\svect_t-B\nabla \Div \rvect=0\qquad &\text{in
$\R\times\Omega$,}\\
\curl\rvect=0\qquad &\text{in
$\R\times\Omega$,}\\
\mu z_t- \DivGamma (\sigma \nabla_\Gamma v)+\delta z+\kappa v-B\Div \rvect =0\qquad
&\text{on
$\R\times \Gamma_1$,}\\
\rvect\cdot{\boldsymbol{\nu}} =0 \quad\text{on $\R\times \Gamma_0$,}\qquad
\rvect\cdot{\boldsymbol{\nu}} =-v\qquad
&\text{on
$\R\times \Gamma_1$,}\\
\rvect(0,x)=\rvect_0(x),\quad \svect(0,x)=\rvect_1(x) &
 \text{in $\Omega$,}\\
v(0,x)=v_0(x),\quad z(0,x)=v_1(x) &
 \text{on $\Gamma_1$.}
\end{cases}
\end{equation}
More formally, working in phase space $\mathcal{H}^1_\Lagr$ in which equations \eqref{4.40}$_4$ and \eqref{4.40}$_6$ implicitly hold, we introduce the unbounded operator $A_\Lagr: D(A_\Lagr)\subset \mathcal{H}^1_\Lagr \to \mathcal{H}^1_\Lagr$ given by
\begin{gather}\label{4.41}
D(A_\Lagr)=D^1_\Lagr=\{(\rvect,v,\svect,z)\in \mathcal{H}^2_\Lagr: (\svect,z)\in \mathbb{H}^1_\Lagr\}\\
\label{4.42}
A_\Lagr\begin{pmatrix}\rvect\\v\\\svect\\z\end{pmatrix} =
\begin{pmatrix}-\svect\\-z\\-(B/\rho_0)\nabla\Div \rvect\\
\frac 1\mu\left[-\DivGamma(\sigma\nabla_\Gamma v)+\delta z+\kappa v-B\Div\rvect_{|\Gamma_1}\right]
\end{pmatrix},
\end{gather}
together with the abstract equation and Cauchy problem
\begin{gather}\label{4.43}U_\Lagr'+A_\Lagr U_\Lagr=0\qquad\text{in $\mathcal{H}^1_\Lagr$,}
\\\label{4.44}
U_\Lagr'+A_\Lagr U_\Lagr=0\qquad\text{in $\mathcal{H}^1_\Lagr$}, \qquad U_\Lagr(0)=U_{0\Lagr}\in \mathcal{H}^1_\Lagr.
\end{gather}
The following result shows that problem \eqref{4.44}  is essentially equivalent to problem \eqref{3.32} restricted to $\dot{\mathcal{H}}^1_\Potc$ and  hence it is well--posed.
\begin{thm}[\bf Well--posedness for \eqref{4.44}]\label{Theorem 4.1}
The operator $-A_\Lagr$ is densely defined and it generates on $\mathcal{H}^1_\Lagr$ the strongly continuous
group $\{\Lagr^1(t),t\in\R\}$ given by
\begin{equation}\label{4.45}
  \Lagr^1(t)=F_\PcL\Pot_c^1(t) F_\LPc\qquad\text{for all $t\in\R$,}
\end{equation}
and hence similar to the group $\{\Pot_c^1(t),t\in\R\}$ in Corollary~\ref{Corollary 3.3}.
Consequently, for any $U_{0\Lagr}\in \mathcal{H}^1_\Lagr$, problem \eqref{4.44} has a unique generalized solution $U_\Lagr\in Y^1_\Lagr$ given by $U_\Lagr(t)=\Lagr^1(t)[U_{0\Lagr}]$ for all $t\in\R$ and hence continuously depending on $U_{0\Lagr}$ in the topologies of the respective spaces. Moreover $U_\Lagr$ is a strong solution if and only if $U_{0\Lagr}\in D(A_\Lagr)$.
Finally, if $\dot{U}_\Pot\in \dot{Y}^1_\Potc$ is the unique generalized solution of problem \eqref{3.32} with data $\dot{U}_{0\Pot}=F_\LPc U_{0\Lagr}$, one has $U_\Lagr=\Phi_\PcL \dot{U}_\Potc$.
\end{thm}
\begin{proof}By standard semigroup theory  formula \eqref{4.45} defines on $\mathcal{H}^1_\Lagr$ a strongly continuous group $\{\Lagr^1(t),t\in\R\}$, similar to $\{\Pot_c^1(t),t\in\R\}$, having as generator the operator $-B_1$ defined by $D(B_1)=F_\PcL D(\dot{A}_\Potc)$ (so $B_1$ is densely defined) and $B_1=F_\PcL\dot{A}_\Potc F_\LPc$. We now claim that $B_1=A_\Lagr$. By combining \eqref{3.38}, \eqref{4.41} and \eqref{4.34} we get $D(B_1)=D^1_\Lagr=D(A_\Lagr)$. Moreover, using  \eqref{4.26} and \eqref{4.27}, for any $U=(\rvect,v,\svect,z)\in D(A_\Lagr)$ we have
$$B_1U=F_\PcL\dot{A}_\Pot F_\LPc U=F_\PcL\dot{A}_\Pot (\dot u,v,-B\Div\rvect/\rho_0,z)^t$$
where $\dot{u}\in \dot{H}^2(\Omega)$ is the unique solution of $-\nabla\dot u=\svect$.
Consequently, using also \eqref{3.30}, \eqref{3.38} and \eqref{4.29}
one gets
\begin{equation}\label{4.46}
B_1 U=\begin{pmatrix}\tvect\\-z\\-\tfrac B{\rho_0}\nabla \Div \rvect\\
\frac 1\mu\left[-\DivGamma(\sigma\nabla_\Gamma v)+\delta z+\kappa v-B\Div\rvect_{|\Gamma_1}\right]
\end{pmatrix}
\end{equation}
where $\tvect\in H^1_{\curl 0}(\Omega)$ is the unique solution of the problem
$-B\Div \tvect=B\Div\svect$ in $\Omega$, $\curl\tvect=0$  in $\Omega$, $\tvect\cdot\boldsymbol{\nu}=z$ on $\Gamma_1$, $\tvect\cdot\boldsymbol{\nu}=0$ on $\Gamma_0$.
Since, by \eqref{4.41}, $(\svect,z)\in\mathbb{H}^1_\Lagr$, using Lemma~\ref{Lemma 4.2} we get $\tvect=-\svect$, so by \eqref{4.42} and \eqref{4.46} we have $B_1U=A_\Lagr U$, proving our claim. The proof can then be completed by using Lemma~\ref{Lemma 2.1}-- i) and iii).
\end{proof}
We now set the Fr\'{e}chet spaces of generalized solutions of \eqref{4.43} that is, for $n\le r$,
\begin{equation}\label{4.47}
\mathcal{T}^n_\Lagr=\{U_\Lagr\in Y^n_\Lagr: \quad \text{$U_\Lagr$ is a generalized solution of \eqref{4.43}}\},
\end{equation}
endowed with the topology inherited from $Y^n_\Lagr$. We point out that, as a trivial consequence of Theorem~\ref{Theorem 4.1}, the following result holds.
\begin{cor}\label{Corollary 4.1}
The operators $\Phi_\LPc$ and $\Phi_\PcL$ restrict to bijective isomorphisms between $\mathcal{T}^n_\Lagr$ and $\dot{\mathcal{T}}^n_\Potc$ for all $n\in\widetilde{\N}$, $n\le r$, being the inverse of each other.
\end{cor}
Abstract regularity properties of solutions of \eqref{4.44} are then given as follows.
\begin{thm}[\bf Regularity for \eqref{4.44}]\label{Theorem 4.2}
For all $n\in\widetilde{\N}$, $2\le n\le r$, one has $D(A^{n-1}_\Lagr)=D^{n-1}_\Lagr$, the respective norms being equivalent, so the operator $-\,\, {_n}A_\Lagr$ given by \eqref{2.18} generates on $D^{n-1}_\Lagr$ the strongly continuous group $\{\Lagr^n(t), t\in\R\}$ given by
\begin{equation}\label{4.48}
\Lagr^n(t)=\Lagr^1(t)_{|D^{n-1}_\Lagr}= F_\PcL\Pot^n_c(t)F_\LPc\qquad\text{for all $t\in\R$,}
\end{equation}
and hence similar to the group $\{\Pot^n_c(t),t\in\R\}$ in Corollary~\ref{Corollary 3.3}.

Moreover, for any $U_{0\Lagr}\in\mathcal{H}^1_\Lagr$, denoting by $U_\Lagr$ the generalized solution of \eqref{4.44}, for any $n\in\widetilde{\N}$, $2\le n\le r$,  $U_{0\Lagr}\in D^{n-1}_\Lagr$ if and only if $U_\Lagr\in Y^n_\Lagr$
 and in this case $U_\Lagr$  continuously depends on  $U_{0\Lagr}$ in the topologies of the respective spaces.
\end{thm}
\begin{proof} In the proof of Theorem~\ref{Theorem 4.1} we got that $A_\Lagr=F_\PcL\dot{A}_\Potc F_\LPc$. Since, by Corollary~\ref{Corollary 3.3}, $D(\dot{A}^{n-1}_\Potc)=\dot{D}^{n-1}_\Potc$, by Proposition~\ref{Proposition 4.1} and \eqref{2.16} one easily gets by induction that $D(A^{n-1}_\Lagr)=D^{N-1}_\Lagr$, with equivalence of norms.  By Lemma~\ref{Lemma 2.1}--ii) then $-{_n}A_\Lagr$ generates on $D^{n-1}_\Lagr$ the strongly continuous group  $\{\Lagr^n(t), t\in\R\}$ given by $\Lagr^n(t)=\Lagr^1(t)_{|D^{n-1}_\Lagr}$ for all $t\in\R$. By \eqref{4.45} and Proposition~\ref{Proposition 4.1} then, since $\Pot^n_c(t)=\Pot^1_c(t)_{|\dot{D}^{n-1}_\Potc}$, we get formula \eqref{4.48} and the group similarity. The proof is then completed by using Corollary~\ref{Corollary 4.1} and Lemma~\ref{Lemma 2.1}--iv) to get the asserted continuous dependence when $n\in\N$, while the case $n=\infty$ follows from  the case $n\in\N$ and \eqref{4.9}.
\end{proof}
\subsection{Solutions of $(\Lagr)$ and $(\Lagr_0)$}\label{Section 4.4} To apply the abstract results in \S~\ref{Section 4.3} to problems $(\Lagr)$ and $(\Lagr_0)$ we first make precise which type of solutions we shall consider in the sequel, recalling that $(\Lagr)_2$ and $(\Lagr)_4$ are implicit in  $X^1_\Lagr$.
\begin{definition}\label{Definition 4.1}
We say that
\renewcommand{\labelenumi}{{\roman{enumi})}}
\begin{enumerate}
\item $(\rvect,v)\in X^2_\Lagr$ is a {\em strong solution} of $(\Lagr)$ provided $(\Lagr)_1$ holds a.e. in $\R\times\Omega$ and
$(\Lagr)_3$ hold a.e. on $\R\times\Gamma_1$, where $\Div \rvect$ on $\R\times\Gamma_1$ is taken in the pointwise trace sense given in \S~\ref{Section 2.2};
\item $(\rvect,v)\in X^1_\Lagr$ is a {\em generalized solution} of $(\Lagr)$ provided it is the limit in $X^1_\Lagr$ of a sequence of strong solutions of it;
\item $(\rvect,v)\in X^1_\Lagr$ is a {\em weak solution} of $(\Lagr)$ provided the distributional identity
\begin{multline}\label{4.49}
 \int_{-\infty}^\infty\int_\Omega [\rho_0\rvect_t\cdot \phivect_t-B\Div \rvect\Div \phivect]\\
+ \int_{-\infty}^\infty\int_{\Gamma_1}\left[\mu v_t\psi_t-\sigma(\nabla_\Gamma v,\nabla_\Gamma \overline{\psi})_\Gamma-\delta v_t\psi
 -\kappa v\psi\right]=0
 \end{multline}
holds for all $\phivect\in C^{r-1}_c(\R\times\R^3)^3$ such that $\phivect\cdot\boldsymbol{\nu}=0$ on $\R\times\Gamma_0$, where $\psi=-\phivect\cdot\boldsymbol{\nu}$ on $\R\times\Gamma_1$;
\item $(\rvect,v)$ is a strong, generalized or weak solution of problem $(\Lagr_0)$ provided it is a solution of $(\Lagr)$ of the same type and $(\Lagr)_5$ -- $(\Lagr)_6$ hold in the space $X^1_\Lagr$.
\end{enumerate}
\end{definition}
\begin{rem}\label{Remark 4.0} Definition~\ref{Definition 4.1}--iii) seems to depend on $r$ (i.e. on the regularty of $\Gamma$), but its independence on it is a trivial consequence of Theorem~\ref{Theorem 1.2}.
\end{rem}
Trivially strong solutions in Definition~\ref{Definition 4.1} are also generalized ones. Moreover strong and generalized solutions correspond to the homologous solutions of \eqref{4.43} or \eqref{4.4}, as the following results shows.
\begin{lem}\label{Lemma 4.5}
The couple $(\rvect,v)$ is a strong or generalized solution of $(\Lagr)$ if and only if $\rvect$ and $v$ are the first two components of a solution $U_\Lagr=(\rvect,v,\svect,z)$ of \eqref{4.43} of the same type. In this case $\rvect_t=\svect$ and $v_t=z$. The same relation occurs between solutions of $(\Lagr_0)$ and \eqref{4.44}.
\end{lem}
\begin{proof}
Since, by Theorems~\ref{Theorem 4.1} and \ref{Theorem 4.2}, strong solutions of \eqref{4.42} belong to $Y^2_\Lagr$, using \eqref{4.42}, the relation between strong solutions of $(\Lagr)$ and \eqref{4.43} is trivial.
The relation then extends  to generalized solutions by a standard density argument.
\end{proof}
By Lemma~\ref{Lemma 4.5} generalized and strong solutions of $(\Lagr)$ naturally arise from Theorem~\ref{Theorem 4.1}.
While strong solutions are a.e. classical solutions, generalized solutions solve $(\Lagr)_1$ and $(\Lagr)_3$ in a quite indirect sense. Since couples $(\rvect,v)\in X^1_\Lagr$ are not regular enough to be a.e. solutions, a classical way to characterize them would be to consider equations $(\Lagr)_1$ and $(\Lagr)_3$ in a distributional sense, that is in $\mathcal{D}'(\R\times\Omega)^3$ and in $[C^r_c(\R\times\Gamma_1)]'$.
On the other hand, while for $(\rvect,v)\in X^1_\Lagr$ equation $(\Lagr)_1$ has in $\mathcal{D}'(\R\times\Omega)^3$ the natural form
\begin{equation}\label{4.50}
 \int_{-\infty}^\infty\int_\Omega \rho_0\rvect_t\cdot \phivect_t-B\Div \rvect\Div \phivect=0\qquad\text{for all $\phivect\in [\mathcal{D}(\R\times\Omega)]^3$,}
 \end{equation}
 equation $(\Lagr)_3$ can not be written in a distributional sense unless the term $\Div\rvect$ in it, merely belonging to $C(\R;L^2(\Omega))$, has some trace sense on $\R\times\Gamma_1$. This type of difficulty also arises when dealing with the somehow related wave equation with hyperbolic dynamical boundary conditions, studied by the author in a series of paper, see \cite{Dresda1,Dresda2,Dresda3}, and the natural solution of it consists in combining equations $(\Lagr)_1$ and $(\Lagr)_3$ in a single distributional identity, as done in Definition~\ref{Definition 4.1}.

 The following result, which will be also useful in the sequel, shows that Definition~\ref{Definition 4.1}--iii) is the closest possible approximation of the notion of distributional solutions of $(\Lagr)_1$ and $(\Lagr)_3$.
 \begin{prop}\label{Proposition 4.2} Let $(\rvect,v)\in X^1_\Lagr$ be such that $\Div\rvect\in L^1_\loc(\R; H^1(\Omega))$. Then
 $(\rvect,v)$ is a weak solution of $(\Lagr)$ if and only if it satisfies \eqref{4.50} and the further distributional identity
 \begin{equation}\label{4.51}
 \int_{-\infty}^\infty\int_{\Gamma_1}\mu v_t\psi_t-\sigma(\nabla_\Gamma v,\nabla_\Gamma \overline{\psi})_\Gamma-\delta v_t\psi
 -\kappa v\psi+B\Div\rvect\psi=0
 \end{equation}
 for all $\psi\in C^r_c(\R\times\Gamma_1)$.
 \end{prop}
 \begin{proof} Let $(\rvect,v)$ be as in the statement. We first claim that \eqref{4.50} holds if and only if
 \begin{equation}\label{4.52}
 \rvect_t\in W^{1,1}_\loc (\R;H^0_{\curl 0}(\Omega))\quad\text{and}\quad \rho_0\rvect_{tt}=B\nabla\Div \rvect
 \quad\text{in $L^1_\loc (\R;H^0_{\curl 0}(\Omega))$.}
 \end{equation}
 Indeed, if \eqref{4.52} holds one gets \eqref{4.50} simply  multiplying \eqref{4.52} by $\phivect$ and integrating \eqref{4.52} by parts in space and time. Conversely, if \eqref{4.50} holds, by integrating by parts in $\Omega$ one gets
 $\int_{\R\times\Omega} \rho_0\rvect_t\cdot \,\phivect_t+B\nabla\Div \rvect\cdot\phivect=0$ for all $\phivect\in\mathcal{D}(\R\times\Omega)$. Consequently, taking test functions $\phivect(t,x)=\varphi(t)\phivect_0(x)$ with $\varphi\in\mathcal{D}(\R)$ and $\phivect_0\in \mathcal{D}(\Omega)^3$ we get
   $\int_\Omega\left[\int_{-\infty}^\infty \rho_0\rvect_t\varphi'+B\nabla\Div \rvect \varphi\right]\cdot\phivect_0=0$.
 Since the last equation trivially extends by density to $\phivect_0\in H^1_0(\Omega)^3$, we have that $\rvect_t\in W^{1,1}_\loc(\R;H^{-1}(\Omega)^3)$ and equation \eqref{4.52} holds true in $L^1_\loc(\R;H^{-1}(\Omega)^3)$. Since $\nabla\Div \rvect\in L^1_\loc(\R; H^0_{\curl 0}(\Omega))$ and $H^0_{\curl 0}(\Omega)\hookrightarrow L^2(\Omega)^3\simeq [L^2(\Omega)']^3\hookrightarrow H^{-1}(\Omega)^3$ we then get that $\rvect_t\in W^{1,1}_\loc (\R;H^0_{\curl 0}(\Omega))$  and equation \eqref{4.52} hold in the space $L^1_\loc(\R;H^0_{\curl 0}(\Omega))$, so proving \eqref{4.52} and our claim.

 Now let $(\rvect,v)$ be a weak solution of $(\Lagr)$. Taking in \eqref{4.49} test functions $\phivect\in \mathcal{D}(\R\times\Omega)^3$ one obtains \eqref{4.50}. By the previous claim then \eqref{4.52} holds and then, integrating by parts in space and time we can rewrite \eqref{4.49} as \eqref{4.51} but with test functions $\psi$ in a different class, that is $\psi=-\phivect\cdot\boldsymbol{\nu}$, where $\phivect\in C^{r-1}_c(\R\times\R^3)^3$ is such that $\phivect\cdot\boldsymbol{\nu}=0$ on $\R\times\Gamma_0$. Trivially such $\psi$ restrict to $\psi\in C^{r-1}(\R\times\Gamma_1)$, and we actually claim that \eqref{4.51} hold true {\em for all $\psi\in C^{r-1}(\R\times\Gamma_1)$}, and then by standard density properties  for all $\psi\in C^r(\R\times\Gamma_1)$, proving the direct implication in the assertion. Our claim follows since for any $\psi\in C^{r-1}_c(\R\times\Gamma_1)$ its trivial extension $\widetilde{\psi}$ on the whole of $\R\times\Gamma$ belongs, as $\overline{\Gamma_0}\cap\overline{\Gamma_1}=\emptyset$, to $C^{r-1}_c(\R\times\Gamma)$. Defining then $\phivect\in C^{r-1}_c(\R\times\Gamma)^3$ by $\phivect=-\widetilde{\psi}\boldsymbol{\nu}$, we have $\psi=-\phivect\cdot\boldsymbol{\nu}$ on $\R\times\Gamma_1$ and  $\phivect\cdot\boldsymbol{\nu}=0$ on $\R\times\Gamma_0$. Using the compactness of $\Gamma$ and \cite[Definition 1.2.1.1, Chapter 1, p. 5]{grisvard} it is then straightforward to extend $\phivect$ to $\phivect\in C^{r-1}_c(\R\times\R^3)^3$ by using local equations, cut--off arguments and partitions of the unity, so proving our claim.

 Conversely, now let $(\rvect,v)$ satisfies \eqref{4.50} and \eqref{4.51}. We remark that, by density, \eqref{4.51} holds for all
 $\psi\in C^{r-1}(\R\times\Gamma_1)$. By our claim we can rewrite \eqref{4.50} as \eqref{4.52}. We multiply it by a test function $\phivect$ as in Definition~\ref{Definition 4.1}--iii) and we integrate in space and time the resulting equation. In this way we get
 $\int_{-\infty}^\infty\int_\Omega \rho_0\rvect_t\cdot\phivect_t-B\Div\rvect\Div\phivect=B\int_{-\infty}^\infty\int_{\Gamma_1}\Div\rvect \psi$
 which  with \eqref{4.51} gives \eqref{4.49}.
 \end{proof}
 Proposition~\ref{Proposition 4.2} allows to point out the trivial relations among the three types of solutions of $(\Lagr)$ introduced in Definition~\ref{Definition 4.1} .
\begin{lem}\label{Lemma 4.6}Let $(\rvect,v)\in X^1_\Lagr$ be a solution of $(\Lagr)$ according to Definition~\ref{Definition 4.1}. Then strong $\Rightarrow$ generalized $\Rightarrow$ weak and, if $(\rvect,v)\in X^2_\Lagr$, weak $\Rightarrow$ strong.
\end{lem}
\begin{proof} Strong solutions are also generalized ones and, by Proposition~\ref{Proposition 4.2}, also weak. Since \eqref{4.49} is stable with respect to the convergence in $X^1_\Lagr$, we then get that generalized $\Rightarrow$ weak. To prove the  final conclusion
let $(\rvect,v)\in X^2_\Lagr$ be a weak solution. By  Proposition~\ref{Proposition 4.2} it also satisfies  $(\Lagr)_1$ and $(\Lagr)_3$  in a distributional sense. Being regular enough we can integrate by parts to show that it  also satisfies them a.e., so $(\rvect,v)$ is a strong solution.
\end{proof}
\begin{rem}\label{Remark 4.1}\label{Remark 4.1}
By the last result the same conclusions in Remark~\ref{Remark 3.3} apply as well to problem $(\Lagr)$, so  solutions $(\rvect,v)\in X^2_\Lagr$ are equivalently strong, generalized or weak.
\end{rem}
In view of Theorem~\ref{Theorem 4.1} and Lemmas~\ref{Lemma 4.5} and \ref{Lemma 4.6}, which provide the existence of a generalized and thus weak solution of $(\Lagr_0)$, a natural question arising is the uniqueness of weak solutions. To positively answer the question it is useful to characterize weak solutions of $(\Lagr)$ as solutions of an abstract second order ODE.

Specifically for this purpose we shall identify in the sequel $L^2(\Omega)^3$ and $L^2(\Gamma_1)$ with their duals $[L^2(\Omega)^3]'$ and $[L^2(\Gamma_1)]'$, coherently with the distributional identification. Moreover we shall also identify
$L^2(\Omega)^3$ and $L^2(\Gamma_1)$ with their isometric copies $L^2(\Omega)^3\times\{0\}$ and $\{0\}\times L^2(\Gamma_1)$ contained in $L^2_\Lagr=L^2(\Omega)^3\times L^2(\Gamma_1)$. As a consequence we shall identify $L^2_\Lagr$ with its dual ${L^2_\Lagr}'$ according to the identity
$\langle \mathbf{u},\mathbf{w}\rangle_{L^2_\Lagr}=(\mathbf{u},\overline{\mathbf{w}})_{L^2_\Lagr}$  for all $\mathbf{u,w}\in L^2_\Lagr$.
By Lemma~\ref{Lemma 4.4} we can introduce the chain of dense embeddings, or Gel'fand triple,
$V_\Lagr\hookrightarrow L^2_\Lagr\simeq {L^2_\Lagr}'\hookrightarrow V_\Lagr'$,
in view of which the last identity particularizes to
\begin{equation}\label{4.56}
 \langle \mathbf{u},\mathbf{w}\rangle_{V_\Lagr}=(\mathbf{u},\overline{\mathbf{w}})_{L^2_\Lagr}\qquad\text{for all $\mathbf{u}\in L^2_\Lagr$ and $\mathbf{w}\in V_\Lagr$.}
\end{equation}
We also introduce the multiplication operator $R\in\mathcal{L}(L^2_\Lagr)$, the projection operator $\Pi_{\Gamma_1}\in\mathcal{L}(L^2_\Lagr; L^2(\Gamma_1))$ and
$M\in\mathcal{L}(V_\Lagr;V_\Lagr')$ given by
\begin{equation}\label{4.57}
R(\rvect,v)=(\rho_0 \rvect,\mu v),\qquad   \Pi_{\Gamma_1}(\rvect,v)=(0,v)\qquad\text{for all $(\rvect,v)\in L^2_\Lagr$,}
\end{equation}
and, for all $(\rvect,v), (\phivect,\psi)\in V_\Lagr$, by
\begin{equation}\label{4.58}
  \langle M(\rvect,v),(\phivect,\psi)\rangle_{V_\Lagr}=-B\int_\Omega \Div\rvect\Div\phivect-\int_{\Gamma_1}\sigma (\nabla_\Gamma v,\nabla_\Gamma \overline{\psi})_\Gamma-\int_{\Gamma_1}kv\psi.
\end{equation}
We then get the following result.
\begin{prop}\label{Proposition 4.3}
For any $\mathbf{u}=(\rvect,v)\in X^1_\Lagr$ the following properties are equivalent:
\renewcommand{\labelenumi}{{\roman{enumi})}}
\begin{enumerate}
\item $(\rvect,v)$ is a weak solution of $(\Lagr)$;
\item $R\mathbf{u}'\in C^1(\R;V_\Lagr')$ and
\begin{equation}\label{4.59}
(R\mathbf{u}')'  +\delta\Pi_{\Gamma_1}\mathbf{u}'-M\mathbf{u}=0\qquad\text{in $C(\R;V_\Lagr')$;}
\end{equation}
\item $(\rvect,v)$ satisfies the generalized distributional identity
\begin{multline}\label{4.60}
 \int_s^t\left\{\int_\Omega \rho_0\rvect_t\cdot \phivect_t-B\Div \rvect\Div \phivect
+ \int_{\Gamma_1}\mu v_t\psi_t-\sigma(\nabla_\Gamma v,\nabla_\Gamma \overline{\psi})_\Gamma\right.\\
\left.-\int_{\Gamma_1}\delta v_t\psi +\kappa v\psi\right\}=\left[\rho_0\int_\Omega\rvect_t\cdot\phivect+\int_{\Gamma_1}\mu v_t\psi\right]_s^t
 \end{multline}
for all $s,t\in\R$ and $(\phivect,\psi)\in C(\R;V_\Lagr)\cap C^1(\R; L^2_\Lagr)$.
\end{enumerate}
\end{prop}
\begin{proof} We shall prove the implications i)$\Rightarrow$ ii)$\Rightarrow$ iii) $\Rightarrow$ i), starting from the trivial one.
Indeed \eqref{4.60} reduces to \eqref{4.49} when taking test functions as in Definition~\ref{Definition 4.1}--iii), hence iii) $\Rightarrow$ i).

Moreover ii) $\Rightarrow$ iii). Indeed,  by ii), we have $R\mathbf{u}'\in C(\R;V_\Lagr)\cap C^1(\R;V_\Lagr')$. By a trivial extension of the Leibnitz rule then, for any $\boldsymbol{\psi}=(\phivect,\psi)\in C(\R;V_\Lagr)\cap C^1(\R;L^2_\Lagr)$ we have
$\frac d{dt}(R\mathbf{u}',\overline{\boldsymbol{\Psi}})_{L^2_\Lagr}=\langle (R\mathbf{u}')',\boldsymbol{\psi}\rangle_{V_\Lagr}
+(R\mathbf{u'},\overline{\boldsymbol{\psi}}')_{L^2_\Lagr}\in C(\R)$,
so $(R\mathbf{u}',\overline{\boldsymbol{\Psi}})_{L^2_\Lagr}\in C^1(\R)$. Consequently, integrating from $s$ to $t$ and using \eqref{4.59}, we get
$\left[(R\mathbf{u}',\overline{\boldsymbol{\Psi}})_{L^2_\Lagr}\right]_s^t=\int_s^t -(\delta\Pi_{\Gamma_1}\mathbf{u}',\overline{\boldsymbol{\psi}})_{L^2_\Lagr}+\langle M\mathbf{u},\boldsymbol{\psi}\rangle_{V_\Lagr}
+(R\mathbf{u}', \overline{\boldsymbol{\psi}}')_{L^2_\Lagr}
$
which, using \eqref{4.57}--\eqref{4.58}, is \eqref{4.60}.

To complete the proof we then just have to prove that i)$\Rightarrow$ ii), so let $(\rvect,v)$ be a weak solution of $(\Lagr)$. Taking in \eqref{4.49} test functions $\phivect(t,x)=\varphi(t)\mathbf{w}(x)$, with $\varphi\in \mathcal{D}(\R)$ and $\mathbf{w}\in C^{r-1}_c(\R^3)^3$ such that $\mathbf{w}\cdot \boldsymbol{\nu}=0$ on $\Gamma_0$, and denoting $\chi=-\mathbf{w}\cdot \boldsymbol{\nu}_{|\Gamma_1}$, we get
\begin{multline}\label{4.61}
 \rho_0\int_{-\infty}^\infty\varphi'\int_\Omega \rvect_t\cdot \mathbf{w}-B
\int_{-\infty}^\infty \varphi \int_\Omega \Div \rvect\Div \mathbf{w}+ \int_{-\infty}^\infty\varphi'\int_{\Gamma_1} \mu v_t\chi\\
-\int_{-\infty}^\infty\varphi\int_{\Gamma_1}\sigma(\nabla_\Gamma v,\nabla_\Gamma \overline{\chi})_\Gamma+(\delta v_t+\kappa v)\chi=0,
 \end{multline}
and we remark that $(\mathbf{w},\chi)\in V_\Lagr$.

We claim that \eqref{4.61} actually holds {\em for all $(\mathbf{w},\chi)\in V_\Lagr$}. Fixing such a $(\mathbf{w},\chi)$ we remark that, by \eqref{4.23} and the Divergence Theorem, the couple
$\left(-\frac B{\rho_0}\Div\mathbf{w},\chi\right)$ satisfies \eqref{4.11}, so by Lemma~\ref{Lemma 4.2} there is a unique $\mathbf{w}_1\in H^1_{\curl 0}(\Omega)$ such that $\Div \mathbf{w}_1=\Div\mathbf{w}$ in $\Omega$ and $(\mathbf{w}_1,\chi)\in \mathbb{H}^1_\Lagr$. We now set $\mathbf{w}_2=\mathbf{w}-\mathbf{w}_1$, so
\begin{equation}\label{4.62}
 \mathbf{w}= \mathbf{w}_1+\mathbf{w}_2,\quad \Div\mathbf{w}_2=0\quad\text{in $\Omega$, and}\quad \mathbf{w}_2\cdot\boldsymbol{\nu}=0\quad\text{on $\Gamma$.}
\end{equation}
By linearity we shall prove \eqref{4.61} for $\mathbf{w}=\mathbf{w}_i$, $i=1,2$. When $\mathbf{w}=\mathbf{w}_2$, by \eqref{4.62} we can rewrite \eqref{4.61} as
\begin{equation}\label{4.63}
 \rho_0\int_{-\infty}^\infty\varphi'\int_\Omega \rvect_t\cdot \mathbf{w}_2=0,
 \end{equation}
Now, since $(\rvect,v)\in X^1_\Lagr$, by \eqref{1.7}--\eqref{1.11} and Lemma~\ref{Lemma 4.1} for any $t\in\R$ there is
$\varphi_2(t)\in\dot{H}^2(\Omega)$ such that $\rvect(t)=-\nabla \varphi_2(t)$. Consequently, by \eqref{4.62} and the Divergence Theorem, we have
$$\int_\Omega \rvect(t)\cdot\mathbf{w}_2=-\int_\Omega \nabla\varphi_2(t)\cdot\mathbf{w}_2=\int_\Omega \varphi_2(t)\Div \mathbf{w}_2-\int_\Gamma \varphi_2(t)\mathbf{w}_2\cdot \boldsymbol{\nu}=0\quad\text{for all $t\in\R$,}
$$
from which \eqref{4.63} trivially follows.

When $\mathbf{w}=\mathbf{w}_1$, since $(\mathbf{w}_1,\chi)\in \mathbb{H}^1_\Lagr$, using Lemma~\ref{Lemma 4.3} we can suppose by density that $(\mathbf{w}_1,\chi)\in \mathbb{Y}$, i.e. that $\mathbf{w}_1$ extends to $\mathbf{w}_1\in C^{r-1}(\R^3)^3$. Moreover, since $\Omega$ is bounded, a standard truncation argument allows also to suppose that $\mathbf{w}_1\in C^{r-1}_c(\R^3)^3$. For such a $\mathbf{w}_1$ we can take, as a test function in \eqref{4.49}, $\phivect(t,x)=\varphi(t)\mathbf{w}_1(x)$, from which \eqref{4.61} trivially follows when $\mathbf{w}=\mathbf{w}_1$, completing the proof of our claim.

Recalling \eqref{4.57} and \eqref{4.58} we can then rewrite \eqref{4.61}, for test functions $\phivect_0=(\mathbf{w},\chi)\in V_\Lagr$, as
$\int_{-\infty}^\infty \varphi'(R\mathbf{u}',\overline{\phivect_0})_{L^2_\Lagr}+\varphi\left[\langle M\mathbf{u},\phivect_0\rangle_{V_\Lagr}-(\delta\Phi_{\Gamma_1}\mathbf{u}',\overline{\phivect_0})_{L^2_\Lagr}\right]=0
$
or, using \eqref{4.56} and Bochner integrals,
$\int_{-\infty}^\infty \left[\varphi'R\mathbf{u}'+\varphi\left(M\mathbf{u}-\delta\Phi_{\Gamma_1}\mathbf{u}'\right)\right]=0$
 in $V_\Lagr'$,  for all $\varphi\in\mathcal{D}(\R)$.
Consequently $\mathbf{u}$ solves equation \eqref{4.59} in $\mathcal{D}'(\R)$. Since $M\mathbf{u}-\delta\Phi_{\Gamma_1}\mathbf{u}'\in C(\R;V_\Lagr')$ we consequently get that $R\mathbf{u}'\in C^1(\R;V_\Lagr')$ and \eqref{4.59} holds.
\end{proof}
We can give the promised uniqueness result.
\begin{thm}[\bf Uniqueness]  Weak solutions of $(\Lagr_0)$ are unique.\label{Theorem 4.3}\end{thm}
\begin{proof}By linearity it is enough to prove that $\rvect_0=\rvect_1=0$ and $v_0=v_1=0$ yield $(\rvect,v)\equiv 0$ in $\R$. Moreover, since the couple $(\widetilde{\rvect},\widetilde{v})$ defined by $\widetilde{\rvect}(t)=\rvect(-t)$ and $\widetilde{v}(t)=v(-t)$ for $t\in\R$ is still a weak solution of $(\Lagr_0)$ provided $\delta$ is replaced by$-\delta$, we just have to prove that $(\rvect,v)\equiv 0$ in $[0,\infty)$.

Adapting the arguments in  the proof of \cite[Chapter 4,~Lemma~4.2.5]{mugnvit}, we fix $t>0$ and  test functions $(\phivect,\psi)\in C^1(\R;V_{\Lagr})$, depending on $t$, given by
\begin{equation}\label{4.64}
\phivect(s)=
\begin{cases}
\int_s^t \overline{\rvect}(\tau)\,d\tau\quad &\text{if $s\le t$},\\
\overline{\rvect}(t)(t-s) &\text{if $s\ge t$},
\end{cases}
\qquad
\psi(s)=
\begin{cases}
\int_s^t \overline{v}(\tau)\,d\tau\quad &\text{if $s\le t$},\\
\overline{v}(t)(t-s) &\text{if $s\ge t$},
\end{cases}
\end{equation}
which trivially verifies
\begin{equation}\label{4.65}
\phivect(t)=0,\quad \psi(t)=0,\quad\text{and}\quad
\phivect_t=-\overline{\rvect},\quad  \psi_t=-\overline{v}\quad\text{in $[0,t]$.}
\end{equation}
By Proposition~\ref{Proposition 4.3} we can take $(\phivect,\psi)$ as a test function in \eqref{4.60}, writing it for $s=0$. Since $\rvect_1=0$ and $v_0=v_1=0$, using also \eqref{4.65} we consequently get
\begin{multline*}
  \int_0^t\left\{\int_\Omega [\rho_0\rvect_t\cdot\overline{\rvect}-B\Div \overline{\phivect}_t\Div\phivect]+\int_{\Gamma_1} \left[\mu  v_t\overline{v}-\sigma (\nabla_\Gamma (\overline{\psi}_t),\nabla_\Gamma\overline{\psi})_\Gamma-\kappa \psi\overline{\psi}_t]\right]\right\}\\
  =-\int_0^t\int_{\Gamma_1}\delta v_t\psi=\int_0^t \int_{\Gamma_1}\delta v\psi_t=-\int_0^t\int_{\Gamma_1} \delta |v|^2,
\end{multline*}
were we also integrated by parts in time the right--hand side. Taking the real part we then get
$$\int_0^t\frac 12\frac d{dt}\left[\rho_0\|\rvect\|_2^2-B\|\Div\phivect\|_2^2+\int_{\Gamma_1}\mu |v|^2- \int_{\Gamma_1}\sigma |\nabla_\Gamma \psi|_\Gamma^2-\int_{\Gamma_1}\kappa |\psi|^2\right] =
-\int_0^t\int_{\Gamma_1}\delta |v|^2$$
which, since $\rvect_0=\phivect(t)=0$ and $v_0=\psi(t)=0$, can be rewritten as
\begin{equation}\label{4.65BIS}
\tfrac {\rho_0}2\|\rvect(t)\|_2^2+\tfrac B2\|\Div\phivect(0)\|_2^2+\int_{\Gamma_1}\mu |v|^2+ \int_{\Gamma_1}\sigma |\nabla_\Gamma \psi(0)|_\Gamma^2=-\int_{\Gamma_1}\kappa |\psi(0)|^2 -\int_0^t\int_{\Gamma_1}\delta |v|^2.
\end{equation}
We now  estimate the right--hand side of \eqref{4.65BIS} by using  assumption (A), \eqref{4.64} and H\"{o}lder inequality to get
$\|v(t)\|_{2,\Gamma_1}^2\le \tfrac 1{\mu_0}(\|\delta\|_{\infty,\Gamma_1}+\|\kappa\|_{\infty,\Gamma_1} t) \int_0^t\|v\|_{2,\Gamma_1}^2$
for all $t\ge 0$.
Denoting $\Upsilon(t)=\int_0^t\|v\|_{2,\Gamma_1}^2$ for $t\in [0,\infty)$ we can rewrite the last estimate as
$\Upsilon'(t)\le \tfrac 1{\mu_0}(\|\delta\|_{\infty,\Gamma_1}+\|\kappa\|_{\infty,\Gamma_1} t)\Upsilon(t)$ for all $t\ge 0$,
so a standard integration yields, since $\Upsilon(0)=0$, that $\Upsilon\equiv 0$ in $[0,\infty)$, from which $v\equiv 0$ in $[0,\infty)$.

Fixing $t>0$ again, and recalling the test functions $(\phivect,\psi)$ in \eqref{4.64}, we then get $\psi(0)=0$.
Plugging it in the consequent identity  \eqref{4.65BIS}, also recalling that $v\equiv 0$, by assumption (A) we then get $\rvect(t)=0$. Being $t>0$ arbitrary we then have $(\rvect,v)\equiv 0$ in $[0,\infty)$, concluding the proof.
\end{proof}
\subsection{Main results for problems $(\Lagr)$ and $(\Lagr_0)$}\label{Section 4.5} We can finally prove Theorems~\ref{Theorem 1.2} and \ref{Theorem 1.4} and deal with  optimal regularity issues.
\begin{proof}[\bf Proof of Theorem~\ref{Theorem 1.2}]
By combining Theorem~\ref{Theorem 4.1} and Lemma~\ref{Lemma 4.5} we get that for all $U_{0\Lagr}\in \mathcal{H}^1_\Lagr$ problem
$(\Lagr_0)$ has a unique generalized solution $(\rvect,v)\in X^1_\Lagr$, continuously depending on  $U_{0\Lagr}$ in the topologies of the respective spaces. By Lemma~\ref{Lemma 4.6} the couple $(\rvect,v)$ is also a weak solution of $(\Lagr_0)$ which, by Theorem~\ref{Theorem 4.3}, is unique among weak solutions. By Theorems~\ref{Theorem 4.1} -- \ref{Theorem 4.2} and Lemma~\ref{Lemma 4.5} the couple  $(\rvect,v)$ is a strong solution if and only if $U_{0\Lagr}\in D(A_\Lagr)$, which by \eqref{4.41} exactly reads in the statement, the continuous dependence following by Theorem~\ref{Theorem 4.2}.

To complete the proof we just have to get the energy identity \eqref{1.12}, which can be obtained either by re--deriving it for strong solutions and hence using a density argument or deducing it from the energy identity \eqref{1.5} for problem $(\Pot_0)$. Choosing the second alternative, we remark that by Lemma~\ref{Lemma 4.5} and Theorem~\ref{Theorem 4.1}, $\rvect$ and $v$ are the first two components of $\Phi_\PcL\dot{U}_\Pot$, where $\dot{U}_\Pot$ is the solution of \eqref{3.32} corresponding to $\dot{U}_{0\Pot}=F_\LPc U_{0\Lagr}$ and, by \eqref{3.34}, we have $\dot{U}_\Pot=(\pi_0 u,v,u_t,v_t)$ for the solution $(u,v)$ of $(\Pot_0)$ corresponding to
$U_{0\Pot}\in Q^{-1}\dot{U}_{0\Pot}$. Using \eqref{4.29} and \eqref{4.32} we have $-\nabla u=\rvect_t$ and $-B\Div\rvect=\rho_0u_t$. Plugging them into \eqref{1.5} we then get \eqref{1.12} and complete the proof.
\end{proof}
The uniqueness of weak solution of $(\Lagr_0)$ in Theorem~\ref{Theorem 4.3} together with the existence of a generalized and hence weak solution of it shows that solutions of $(\Lagr)$ are {\em equivalently} weak or generalized. Combining this remark with Remark~\ref{Remark 4.1} we thus obtain that, as for problem $(\Pot)$,
{\em all type of solutions in Definition~\ref{Definition 4.1} coincide, strong solutions being defined only in the class $X^2_\Lagr$}. Consequently in the sequel we  shall only deal  with weak solutions of $(\Lagr)$, i.e. with elements of the spaces $\mathcal{S}^n_\Lagr$ defined in \eqref{1.18}. Recalling the spaces $\mathcal{T}_\Lagr^n$ defined in \eqref{4.47} it is useful to point out the following result, which is the equivalent of Lemma~\ref{Lemma 3.3} for problem $(\Lagr)$.

\begin{lem}\label{Lemma 4.7}
The operator $\mathcal{I}'\in\mathcal{L}(X^1_\Lagr;Y^1_\Lagr)$ defined by
$\mathcal{I}'(\rvect,v)=(\rvect,v,\rvect_t,v_t)$restricts for each $n\in\widetilde{\N}$, $n\le r$, to a bijective isomorphism $\mathcal{I}^n_\Lagr\in\mathcal{L}(\mathcal{S}^n_\Lagr;\mathcal{T}^n_\Lagr)$ with inverse $\left(\mathcal{I}^n_\Lagr\right)^{-1}\in\mathcal{L}(\mathcal{T}^n_\Lagr;\mathcal{S}^n_\Lagr)$ given by
$\left(\mathcal{I}^n_\Lagr\right)^{-1}(\rvect,v,\svect,z)=(\rvect,v)$.
\end{lem}
\begin{proof} It  follows by \eqref{1.7} -- \eqref{1.11}, Lemma~\ref{Lemma 4.5} and the remarked equivalence.
\end{proof}
\begin{proof}[\bf Proof of Theorem~\ref{Theorem 1.4}]
By Proposition~\ref{Proposition 3.2}, Corollary~\ref{Corollary 3.4} and Corollary~\ref{Corollary 4.1}, for any
$(u,v)\in \mathcal{S}^1_\Potc$ we have
$[\Phi_\PcL\mathcal{J}^1_\Potc (u,v)](t)= (\rvect(t),v(t),-\nabla u(t),v_t(t))$  for all $t\in\R$,
where, by using Lemma~\ref{Lemma 4.2}, $\rvect(t)$ is defined by \eqref{1.21}, and \eqref{1.22} holds, so by Lemma~\ref{Lemma 4.7} we get the operator $\Psi_\PcL$ as in the statement, that is
\begin{equation}\label{4.72}
  \Psi_\PcL=(\mathcal{I}_\Lagr^1)^{-1}\cdot \Phi_\PcL\cdot \mathcal{J}^1_\Potc\in\mathcal{L}(\mathcal{S}^1_\Potc; \mathcal{S}^1_\Lagr).
\end{equation}
By Proposition~\ref{Proposition 3.2}, Corollary~\ref{Corollary 3.4} and Corollary~\ref{Corollary 4.1} we also get  that $\text{Ker } \Psi_\PcL= \C_{X_\Pot}$, and the bijective isomorphisms
\begin{equation}\label{4.73}
  \dot{\Psi}_\PcL=(\mathcal{I}_\Lagr^1)^{-1}\!\!\!\cdot \Phi_\PcL\!\!\!\cdot \dot{\mathcal{J}}^1_\Potc\!\in\mathcal{L}(\dot{\mathcal{S}}^1_\Potc; \mathcal{S}^1_\Lagr),\quad
  \dot{\Psi}_\LPc=(\dot{\mathcal{J}}^1_\Potc)^{-1}\!\!\cdot \Phi_\LPc\!\!\!\cdot\mathcal{I}_\Lagr^1
  \!\in\mathcal{L}(\mathcal{S}^1_\Lagr;\dot{\mathcal{S}}^1_\Potc ),
\end{equation}
trivially being one the inverse of the other. Moreover, by \eqref{4.32} and \eqref{4.29}, for any $(\rvect,v)\in\mathcal{S}^1_\Lagr$ we have
$$ [\Phi_\LPc\mathcal{I}^1_\Lagr (\rvect,v)](t)= F_\LPc (\rvect(t),v(t),\rvect_t(t),v_t(t))=\left(\dot u(t),v(t),-\tfrac B{\rho_0}\Div \rvect(t),v_t(t)\right)$$
for all $t\in\R$, where $\dot u(t)\in\dot{H}^1(\Omega)$ is the unique solution of $-\nabla \dot u(t)=\rvect_t(t)$. Consequently, since $\mathcal{J}^1_\Potc$ is surjective, there is $(u,v)\in \mathcal{S}^1_\Potc$ such that
\begin{equation}\label{4.73BIS}
\left(\dot u(t), v(t), -\tfrac B{\rho_0}\Div \rvect(t),v_t(t\right)=(\pi_0 u(t), v(t), u_t(t), v_t(t))\quad\text{for all $t\in\R$,}
\end{equation}
that is
\begin{equation}\label{4.74}
\pi_0 u(t)=  \dot u(t)\quad\text{and}\quad -B\Div \rvect(t)=\rho_0 u_t(t)\qquad\text{for all $t\in\R$,}
\end{equation}
which yields \eqref{1.26} since, by \eqref{1.22}, we have $-\nabla \dot u(0)=\rvect_t(0)$. Moreover trivially \eqref{1.26} defines $u$ up to a space--time constant and then also completely defines $\dot{\Psi}_\LPc(\rvect,v)$.
Finally, to prove \eqref{1.27}, we remark that, as we just proved,  $(\rvect,v)=\Psi_\PcL(u,v)$  is equivalent to \eqref{4.73BIS}, that is  \eqref{4.74}. By uniqueness then it is equivalent to \eqref{4.74} at $t=0$, that is to $-\nabla u_0=\rvect_1$ and $-B\Div \rvect_0=\rho_0u_1$.
\end{proof}
We can also give an optimal regularity result for problem $(\Lagr_0)$ and show that the operators in Theorem~\ref{Theorem 1.4} have a good behavior with respect to regularity classes.
\begin{thm}[\bf Optimal regularity for $(\Lagr_0)$ and $(\Potc) \rightleftarrows (\Lagr)$ ]\label{Theorem 4.4}
For all data $U_{0\Lagr}=(\rvect_0, v_0,\rvect_1,v_1)\in\mathcal{H}^1_\Lagr$ and $n\in\widetilde{\N}$, $2\le n\le r$, the weak solution $(\rvect,v)$ given in Theorem~\ref{Theorem 1.2} belongs to $X_\Lagr^n$ if and only if $U_{0\Lagr}\in D^{n-1}_\Lagr$, and in this case
$(\rvect,v)$ continuously depends on the data $U_{0\Lagr}$ in the topologies of the respective spaces.
Moreover the operator $\Psi_\PcL$ in Theorem~\ref{Theorem 1.4} restricts to a surjective operator $\Psi^n_\PcL\in\mathcal{L}(\mathcal{S}^n_\Potc;\mathcal{S}^n_\Lagr)$ such that $\ker \Psi^n_\PcL=\C_{X_\Pot}$.
Consequently the bijective isomorphisms $\dot{\Psi}_\PcL$ and $\dot{\Psi}_\LPc$ in the same result respectively restrict to bijective isomorphisms
\begin{equation}\label{4.75}
  \dot{\Psi}_\PcL^n\in\mathcal{L}(\dot{\mathcal{S}}^n_\Potc;\mathcal{S}^n_\Lagr),\quad\text{and}\quad
  \dot{\Psi}_\LPc^n\in\mathcal{L}(\mathcal{S}^n_\Lagr;\dot{\mathcal{S}}^n_\Potc).
\end{equation}
\end{thm}
\begin{proof} The optimal regularity for problem $(\Lagr_0)$ immediately follows by combining Theorem~\ref{Theorem 4.2} and Lemma~\ref{Lemma 4.7}. Moreover the restriction properties of the operators $\Psi_\PcL$, $\dot{\Psi}_\PcL$ and $\dot{\Psi}_\LPc$,  by formulas \eqref{4.72} and \eqref{4.73}, are simply consequences of Corollaries~\ref{Corollary 3.4},  \ref{Corollary 4.1}  and Lemma~\ref{Lemma 4.7}.
\end{proof}
\section{The Eulerian models}\label{Section 5}
\subsection{Preliminaries}\label{Section 5.1}
To deal with problems $(\Eul)$ and $(\Eulc)$ it is convenient, beside the spaces already defined in \S~\ref{intro},  to introduce some additional functional spaces. We anticipate the compatibility conditions for problem $(\Eul_0)$ needed in the sequel, that is, for $n\in\N$, $2\le n\le r$ and data $U_{0\Eul}=(p_0,\vvect_0,v_0,v_1)\in \mathcal{H}^n_\Eul$,
\begin{equation}\label{5.1}\left\{
\begin{aligned}
&\vvect_0\cdot \boldsymbol{\nu}=-v_1,\,\,\qquad \text{on $\Gamma_1$,}\qquad \vvect_0\cdot \boldsymbol{\nu}=0,\quad \text{on $\Gamma_0$,} \\
&\partial_{\boldsymbol{\nu}}\Delta^i\Div\vvect_0 =0,\quad \text{on $\Gamma_0$ \qquad for $i=0,\ldots, \lfloor n/2\rfloor -2$, \,\,\,\,\qquad when $n\ge 4$,}\\
&\partial_{\boldsymbol{\nu}}\Delta^i p_0 =0,\quad \text{on $\Gamma_0$ \qquad for $i=0,\ldots, \lfloor (n-1)/2\rfloor -1$, \qquad when $n\ge 3$,}\\
&\mu\partial_{\boldsymbol{\nu}}p_0=\rho_0[\,\DivGamma (\sigma\nabla_\Gamma v_0)-\delta v_1-\kappa v_0-p_0\,],\quad \text{on $\Gamma_1$, \qquad \,\,\, when $n\ge 3$,}\\
&\tfrac{B\mu}{\rho_0}\partial_{\boldsymbol{\nu}}\Div\vvect_0\!\!=\!-\DivGamma (\sigma\nabla_\Gamma v_1)\!+\!\tfrac\delta{\rho_0}\partial_{\boldsymbol{\nu}}p_0\!+\!\kappa v_1\!+\!\tfrac B{\rho_0}\Delta p_0,\,\text{on $\Gamma_1$, when $n=4$,}\\
&\begin{split}
\tfrac{B\mu}{\rho_0} \partial_{\boldsymbol{\nu}}\Div\Delta^i\vvect_0\,\,\,=\,\,\,\DivGamma (\sigma\nabla_\Gamma
\Div\Delta^{i-1}\vvect_0)+\tfrac\delta{\rho_0} \partial_{\boldsymbol{\nu}}\Delta^i p_0-\kappa\partial_{\boldsymbol{\nu}}\Div\Delta^{i-1}\vvect_0\\
\qquad\- B\Div\Delta^i \vvect_0 \quad\text{on $\Gamma_1$,  \quad\,  for $i=1,\ldots, \lfloor n/2\rfloor -2$, \quad when $n\ge 6$,}
\end{split}\\
&\begin{split}
\tfrac{B\mu}{\rho_0}\partial_{\boldsymbol{\nu}}\Delta^i p_0\quad\,=\,\,\DivGamma (\sigma\nabla_\Gamma \partial_{\boldsymbol{\nu}}\Delta^{i-1}p_0)
+B\delta\partial_{\boldsymbol{\nu}}\Div\Delta^{i-1}\vvect_0-\kappa\partial_{\boldsymbol{\nu}}\Delta^{i-1}p_0\quad\\
\qquad- B\Div\Delta^i p_0 \quad \text{on $\Gamma_1$,  \quad for $i=1,\ldots, \lfloor (n-1)/2\rfloor -1$,  when $n\ge 5$.}
\end{split}\end{aligned}\right.
\end{equation}
By \eqref{2.9} and the Trace Theorem one gets that
$D^{n-1}_\Eul:=\{U_{0\Eul}\in \mathcal{H}^n_\Eul: \quad \text{\eqref{5.1} hold}\}$
and $D^{n-1}_\Eulc:=D^{n-1}_\Eul\cap\mathcal{H}^n_\Eulc$
are closed subspaces, respectively of $\mathcal{H}^n_\Eul$ and $\mathcal{H}^n_\Eulc$, and then Hilbert spaces.
When $r=\infty$ we also set the Fr\'{e}chet spaces
$$
 \mathcal{H}_\Eul^\infty = \left[C^\infty(\overline{\Omega})\times  C^\infty(\overline{\Omega})^3\times C^\infty(\Gamma_1)\times C^\infty(\Gamma_1)\right]\cap\mathcal{H}^1_\Eul,\quad \mathcal{H}_\Eulc^\infty=\mathcal{H}_\Eul^\infty\cap\mathcal{H}^1_\Eulc
$$
and their closed subspaces
$$ D^\infty_\Eul=\{U_{0\Eul}\in \mathcal{H}^n_\Eul: \quad \text{\eqref{5.1} hold for all $n\in\N$}\},
 \qquad D^\infty_\Eulc=D^\infty_\Eul\cap\mathcal{H}^n_\Eulc.
$$
Moreover, for $n\in\widetilde{N}$, $n\le r$ we set the Fr\'{e}chet spaces
\begin{equation}\label{5.5}
  Y^n_\Eul=\bigcap_{i=0}^{n-1} C^i(\R;\mathcal{H}^{n-i}_\Eul)\quad\text{and}\quad
  Y^n_\Eulc=\bigcap_{i=0}^{n-1} C^i(\R;\mathcal{H}^{n-i}_\Eulc).
\end{equation}
In this section we shall also use the Hilbert space
\begin{equation}\label{5.6}
 H(\Div,\Omega)=\{\vvect\in L^2(\Omega)^3: \Div \vvect\in L^2(\Omega)\},
\end{equation}
where $\Div\vvect$ is taken in the distributional sense, introduced in \cite[Chapter IX, \S 2, p. 203]{dautraylionsvol3}, equipped with the  norm $\|\cdot\|_{H(\Div,\Omega)}=\left(\|\cdot\|_2^2+\|\Div(\cdot)\|_2^2\right)^{1/2}$.
We recall, see \cite[Chapter IX, \S 2, Theorem 1, p. 204 and (1.1), p. 206]{dautraylionsvol3}, that the map $\vvect\mapsto\vvect\cdot\boldsymbol{\nu}$ defined from $C^\infty(\overline{\Omega})^3$ to $C^1(\Gamma)$ extends by density as a bounded linear operator from $H(\Div,\Omega)$ to $H^{-1/2}(\Gamma)$, for simplicity denoted by the same symbol, and that
\begin{equation}\label{5.8}
\int_\Omega \vvect\cdot\nabla \varphi+\int_\Omega \Div\vvect\varphi=\langle \vvect\cdot\boldsymbol{\nu},\varphi_{|\Gamma}\rangle_{H^{1/2}(\Gamma)}\quad\text{for all $\vvect\in H(\Div,\Omega),\, \varphi\in H^1(\Omega)$.}
\end{equation}
\subsection{The main isomorphisms}\label{Section 5.2} To deal with problems $(\Eul)$ and $(\Eulc)$ we first introduce the main isomorphisms between the phase space of problem  $(\Eul)$ and the space $\dot{\mathcal{H}}^1_\Pot$. We then set, for all $(\dot u,v,w,z)\in \dot{\mathcal{H}}^1_\Pot$ and $(p,\vvect,v,z)\in\mathcal{H}^1_\Eul$,
\begin{equation}\label{5.9}
F_\PE  (\dot u,v,w,z)=(\rho_0,w,-\nabla\dot u, v,z),\qquad F_\EP(p,\vvect,v,z)=\left(\dot u, v, p/\rho_0,z\right)
\end{equation}
where $\dot u\in\dot{H}^1(\Omega)$ is the unique solution of $-\nabla \dot u=\vvect$ given by Lemma~\ref{Lemma 4.1}.

Trivially $F_\PE\in\mathcal{L}(\dot{\mathcal{H}}^1_\Pot;\mathcal{H}^1_\Eul)$ and  $F_\EP\in\mathcal{L}(\mathcal{H}^1_\Eul;\dot{\mathcal{H}}^1_\Pot)$, so they induce the operators
$\Phi_\PE \in\mathcal{L}(\dot{Y}^1_\Pot; Y^1_\Eul)$  and $\Phi_\EP \in\mathcal{L}(Y^1_\Eul; \dot{Y}^1_\Pot)$
given by
\begin{equation}\label{5.11}
 (\Phi_\PE \dot{U}_\Pot)(t)=F_\PE \dot{U}_\Pot (t),\quad  (\Phi_\EP U_\Eul)(t)=F_\EP {U}_\Eul (t)\quad \text{for all $t\in\R$.}
\end{equation}
The following result points out all properties of these operators needed in the sequel
\begin{lem}\label{Lemma 5.1}
The operator $F_\PE$ is a bijective isomorphism between $\dot{\mathcal{H}}^1_\Pot$ and $\mathcal{H}^1_\Eul$ having
$F_\EP$ as inverse. Consequently $\Phi_\PE$ is a bijective isomorphism between $\dot{Y}^1_\Pot$ and  $Y^1_\Eul$ with inverse
$\Phi_\EP$. They restrict to operators
\begin{equation}\label{5.12}
\begin{aligned}
F_\PcEc\in\mathcal{L}(\dot{\mathcal{H}}^1_\Potc;\mathcal{H}^1_\Eulc), \qquad & F_\EcPc\in \mathcal{L}(\mathcal{H}^1_\Eulc;\dot{\mathcal{H}}^1_\Potc),\\
\Phi_\PcEc \in\mathcal{L}(\dot{Y}^1_\Potc;Y^1_\Eulc),\qquad &
\Phi_\EcPc \in\mathcal{L}(Y^1_\Eulc; \dot{Y}^1_\Potc),
\end{aligned}
\end{equation}
and, for $n\in\N$, $2\le n\le r$, to operators
\begin{equation}\label{5.13}
\begin{aligned}
F_\PE\in\mathcal{L}(\dot{\mathcal{H}}^n_\Pot;\mathcal{H}^n_\Eul), \qquad & F_\EP\in \mathcal{L}(\mathcal{H}^n_\Eul;\dot{\mathcal{H}}^n_\Pot),\\
\Phi_\PE \in\mathcal{L}(\dot{Y}^n_\Pot;Y^n_\Eul),\qquad &
\Phi_\EP \in\mathcal{L}(Y^n_\Eul; \dot{Y}^n_\Pot),\\
F_\PcEc\in\mathcal{L}(\dot{\mathcal{H}}^n_\Potc;\mathcal{H}^n_\Eulc), \qquad & F_\EcPc\in \mathcal{L}(\mathcal{H}^n_\Eulc;\dot{\mathcal{H}}^n_\Potc),\\
\Phi_\PcEc \in\mathcal{L}(\dot{Y}^n_\Potc;Y^n_\Eulc),\qquad &
\Phi_\EcPc \in\mathcal{L}(Y^n_\Eulc; \dot{Y}^n_\Potc),
\end{aligned}
\end{equation}
enjoying the same properties. Moreover, or $n\in\N$, $2\le n\le r$, we have
\begin{equation}\label{5.14}
  F_\PE \dot{D}^{n-1}_\Pot=D^{n-1}_\Eul,\qquad\text{and}\quad F_\PcEc \dot{D}^{n-1}_\Potc=D^{n-1}_\Eulc,
\end{equation}
so $F_\PE$,  $F_\EP$, $F_\PcEc$ and $F_\EcPc$ further restrict to
bijective isomorphisms.
\end{lem}
\begin{proof} The first two sentences in the statements, including \eqref{5.12} and \eqref{5.13}, trivially follow from \eqref{1.3}, \eqref{1.6}, \eqref{1.13}, \eqref{3.15}, \eqref{5.5}, \eqref{5.9}, \eqref{5.11} and Lemma~\ref{Lemma 4.1}. Moreover the last restriction property trivially follows from \eqref{5.12} -- \eqref{5.14}. To complete the proof we then have to show that \eqref{5.14} holds and, in particular, that $F_\PE \dot{D}^{n-1}_\Pot=D^{n-1}_\Eul$, since the second identity in \eqref{5.14} follows from the first one  by \eqref{5.12}. To prove that $F_\PE \dot{D}^{n-1}_\Pot=D^{n-1}_\Eul$ we just have to prove that, for all $(u_0,v_0,u_1,v_1)\in \mathcal{H}^n_\Pot$ and $(p_0,\vvect_0,v_0,v_1)\in\mathcal{H}^1_\Eul$ such that
\begin{equation}\label{5.15}
 p_0=\rho_0 u_1\qquad\text{and}\quad -\nabla u_0=\vvect_0,
\end{equation}
the compatibility conditions \eqref{3.4} and \eqref{5.1} are equivalent. Using \eqref{5.15} one immediately gets that $\vvect_0\cdot\boldsymbol{\nu}=-v_1$  is equivalent to $\partial_{\boldsymbol{\nu}}u_0=v_1$ on $\Gamma_1$. It is then straightforward to get the required equivalence.
\end{proof}
\subsection{Abstract analysis of problems $(\Eul)$, $(\Eul_0)$, $(\Eulc)$ and $(\Eul_0^c)$}\label{Section 5.3}
To deal with problem $(\Eul_0)$ in a semigroup setting we  reduce it to a first order problem, that is to the problem
\begin{equation}\label{5.16}
\begin{cases}
p_t+B \Div \vvect=0\qquad &\text{in
$\R\times\Omega$,}\\
\rho_0\vvect_t+\nabla p=0\qquad &\text{in
$\R\times\Omega$,}\\
\curl\vvect=0\qquad &\text{in
$\R\times\Omega$,}\\
v_t =z\qquad
&\text{on
$\R\times \Gamma_1$,}\\
\mu z_t- \DivGamma (\sigma \nabla_\Gamma v)+\delta z+\kappa v+p =0\qquad
&\text{on
$\R\times \Gamma_1$,}\\
\vvect\cdot{\boldsymbol{\nu}} =0 \quad\text{on $\R\times \Gamma_0$,}\qquad
\vvect\cdot{\boldsymbol{\nu}} =-z\qquad
&\text{on
$\R\times \Gamma_1$,}\\
p(0,x)=p_0(x),\quad \vvect(0,x)=\vvect_0(x) &
 \text{in $\Omega$,}\\
v(0,x)=v_0(x),\quad z(0,x)=v_1(x) &
 \text{on $\Gamma_1$.}
\end{cases}
\end{equation}
More formally, working in phase space $\mathcal{H}^1_\Eul$ in which \eqref{5.16}$_3$  holds, we introduce the unbounded operator $A_\Lagr: D(A_\Eul)\subset \mathcal{H}^1_\Eul \to \mathcal{H}^1_\Eul$ given by
\begin{gather}\label{5.17}
D(A_\Eul)=D^1_\Eul=\{(p,\vvect,v,z)\in \mathcal{H}^2_\Eul: (\vvect,z)\in \mathbb{H}^1_\Lagr\}\\
\label{5.18}
A_\Eul\begin{pmatrix}p\\\vvect\\v\\z\end{pmatrix} =
\begin{pmatrix}B\Div\vvect\\\tfrac 1{\rho_0}\nabla p\\-z\\
\frac 1\mu\left[-\DivGamma(\sigma\nabla_\Gamma v)+\delta z+\kappa v+p_{|\Gamma_1}\right]
\end{pmatrix},
\end{gather}
together with the abstract equation and Cauchy problem
\begin{gather}\label{5.19}U_\Eul'+A_\Eul U_\Eul=0\qquad\text{in $\mathcal{H}^1_\Eul$,}
\\\label{5.20}
U_\Eul'+A_\Eul U_\Eul=0\qquad\text{in $\mathcal{H}^1_\Eul$}, \qquad U_\Eul(0)=U_{0\Eul}\in \mathcal{H}^1_\Eul.
\end{gather}
The following result shows that  problems \eqref{3.32} and \eqref{5.20} are essentially equivalent,  simultaneously restricting to the spaces $\dot{\mathcal{H}}^1_\Potc$ and $\mathcal{H}^1_\Eulc$.
\begin{thm}[\bf Well--posedness for \eqref{5.20} and its restriction to $\mathcal{H}^1_\Eulc$]\label{Theorem 5.1}
\phantom{A}
\renewcommand{\labelenumi}{{\Roman{enumi})}}
\begin{enumerate}
\item The operator $-A_\Eul$ is densely defined and it generates on $\mathcal{H}^1_\Eul$ the strongly continuous
group $\{\Eul^1(t),t\in\R\}$ given by
\begin{equation}\label{5.21}
  \Eul^1(t)=F_\PE\Pot^1(t) F_\EP\qquad\text{for all $t\in\R$,}
\end{equation}
and hence similar to the group $\{\Pot^1(t),t\in\R\}$ in Proposition~\ref{Proposition 3.1}.
Consequently, for any $U_{0\Eul}\in \mathcal{H}^1_\Eul$, problem \eqref{5.20} has a unique generalized solution $U_\Eul\in Y^1_\Eul$ given by $U_\Eul(t)=\Eul^1(t)[U_{0\Eul}]$ for all $t\in\R$ and hence continuously depending on $U_{0\Eul}$ in the topologies of the respective spaces. Moreover $U_\Eul$ is a strong solution if and only if $U_{0\Eul}\in D(A_\Eul)$.
Next, if $\dot{U}_\Pot\in \dot{Y}^1_\Pot$ is the unique generalized solution of problem \eqref{3.32} with data $\dot{U}_{0\Pot}=F_\EP U_{0\Eul}$, one has $U_\Eul=\Phi_\PE \dot{U}_\Pot$.
\item The subspace $\mathcal{H}^1_\Eulc$ is invariant under the  flow of $\{\Eul^1(t),t\in\R\}$, so the unbounded operator $-A_\Eulc: D(A_\Eulc)\subset\mathcal{H}^1_\Eulc\to\mathcal{H}^1_\Eulc$ defined by
    \begin{equation}\label{5.21BIS}
     D(A_\Eulc)=D(A_\Eul)\cap\mathcal{H}^1_\Eulc,\qquad A_\Eulc U=A_\Eul U\quad\text{for all $U\in D(A_\Eulc)$,}
    \end{equation}
    generates on $\mathcal{H}^1_\Eulc$ the subspace group $\{\Eul_c^1(t),t\in\R\}$ given by $\Eul_c^1(t)=\Eul^1(t)_{|\mathcal{H}^1_\Eulc}$ for all $t\in\R$, which is similar to the group $\{\Pot_c^1(t),t\in\R\}$
    in Corollary~\ref{Corollary 3.3} since we have
    \begin{equation}\label{5.22}
  \Eul_c^1(t)=F_\PcEc\Pot_c^1(t) F_\EcPc\qquad\text{for all $t\in\R$.}
\end{equation}
Consequently part I) continues to hold when replacing the spaces $\mathcal{H}^1_\Eul$, $Y^1_\Eul$, $\dot{Y}^1_\Pot$ with the spaces
$\mathcal{H}^1_\Eulc$, $Y^1_\Eulc$, $\dot{Y}^1_\Potc$, the operator $A_\Eul$ with $A_\Eulc$  and restricting problems \eqref{3.32} and \eqref{5.20} to the spaces $\dot{\mathcal{H}}^1_\Potc$ and $\mathcal{H}^1_\Eulc$.
\end{enumerate}
\end{thm}
\begin{proof} To prove part I) we use  arguments similar to those in the proof of Theorem~\ref{Theorem 4.1}. Indeed, by standard semigroup theory, \eqref{5.21} defines on $\mathcal{H}^1_\Eul$ the  strongly continuous
group $\{\Eul^1(t),t\in\R\}$ similar to $\{\Pot^1(t),t\in\R\}$, having as generator the operator $-B_2$ defined by $D(B_2)=F_\PE D(\dot A_\Pot)$ and $B_2=F_\PE \dot A_\Pot F_\EP$. By \eqref{3.29}, \eqref{5.17} and Lemma~\ref{Lemma 5.1} we then get $D(B_2)=F_\PE \dot{D}^1_\Pot=D^1_\Eul=D(A_\Eul)$, dense in $\mathcal{H}^1_\Eul$. Moreover, for any $U=(p,\vvect,v,z)\in D(A_\Eul)$, using \eqref{3.29} and \eqref{5.9} one has
$$B_2U=F_\PE\dot{A}_\Pot \begin{pmatrix}\dot u\\v\\\tfrac 1{\rho_0}p\\z\end{pmatrix}
=\begin{pmatrix}-\tfrac 1{\rho_0}\pi_0 p\\-z\\\tfrac B{\rho_0}\Div \vvect\\
\frac 1\mu\left[-\DivGamma(\sigma\nabla_\Gamma v)+\delta z+\kappa v+p_{|\Gamma_1}\right]
\end{pmatrix}=A_\Eul U,
$$
where $\dot u\in \dot H^2(\Omega)$ is the unique solution of $-\nabla\dot u=\vvect$. Hence $A_\Eul=B_2$, so $-A_\Eul$ generates the
group $\{\Eul^1(t),t\in\R\}$. The proof of part I) can then be completed using Lemma~\ref{Lemma 2.1}--i) and iii).

To prove part II) we first remark that, by Corollary~\ref{Corollary 3.3} and \eqref{5.12}, formula \eqref{5.22} defines on $\mathcal{H}^1_\Eulc$ the group $\{\Eul_c^1(t),t\in\R\}$ similar to  the group $\{\Pot_c^1(t),t\in\R\}$. So by \eqref{5.21} and the invariance of $\dot{\mathcal{H}}^1_\Potc$ with respect to the group $\{\Pot^1(t),t\in\R\}$ asserted in Lemma~\ref{Lemma 3.4}, we get that $\mathcal{H}^1_\Eulc$ is invariant with respect to $\{\Eul^1(t),t\in\R\}$,
 so $\{\Eul_c^1(t),t\in\R\}$ is the subspace group of $\{\Eul^1(t),t\in\R\}$ which has -$A_\Eulc$ as generator. To complete the proof of part II) is then trivial.
\end{proof}
We now set , for $n\in\widetilde{\N}$, $n\le r$, the Fr\'{e}chet spaces
\begin{equation}\label{5.23}
\mathcal{T}^n_\Eul=\{U_\Eul\in Y^n_\Eul: U_\Eul\text{ is a generalized solution of \eqref{5.19}}\}, \quad \mathcal{T}^n_\Eulc=\mathcal{T}^n_\Eul\cap Y^n_\Eulc,
\end{equation}
endowed with the topology inherited from $Y^n_\Eul$. The following result is a trivial consequence of formulas \eqref{3.15}, \eqref{5.5}, Lemma~\ref{Lemma 5.1} and Theorem~\ref{Theorem 5.1}.
\begin{cor}\label{Corollary 5.1} For all $n\in\widetilde{\N}$, $n\le r$, the operators $\Phi_\PE$ and $\Phi_\EP$ in \eqref{5.11} restrict to bijective isomorphisms between $\dot{\mathcal{T}}^n_\Pot$ and $\mathcal{T}^n_\Eul$, being each one the inverse of the other, and also to bijective isomorphisms between $\dot{\mathcal{T}}^n_\Potc$ and $\mathcal{T}^n_\Eulc$.
\end{cor}
Abstract regularity properties of the solutions of \eqref{5.20} are then given as follows.
\begin{thm}[\bf Regularity for \eqref{5.20}]\label{Theorem 5.2}\phantom{A}
\renewcommand{\labelenumi}{{\Roman{enumi})}}
\begin{enumerate}
\item For all $n\in\widetilde{\N}$, $2\le n\le r$, one has $D(A^{n-1}_\Eul)=D^{n-1}_\Eul$, the respective norms being equivalent, so the operator $-\,\, {_n}A_\Eul$ given by \eqref{2.18} generates on $D^{n-1}_\Eul$ the strongly continuous group $\{\Eul^n(t), t\in\R\}$ given by
\begin{equation}\label{5.23BB}
\Eul^n(t)= F_\PE\Pot^n(t)F_\EP\qquad\text{for all $t\in\R$,}
\end{equation}
and hence similar to the group $\{\Pot^n(t),t\in\R\}$ in Proposition~\ref{Proposition 3.1} --II).
Consequently  for any $U_{0\Eul}\in\mathcal{H}^1_\Eul$, denoting by $U_\Eul$ the generalized solution of \eqref{5.20}, for any $n\in\widetilde{\N}$, $2\le n\le r$,  one has  $U_{0\Eul}\in D^{n-1}_\Eul$ if and only if $U_\Eul\in Y^n_\Eul$,
 and in this case $U_\Eul$  continuously depends on it $U_{0\Eul}$ in the topologies of the respective spaces.
\item All assertions in part I) continue to hold when restricting problem \eqref{5.20} to $\mathcal{H}^1_\Eulc$, provided one replaces the spaces $\mathcal{H}^1_\Eul$, $D^{n-1}_\Eul$, $Y^n_\Eul$ with the spaces $\mathcal{H}^1_\Eulc$, $D^{n-1}_\Eulc$, $Y^n_\Eulc$, the operators $A_\Eul$ and ${_n}A_\Eul$ with $A_\Eulc$ and ${_n}A_\Eulc$ and, finally, the group $\{\Eul^n(t),t\in\R\}$ with the group $\{\Eulc^n(t),t\in\R\}$ defined by $\Eul_c^n(t)=\Eul^n(t)_{|D^{n-1}_\Eulc}$ for $t\in\R$.
\end{enumerate}
\end{thm}
\begin{proof} To prove part I) we remark that, since $A_\Eul=F_\PE \dot A_\Pot F_\EP$ (see the proof of previous result), using Lemma~\ref{Lemma 5.1} and \eqref{2.16} one gets, by induction, that $D(A^{n-1}_\Eul)=D^{n-1}_\Eul$, the respective norms being equivalent, and that \eqref{5.23BB} defines on it the strongly continuous group $\{\Eul^n(t), t\in\R\}$, which has $-\,{_n}A_\Eul$ as generator. The proof of part I) can then be completed by using Corollary~\ref{Corollary 5.1} and Lemma~\ref{Lemma 5.1}, since they yield the implications $U_{0\Eul}\in D^{n-1}_\Eul \Leftrightarrow F_\EP U_{0,\Eul}\in\dot{D}^{n-1}_\Pot \Leftrightarrow \dot{U}_\Pot\in\dot{Y}^n_\Pot
\Leftrightarrow U_\Eul\in Y^n_\Eul$,
where $\dot{U}_\Pot=\Phi_\EP U_\Eul$. The proof of part II) uses similar arguments.
\end{proof}
\subsection{Solutions of $(\Eul)$, $(\Eulc)$, $(\Eul_0)$ and $(\Eul_0^c)$}\label{Section 5.4}
To apply the abstract results in \S~\ref{Section 5.3} to all $(\Eul)$ -- related problems we make precise, at first, which type of solutions of them we shall consider, recalling that  $(\Eul)_3$ is  implicit in the definition of $X^1_\Eul$.
\begin{definition}\label{Definition 5.1}
We say that
\renewcommand{\labelenumi}{{\roman{enumi})}}
\begin{enumerate}
\item $(p,\vvect,v)\in X^2_\Eul$ is a {\em strong solution} of $(\Eul)$ provided $(\Eul)_1$ -- $(\Eul)_2$ hold a.e. in $\R\times\Omega$ and $(\Eul)_4$ -- $(\Eul)_5$ hold a.e. on $\R\times\Gamma_1$, where $p$ and $\vvect$ on $\R\times\Gamma_1$ are taken in the pointwise trace sense given in \S~\ref{Section 2.2};
\item $(p,\vvect,v)\in X^1_\Eul$ is a {\em generalized solution} of $(\Eul)$ provided it is the limit in $X^1_\Eul$ of a sequence of strong solutions of it;
\item $(p,\vvect,v)\in X^1_\Eul$ is a {\em weak solution} of $(\Eul)$ provided the distributional identities
\begin{align}\label{5.24}
\int_{-\infty}^\infty\int_\Omega [p\varphi_t+B\vvect\cdot\nabla\varphi]+&\int_{-\infty}^\infty\int_{\Gamma_1}v_t\varphi=0,\\
\label{5.25}
\int_{-\infty}^\infty\int_\Omega  [\rho_0\vvect\cdot \phivect_t  +p\Div \phivect] &\\
 \notag+ \int_{-\infty}^\infty\int_{\Gamma_1}&\left[\mu v_t\psi_t-\sigma(\nabla_\Gamma v,\nabla_\Gamma \overline{\psi})_\Gamma-\delta v_t\psi
 -\kappa v\psi\right]=0,
 \end{align}
 hold for all $\varphi\in C^r_c(\R\times\R^3)$ and $\phivect\in C^{r-1}_c(\R\times\R^3)^3$ such that $\phivect\cdot\boldsymbol{\nu}=0$ on $\R\times\Gamma_0$, where $\psi=-\phivect\cdot\boldsymbol{\nu}$ on $\R\times\Gamma_1$.
\end{enumerate}
Moreover solutions of $(\Eul)$ of the types i)--iii) above are said to be solutions of the same type of:
j) problem $(\Eulc)$ when also \eqref{1} holds; jj) problem $(\Eul_0)$ when also  $(\Eul_0)_6$--$(\Eul_0)_7$ hold in $X^1_\Eul$;
jjj) problem $(\Eul_0^c)$ when both j) and jj) hold.
\end{definition}
Trivially strong solutions in the definition above are also generalized ones. Moreover strong and generalized solutions correspond to the homologous ones of \eqref{5.19} and \eqref{5.20} as the following result shows.
\begin{lem}\label{Lemma 5.2}
The triple  $(p,\vvect,v)$ is a strong or generalized solution of $(\Eul)$ if and only if $p$, $\vvect$ and $v$ are the first three components of a solution $U_\Eul=(p,\vvect,v,z)$ of \eqref{5.19} of the same type, and in this case $v_t=z$. The same relation occurs
between:  solutions of $(\Eulc)$ and solutions of the restriction of \eqref{5.19} to $\mathcal{H}^1_\Eulc$;
 solutions of $(\Eul_0)$ and of \eqref{5.20};
solutions of $(\Eul_0^c)$ and solutions of the restriction of \eqref{5.20} to $\mathcal{H}^1_\Eulc$.
\end{lem}
\begin{proof}
By Definition~\ref{Definition 5.1} only the first assertion needs a proof. Since, by Theorem~\ref{Theorem 5.2}--I), strong solutions of \eqref{5.19} belong to $Y^2_\Eul$, using \eqref{5.18} the assertion is trivial for strong solutions. Using this fact, given any generalized solution $(p,\vvect,v)\in X^1_\Eul$ of $(\Eul)$ the quadruple $(p,\vvect,v,v_t)\in Y^1_\Eul$ is a generalized solution of \eqref{5.19}. The converse implication, with the identity $v_t=z$, is proved as in the proof of Lemma~\ref{Lemma 4.5}.
\end{proof}
By Lemma~\ref{Lemma 5.2} generalized and strong solutions of $(\Eul)$ --  related problems  naturally arise from Theorems~\ref{Theorem 5.1} and \ref{Theorem 5.2}.
While strong solutions are a.e. classical solutions, generalized solutions of $(\Eul)$ solve all equations (but $(\Eul)_3$) in a quite indirect sense. Since, also in this case, solutions in $X^1_\Eul$ are not regular enough to be a.e. solutions, it is natural also for this problem to consider $(\Eul)_1$ -- $(\Eul)_2$ and $(\Eul)_4$ -- $(\Eul)_5$ in a distributional sense.
On the other hand, while for $(p,\vvect,v)\in X^1_\Eul$ equations $(\Eul)_1$  and  $(\Eul)_2$ have the natural distributional forms \begin{alignat}2\label{5.26}
 &\int_{\R\times\Omega} p\varphi_t+B\vvect\cdot\nabla\varphi=0&&\qquad\text{for all $\varphi\in\mathcal{D}(\R\times\Omega)$,}\\
 \label{5.27}
&\int_{\R\times\Omega}\rho_0\vvect\cdot \phivect_t+p \Div \phivect=0&&\qquad\text{for all $\phivect\in \mathcal{D}(\R\times\Omega)^3$,}
  \end{alignat}
 equations $(\Eul)_4$ and  $(\Eul)_5$ can not be written in the sense of distributions unless the terms $p$ and $\vvect\cdot\boldsymbol{\nu}$ in them have some trace sense on $\R\times\Gamma_1$.
This type of difficulty,  also arising for problem $(\Lagr)$, was solved in Definition~\ref{Definition 5.1}--iii) by combining equations $(\Eul)_1$ and $(\Eul)_5$ in the single distributional identity \eqref{5.24} and $(\Eul)_2$ and $(\Eul)_4$ in the single distributional identity \eqref{5.25}.
Of course one can be more precise when solutions are more regular. As to the term $p$ one can simply ask that $p\in L^1_\loc(\R; H^1(\Omega))$. In this case  equation $(\Eul)_4$, formally  identical to $(\Lagr')_4$, has the distributional form
\begin{equation}
 \label{4.84}
 \int_{-\infty}^\infty\int_{\Gamma_1}\mu v_t\psi_t-\sigma(\nabla_\Gamma v,\nabla_\Gamma \overline{\psi})_\Gamma-\delta v_t\psi
 -\kappa v\psi-p\psi=0
 \end{equation}
 for all $\psi\in C^r_c(\R\times\Gamma_1)$.
As for the term $\vvect\cdot\boldsymbol{\nu}$,  recalling \eqref{5.6}, one can give a sense to it by asking that $\vvect\in L^1_\loc(\R; H(\Div,\Omega))$. Indeed, using \eqref{5.8} in this case one has $\vvect\cdot\boldsymbol{\nu}\in L^1_\loc(\R; H^{-1/2}(\Gamma))$ and equation $(\Eul)_5$ has the natural distributional form
\begin{equation}\label{5.28}
\int_{-\infty}^\infty \langle \vvect\cdot\boldsymbol{\nu}, \varphi_{|\Gamma}\rangle_{H^{1/2}(\Gamma)} +\int_{\R\times\Gamma_1}v_t\varphi=0\qquad\text{for all $\varphi\in C^r_c(\R\times\Gamma)$.}
\end{equation}
After these preliminary remarks we can finally state the following result, which  shows that Definition~\ref{Definition 5.1}--iii) is the closest possible approximation of the notion of distributional solutions of $(\Eul)$. It will be also useful in the sequel.
\begin{prop}\label{Proposition 5.1}
Let $(p,\vvect,v)\in X^1_\Eul$ be such that $p\in L^1_\loc(\R; H^1(\Omega))$ and $\vvect\in L^1_\loc(\R; H(\Div,\Omega))$.
Then $(p,\vvect,v)$ is a weak solution of $(\Eul)$ if and only if it satisfies the distributional equations \eqref{5.26} -- \eqref{5.28}.
\end{prop}
\begin{proof} The proof is organized in five Steps.

\noindent{\bf Step 1.} We claim that \eqref{5.26} holds if and only if
\begin{equation}\label{5.29}
 p\in W^{1,1}_\loc (\R;L^2(\Omega))\quad\text{and}\quad p_t+B\Div \vvect=0
 \quad\text{in $L^1_\loc (\R;L^2(\Omega))$.}
 \end{equation}
 When \eqref{5.29} holds, multiplying it by $\varphi\in \mathcal{D}(\R\times\Omega)$, integrating by parts in time and using \eqref{5.8} one easily gets \eqref{5.26}. Conversely, when \eqref{5.26} holds, taking test functions $\varphi(t,x)=\varphi_1(t)\varphi_0(x)$, $\varphi_0\in \mathcal{D}(\Omega)$, $\varphi_1\in \mathcal{D}(\R)$, and using \eqref{5.8} again we get
 $ \int_\Omega\varphi_0 \int_{-\infty}^\infty p\varphi_1'-B\Div \vvect \varphi_1=0$.
Since $\varphi_0$ is arbitrary by density we can take $\varphi_0\in L^2(\Omega)$ in it, so obtaining
 $\int_{-\infty}^\infty p\varphi_1'-B\Div \vvect \varphi_1=0$ in $L^2(\Omega)$, for all $\varphi_1\in \mathcal{D}(\R)$,
 which gives \eqref{5.29} and proves our claim.

\noindent{\bf Step 2.} We claim that \eqref{5.27} holds if and only if
\begin{equation}\label{5.31}
 \vvect\in W^{1,1}_\loc (\R;H^0_{\curl 0}(\Omega))\quad\text{and}\quad \rho_0\vvect_t+\nabla p=0
 \quad\text{in $L^1_\loc (\R;H^0_{\curl 0}(\Omega))$.}
 \end{equation}
 When \eqref{5.31} holds, multiplying it by $\phivect\in \mathcal{D}(\R\times\Omega)^3$ and using the same arguments in Step 1 one easily gets \eqref{5.27}.  Conversely, when \eqref{5.27} holds, taking test functions $\phivect(t,x)=\varphi_(t)\phivect_0(x)$, $\varphi\in \mathcal{D}(\R)$, $\phivect_0\in \mathcal{D}(\Omega)^3$, using the arguments in Step 1 again  we get
$\int_{-\infty}^\infty \rho_0\vvect\,\varphi'-\nabla p \,\varphi=0$ in $L^2(\Omega)$, for all $\varphi\in \mathcal{D}(\R)$,
which gives \eqref{5.31}, proving our claim.

\noindent{\bf Step 3.} We claim that, when \eqref{5.26} holds, \eqref{5.24} and \eqref{5.28} are equivalent. Indeed, using Step 1 and \eqref{5.8}, equation \eqref{5.24} can be rewritten as \eqref{5.28}, but for test functions $\varphi$, belonging to  $C^r_c(\R\times\R^3)$, instead that to $C^r_c(\R\times\Gamma)$. On the other hand, while elements of $C^r_c(\R\times\R^3)$ trivially restrict to elements of the second space, using the compactness of $\Gamma$ and \cite[Definition 1.2.1.1, p. 5]{grisvard}, it is straightforward to show that elements of  $C^r_c(\R\times\Gamma)$ can be extended to elements of $C^r_c(\R\times\R^3)$, proving our claim.

\noindent{\bf Step 4.} We claim that, when \eqref{5.27} holds, \eqref{5.25} and \eqref{4.84} are equivalent. Indeed, using Step 2 and integrating \eqref{5.31} by parts in time, \eqref{5.25} can be rewritten as \eqref{4.84}, but for test functions $\psi=-\phivect\cdot\boldsymbol{\nu}$ on $\R\times\Gamma_1$, where $\phivect\in C^{r-1}_c(\R\times\R^3)^3$ such that $\phivect\cdot\boldsymbol{\nu}=0$ on $\R\times\Gamma_0$, instead that for $\psi\in C^r_c(\R\times\Gamma_1)$. On the other hand test functions of the first type trivially belong to $\psi\in C^{r-1}_c(\R\times\Gamma_1)$, and by the same argument used in the proof of Proposition~\ref{Proposition 4.2} one can see that any $\psi\in C^{r-1}_c(\R\times\Gamma_1)$ can be obtained as $\psi=-\phivect\cdot\boldsymbol{\nu}$ on $\R\times\Gamma_1$, where $\phivect\in C^{r-1}_c(\R\times\R^3)^3$ is such that $\phivect\cdot\boldsymbol{\nu}=0$ on $\R\times\Gamma_0$. A standard density argument shows that to take $\psi\in C^{r-1}_c(\R\times\Gamma_1)$  in \eqref{4.84} is equivalent to take $\psi\in C^r_c(\R\times\Gamma_1)$, proving our claim.

\noindent{\bf Step 5. Conclusion.} Recalling that $(p,\vvect,v)$ is a weak solution of $(\Eul)$ when \eqref{5.24} and \eqref{5.25} hold, and noticing that \eqref{5.24} and \eqref{5.25} trivially yield \eqref{5.26} and \eqref{5.27}, using Steps 3--4 one trivially concludes the proof.
\end{proof}

Beside its independent interest Proposition~\ref{Proposition 5.1} also allows to point out the trivial relations occurring among the three types of solutions of $(\Eul)$ introduced in Definition~\ref{Definition 5.1}, completely analogous to the ones occurring for problems $(\Pot)$ and $(\Lagr)$.

\begin{lem}\label{Lemma 5.3} Let $(p,\vvect,v)\in X^1_\Eul$ be a solution of $(\Eul)$ according to Definition~\ref{Definition 5.1}. Then strong $\Rightarrow$ generalized $\Rightarrow$ weak and, if $(p,\vvect,v)\in X^2_\Eul$, weak $\Rightarrow$ strong.
\end{lem}
\begin{proof} Strong solutions are also generalized ones and, by Proposition~\ref{Proposition 5.1}, also weak. Since \eqref{5.24} and \eqref{5.25} are stable with respect to the convergence in $X^1_\Eul$, we then get that generalized $\Rightarrow$ weak.
To prove the  final conclusion let $(p,\vvect,v)\in X^2_\Eul$ be a weak solution. By  Proposition~\ref{Proposition 5.1} it also satisfies  $(\Eul)_1$, and $(\Eul)_2$ and   $(\Eul)_4$ in a distributional sense. Being regular enough it  also satisfies them a.e., so it is a strong solution.
\end{proof}

In view of Theorem~\ref{Theorem 5.1}, which by  Lemmas~\ref{Lemma 5.2} and \ref{Lemma 5.3} provides the existence of generalized and hence weak solutions of $(\Eul)$, it is natural to pose the question of uniqueness of weak solutions of $(\Eul_0)$. The answer is given by the following result.
\begin{thm}[\bf Uniqueness]\label{Theorem 5.3}  Weak solutions of $(\Eul_0)$ are unique.
\end{thm}
\begin{proof} By linearity it is enough to prove that $p_0=0$, $\vvect_0=0$ and $v_0=v_1=0$ yield $(p,\vvect,v)\equiv 0$ in $\R$.
At first we remark that taking in \eqref{5.24} test functions $\varphi(t,x)=\varphi_1(t)\varphi_0(x)$,  $\varphi_1\in C^r_c(\R)$,
$\varphi_0\in C^r_c(\R^3)$, we get
\begin{equation}\label{5.32}
 \int_{-\infty}^\infty \varphi_1'\int_\Omega p\varphi_0+ \int_{-\infty}^\infty B \varphi_1\left(\int_\Omega \vvect\cdot \nabla\varphi_0+\int_{\Gamma_1} v_t\varphi_0\right)=0,
\end{equation}
which, by density holds true for all $\varphi_0\in H^1(\Omega)$. Since
$\int_\Omega p\varphi_0, \int_\Omega\vvect\cdot\nabla\varphi_0,\int_{\Gamma_1}v_t\varphi_0\in C(\R)$, by \eqref{5.32} we then get that
$\int_\Omega p\varphi_0\in C^1(\R)$ and
$\frac d{dt} \int_\Omega p\varphi_0=B\int_\Omega \vvect\cdot\nabla\varphi_0+ B \int_{\Gamma_1}v_t\varphi_0$
for all $\varphi_0\in H^1(\Omega)$. Consequently, setting $\rvect\in C^1(\R; H^0_{\curl 0}(\Omega))$ (since $\vvect\in C(\R;H^0_{\curl 0}(\Omega))$ by
$\rvect(t)=\int_0^t \vvect(\tau)\,d\tau$  for all $t\in\R$,
since $p(0)=p_0=0$,  $v(0)=v_0=0$ and $\rvect(0)=0$, we have
\begin{equation}\label{5.34}
\int_\Omega p(t)\varphi_0=B\int_\Omega \rvect(t)\cdot\nabla\varphi_0+ B \int_{\Gamma_1}v(t)\varphi_0  \qquad\text{for all $t\in\R$.}
\end{equation}
Taking in \eqref{5.34} test functions $\varphi_0\in \mathcal{D}(\Omega)$ we get that
$-B\Div\rvect(t)=p(t)$  in $\mathcal{D}'(\Omega)$, for all $t\in\R$. Since $p(t)\in L^2(\Omega)$ for all $t\in\R$ we consequently get that $\rvect(t)\in H(\Div,\Omega)$ for all $t\in\R$ and thus by \eqref{5.8}, for all $t\in\R$ and $\varphi_0\in H^1(\Omega)$ we have
\begin{equation}\label{5.35}
\int_\Omega p(t)\varphi_0=B\int_\Omega \rvect(t)\cdot \nabla\varphi_0- B\langle \rvect(t)\cdot\boldsymbol{\nu},{\varphi_0}_{|\Gamma}  \rangle_{H^{1/2}(\Gamma)}.
\end{equation}
By comparing \eqref{5.34} and \eqref{5.35}, since the trace operator $H^1(\Omega)\to H^{1/2}(\Gamma)$ is surjective, recalling the splitting \eqref{2.8}, we then get that
$\rvect(t)\cdot\boldsymbol{\nu}=-v(t)$ in $H^{-1/2}(\Gamma)$ for all $t\in\R$. Since $v(t)\in H^1(\Gamma_1)$ we then get
$\rvect(t)\cdot\boldsymbol{\nu}\in H^1(\Gamma)\hookrightarrow H^{1/2}(\Gamma)$. Summarizing we have $\curl \rvect(t)=0$, $\Div \rvect(t)\in L^2(\Omega)$ and $\rvect(t)\cdot\boldsymbol{\nu}\in H^{1/2}(\Gamma)$, hence we can apply
\cite[Chapter 9, Corollary 1, p. 212]{dautraylionsvol3} to conclude that $\rvect(t)\in H^1_{\curl 0}(\Omega)$ for all $t\in\R$. Hence $\rvect(t)$ solves problem \eqref{4.12} with $w=p(t)/\rho_0$ and, by Lemma~\ref{Lemma 4.2}, we $\rvect\in C(\R; H^1_{\curl 0}(\Omega))$ and $\rvect(t)\cdot\boldsymbol{\nu}=-v(t)$ on $\Gamma_1$,  $\rvect(t)\cdot\boldsymbol{\nu}=0$ on $\Gamma_0$ in the trave sense. Consequently $(\rvect,v)\in X^1_\Lagr$, so $(\rvect,v)\in X^1_{\Lagr}$ and $\rvect_t=\vvect$. Consequently \eqref{5.25} translates into \eqref{4.49} and $(\rvect,v)$ is a weak solution of $(\Lagr_0)$ with $\rvect_0=\rvect_1=\vvect_0=0$ and $v_0=v_1=0$, by Theorem~\ref{Theorem 4.3} we get that $\rvect\equiv 0$ and $v\equiv 0$ in $\R$, so also $\vvect=\rvect_t\equiv 0$ in $\R$, concluding the proof.
\end{proof}
\subsection{Main results for problems $(\Eul)$, $(\Eul_0)$, $(\Eul^c)$ and $(\Eul^c_0)$ }\label{Section 5.5} We can finally prove Theorems~\ref{Theorem 1.3} and \ref{Theorem 1.5} and deal with  optimal regularity issues.
\begin{proof}[\bf Proof of Theorem~\ref{Theorem 1.3}]
By combining Theorem~\ref{Theorem 5.1} and Lemma~\ref{Lemma 5.2} for all $U_{0\Eul}\in \mathcal{H}^1_\Eul$ problem $(\Eul_0)$ has a unique generalized solution $(p,\vvect,v)\in X^1_\Eul$, continuously depending on  $U_{0\Eul}$, which is the unique generalized solution of problem $(\Eul^c_0)$ if and only if $U_{0\Eul}\in \mathcal{H}^1_\Eulc$. By Lemma~\ref{Lemma 5.3} the triple $(p,\vvect,v)$ is also a weak solution of $(\Eul_0)$ (or of $(\Eul^c_0)$, when $U_{0\Eul}\in \mathcal{H}^1_\Eulc$) which, by Theorem~\ref{Theorem 5.3}, is unique among them. Moreover, again by Theorem~\ref{Theorem 5.1} and  Lemma~\ref{Lemma 5.2}, the solution is strong if and only if $U_{0\Eul}\in \mathcal{H}^2_\Eul$ and $(\vvect_0,v_1)\in \mathbb{H}^1_\Lagr$, i.e. $U_{0\Eul}\in D(A_\Eul)$, which is dense in $\mathcal{H}^1_\Eul$. The same conclusions trivially hold for problem $(\Eul^c_0)$. The continuous dependence in this case follows
by Theorem~\ref{Theorem 5.2}. Finally the energy identity \eqref{1.15} can be obtained as in the proof of Theorem~\ref{Theorem 1.2}  directly from \eqref{1.5}, or re--deriving it for strong solutions and then using a density argument.
\end{proof}

The uniqueness of weak solution of $(\Eul_0)$ in Theorem~\ref{Theorem 5.3} together with the existence of a generalized and hence weak solution of it shows that solutions of $(\Eul)$ are {\em equivalently} weak or generalized. Combining this remark with Lemma~\ref{Lemma 5.3} we thus obtain that, as for all previous problems,
{\em all types of solutions in Definition~\ref{Definition 5.1} coincide, strong solutions being defined only in the class $X^2_\Eul$}. In the sequel we consequently shall deal only with weak solutions of $(\Eul)$ and  $(\Eulc)$, i.e. with elements of the spaces $\mathcal{S}^n_\Eul$ and $\mathcal{S}^n_\Eulc $ defined in \eqref{1.18}. Recalling the spaces
$\mathcal{T}^n_\Eul$ and $\mathcal{T}^n_\Eulc$ defined in \eqref{5.23},
 it is useful to point out the following result.

\begin{lem}\label{Lemma 5.4}
The operator $\mathcal{I}''\in\mathcal{L}(X^1_\Eul;Y^1_\Eul)$ defined by
$\mathcal{I}''(p,\vvect,v)=(p,\vvect,v,v_t)$
restricts for each $n\in\widetilde{\N}$, $n\le r$, to a bijective isomorphism $\mathcal{I}^n_\Eul\in\mathcal{L}(\mathcal{S}^n_\Eul;\mathcal{T}^n_\Eul)$ with inverse $\left(\mathcal{I}^n_\Eul\right)^{-1}\in\mathcal{L}(\mathcal{T}^n_\Eul;\mathcal{S}^n_\Eul)$ simply given by
\begin{equation}\label{5.38}
\left(\mathcal{I}^n_\Eul\right)^{-1}(p,\vvect,v,z)=(p,\vvect,v).
\end{equation}
Moreover $\mathcal{I}^n_\Eul$ and $\left(\mathcal{I}^n_\Eul\right)^{-1}$ further restrict to $$\mathcal{I}^n_\Eulc\in\mathcal{L}(\mathcal{S}^n_\Eulc;\mathcal{T}^n_\Eulc)\qquad\text{and}\quad  \left(\mathcal{I}^n_\Eulc\right)^{-1}\in\mathcal{L}(\mathcal{T}^n_\Eulc;\mathcal{S}^n_\Eulc).$$
\end{lem}
\begin{proof}The first assertion follows by \eqref{5.5}, Lemma~\ref{Lemma 5.2} and the equivalence remarked just before the statement. The second one follows by the first one once since
$X^n_\Eulc=\left\{(p,\vvect,v)\in X^n_\Eul: \int_\Omega \frac{\partial^ip}{\partial t^i}=B\int_{\Gamma_1}\frac{\partial^iv}{\partial t^i}\quad\text{for $i=0,\ldots,n-1$}\right\}$.
\end{proof}
We can now give the
\begin{proof}[\bf Proof of Theorem~\ref{Theorem 1.5}]
By Proposition~\ref{Proposition 3.2}, Lemma~\ref{Lemma 5.1}, Lemma~\ref{Lemma 5.4}, \eqref{5.9} and \eqref{5.11} for any
$(u,v)\in \mathcal{S}^1_\Pot$ we have $\left(\mathcal{I}^1_\Eul\right)^{-1}\cdot \Phi_\PE \cdot \mathcal{J}^1_\Pot(u,v)=(\rho_0 u_t,-\nabla u,v)$, so \eqref{1.28} defines the operators
\begin{equation}\label{5.38}
\Psi_\PE=\left(\mathcal{I}^1_\Eul\right)^{-1}\!\!\!\!\cdot \Phi_\PE \!\!\cdot \mathcal{J}^1_\Pot\in \mathcal{L}(\mathcal{S}^1_\Pot;\mathcal{S}^1_\Eul), \,\,  \Psi_\PcEc=\left(\mathcal{I}^1_\Eulc\right)^{-1}\!\!\!\!\cdot \Phi_\PcEc \!\!\! \cdot \mathcal{J}^1_\Potc\in\mathcal{L}(\mathcal{S}^1_\Potc;\mathcal{S}^1_\Eulc)
\end{equation}
given in the statements, which are surjective and satisfy $\text{Ker } \Psi_\PE = \text{Ker } \Psi_\PcEc =\C_{X_\Pot}$. One consequently gets the bijective isomorphisms
\begin{equation}\label{5.39}
\dot\Psi_\PE=\left(\mathcal{I}^1_\Eul\right)^{-1}\!\!\!\!\cdot \Phi_\PE \!\!\cdot \dot{\mathcal{J}}^1_\Pot\in \mathcal{L}(\dot{\mathcal{S}}^1_\Pot;\mathcal{S}^1_\Eul), \,\, \dot\Psi_\PcEc\!\!=\left(\mathcal{I}^1_\Eulc\right)^{-1}\!\!\!\!\cdot \Phi_\PcEc \! \!\cdot \dot{\mathcal{J}}^1_\Potc\!\!\in \mathcal{L}(\mathcal{S}^1_\Potc;\mathcal{S}^1_\Eulc)
\end{equation}
given in \eqref{1.30}. Their inverses are the operators
\begin{equation}\label{5.40}
\dot\Psi_\EP=\left(\dot{\mathcal{J}}^1_\Pot\right)^{-1}\!\!\!\!\!\!\!\cdot \Phi_\EP \cdot \mathcal{I}^1_\Eul\in \mathcal{L}(\mathcal{S}^1_\Eul;\dot{\mathcal{S}}^1_\Pot), \,\, \dot\Psi_\EcPc=\left(\dot{\mathcal{J}}^1_\Potc\right)^{-1}\!\!\!\!\!\!\!\cdot \Phi_\EcPc \!\!\cdot \mathcal{I}^1_\Eulc\!\in\! \mathcal{L}(\mathcal{S}^1_\Eulc;\dot{\mathcal{S}}^1_\Potc),
\end{equation}
so formulas \eqref{1.30} trivially hold.

To check that $\dot\Psi_\EP$ is given by \eqref{1.31} -- \eqref{1.32} we recall that, by \eqref{5.9}, \eqref{5.11} and Lemma~\ref{Lemma 5.4}, for all $(p,\vvect,v)\in\mathcal{S}^1_\Eul$ and $t\in\R$ we have
$$[\Phi_\EP\mathcal{I}^1_\Eul(p,\vvect,v)](t)=F_\EP(p(t), \vvect(t), v(t),v_t(t))=(\dot u(t), v(t), p(t)/\rho_0, v_t(t))$$
where $\dot u(t)\in \dot H^1(\Omega)$ is the unique solution of $-\nabla \dot u(t)=\vvect(t)$. Since $\mathcal{J}^1_\Pot$ is surjective there is $(u,v)\in\mathcal{S}^1_\Pot$ such that
\begin{equation}\label{5.40BIS}
(\dot u(t), v(t), p(t)/\rho_0, v_t(t))=(\pi_0 u(t),v(t),u_t(t),v_t(t))\qquad\text{for all $t\in\R$,}
\end{equation}
 that is
 \begin{equation}\label{5.41}
  \pi_0 u(t)=\dot u(t),\quad\text{and}\quad p(t)=\rho_0u_t(t)\qquad\text{for all $t\in\R$.}
 \end{equation}
 Clearly \eqref{5.41} implies \eqref{1.31}--\eqref{1.32}. Since $\dot\Psi_\EP$ is bijective, by \eqref{1.28} one gets \eqref{1.33}, while \eqref{1.34} is trivial. To prove \eqref{1.35} we remark that, as we just proved,  $(p,\vvect,v)=\Psi_\PE(u,v)$  is equivalent to \eqref{5.40BIS}, that is  \eqref{5.40}. By uniqueness then it is equivalent to \eqref{5.40} at $t=0$, that is to
$-\nabla u_0=\vvect_0$ and $p_0=\rho_0u_1$.
\end{proof}
We now give an optimal regularity result for problems $(\Eul_0)$,  $(\Eul^c_0)$ and show that the operators in Theorem~\ref{Theorem 1.5} have a good behavior with respect to regularity classes.

\begin{thm}[\bf Regularity for $(\Eul_0)$, $(\Eul^c_0)$ and $(\Eul) \rightleftarrows (\Pot)$,  $(\Eulc) \rightleftarrows (\Potc)$]\label{Theorem 5.4}
For all data $U_{0\Eul}=(p_0,\vvect_0, v_0,v_1)\in\mathcal{H}^1_\Eul$ and $n\in\widetilde{\N}$, $2\le n\le r$, the weak solution $(p,\vvect,v)$ given in Theorem~\ref{Theorem 1.3} belongs to $X_\Eul^n$ if and only if $U_{0\Eul}\in D^{n-1}_\Eul$, and in this case
$(p,\vvect,v)$ continuously depends on the data $U_{0\Eul}$ in the topologies of the respective spaces.
In previous assertions we can replace problem $(\Eul_0)$ with problem $(\Eul^c_0)$ provided we respectively replace $H^1_\Eul$, $D^{n-1}_\Eul$ and $X^n_\Eul$ with $H^1_\Eulc$, $D^{n-1}_\Eulc$ and $X^n_\Eulc$.

Moreover the operators $\Psi_\PE$ and $\Psi_\PcEc$ in Theorem~\ref{Theorem 1.5} respectively restrict to surjective operators $\Psi^n_\PE\in\mathcal{L}(\mathcal{S}^n_\Pot;\mathcal{S}^n_\Eul)$ and  $\Psi^n_\PcEc\in\mathcal{L}(\mathcal{S}^n_\Potc;\mathcal{S}^n_\Eulc)$ such that $\ker \Psi^n_\PE=\ker \Psi^n_\PcEc=\C_{X_\Pot}$.
Consequently the bijective isomorphisms $\dot{\Psi}_\PE$,  $\dot{\Psi}_\PcEc$, $\dot{\Psi}_\EP$ and $\dot{\Psi}_\EcPc$ in the same result respectively restrict to bijective isomorphisms $\dot{\Psi}_\PE^n\in\mathcal{L}(\dot{\mathcal{S}}^n_\Pot;\mathcal{S}^n_\Eul)$,
  $\dot{\Psi}_\PcEc^n\in\mathcal{L}(\dot{\mathcal{S}}^n_\Potc;\mathcal{S}^n_\Eulc)$,  $\dot{\Psi}_\EP^n\in\mathcal{L}(\mathcal{S}^n_\Eul;\dot{\mathcal{S}}^n_\Pot)$,
  and $\dot{\Psi}_\EcPc^n\in\mathcal{L}(\mathcal{S}^n_\Eulc;\dot{\mathcal{S}}^n_\Potc)$.
\end{thm}
\begin{proof} The optimal regularity for problem $(\Eul_0)$, and hence $(\Eul^c_0)$,  follows by combining Theorem~\ref{Theorem 5.2} and Lemma~\ref{Lemma 5.4}. Moreover the restriction properties of the operators $\dot{\Psi}_\PE$,  $\dot{\Psi}_\PcEc$, $\dot{\Psi}_\EP$ and $\dot{\Psi}_\EcPc$ follow, thank to formulas \eqref{5.38} -- \eqref{5.40}, by combining Proposition~\ref{Proposition 3.2},  Lemmas~\ref{Lemma 5.1}  and \ref{Lemma 5.4}.
\end{proof}
\subsection{The relation between $(\Eulc)$ and $(\Lagr)$ or $(\Lagr')$}\label{Section 5.6}
We start by completing the proofs of the main results stated in \S~\ref{intro}.
\begin{proof}[\bf Proof of Corollary~\ref{Corollary 1.2}] The result is obtained by simply combining Theorems~\ref{Theorem 1.4} and \ref{Theorem 1.5}. Indeed we define the bijective isomorphism $\Psi_{\Eulc\Lagr}$ by \eqref{1.38}, so its inverse  $\Psi_{\Lagr\Eulc}$ is given by $\Psi_\LEc=\dot{\Psi}_\PcEc\cdot \dot{\Psi}_\LPc$ as stated. To evaluate the action of these two isomorphism we use Theorems~\ref{Theorem 1.4} and \ref{Theorem 1.5}.
In particular, for any $(p,\vvect,v)\in \mathcal{S}^1_\Eulc$, by Theorem~\ref{Theorem 1.5} we have $\dot{\Psi}_{\Eulc\Potc}(p,\vvect,v)=(u,v)+\C_{X_\Pot}$, with $u$ given by \eqref{1.32}, and \eqref{1.33} holds true.
Hence, by Theorem~\ref{Theorem 1.4}, $\Psi_{\Eulc\Lagr}(p,\vvect,v)=\dot\Psi_{\Potc\Lagr}[(u,v)+\C_{X_\Pot}]=(\rvect,v)$ where $\rvect$ is exactly given by \eqref{1.36}, and \eqref{1.37} holds, also getting \eqref{5.52}.
Conversely, for any $(\rvect,v)\in \mathcal{S}^1_\Lagr$, by Theorem~\ref{Theorem 1.4} we have
$\dot\Psi_{\Lagr\Potc}(\rvect,v)=(u,v)+\C_{X_\Pot}$, where $u$ is given by \eqref{1.26} and \eqref{1.22} holds true.
Hence, by Theorem~\ref{Theorem 1.5}, we get \eqref{1.39}. Finally \eqref{1.40} follows by \eqref{1.27} and \eqref{1.35}.
\end{proof}
Although, as we have just seen, the relation between solutions of $(\Eulc)$ and $(\Lagr)$ follow by Theorems~\ref{Theorem 1.4} and \ref{Theorem 1.5}, this relation can be seen at a more abstract level by considering the groups $\{\Eul^1_c(t),t\in\R\}$ and $\{\Lagr^1(t),t\in\R\}$.
Indeed we combine the main isomorphisms between phase spaces in \S~\ref{Section 4.2} and \S~\ref{Section 5.2} by setting
\begin{equation}\label{5.43}
 F_\EcL=F_\PcL\cdot F_\EcPc\in\mathcal{L}(\mathcal{H}^1_\Eulc;\mathcal{H}^1_\Lagr),\quad
 F_\LEc=F_\PcEc\cdot F_\LPc\in\mathcal{L}(\mathcal{H}^1_\Lagr;\mathcal{H}^1_\Eulc),
\end{equation}
trivially being the inverse of each other. By \eqref{4.26} -- \eqref{4.27}, \eqref{4.29} and \eqref{5.9} they are given by
\begin{equation}\label{5.44}
F_\EcL(p,\vvect,v,z)=(\rvect,v,\vvect,z),\qquad   F_\LEc(\rvect,v,\svect,z)=(-B\Div\rvect,\svect,v,z),
\end{equation}
where, in the first identity, $\rvect\in H^1_{\curl 0}(\Omega)$ is the unique solution of \eqref{4.12} with $w=p/\rho_0$.
We then get the following result.

\begin{thm}\label{Theorem 5.5} The groups $\{\Lagr^1(t),t\in\R\}$ and $\{\Eul^1_c(t),t\in\R\}$ given in Theorems~\ref{Theorem 4.1} and \ref{Theorem 5.1} are similar since
\begin{equation}\label{5.45}
 \Lagr^1(t)=F_\EcL\cdot\Eul^1_c(t)\cdot F_\LEc\qquad\text{for all $t\in\R$.}
\end{equation}
Moreover, for all $n\in\N$, $2\le n\le r$, the similarity \eqref{5.45} extends to the groups
$\{\Lagr^n(t),t\in\R\}$ and $\{\Eul^n_c(t),t\in\R\}$ given in Theorems~\ref{Theorem 4.2} and \ref{Theorem 5.2} since
\begin{equation}\label{5.46}
 \Lagr^n(t)=F_\EcL\cdot\Eul^n_c(t)\cdot F_\LEc\qquad\text{for all $t\in\R$.}
\end{equation}
\end{thm}
\begin{proof} Formula \eqref{5.45} follows by \eqref{4.45}, \eqref{5.22} and \eqref{5.43}. Moreover by \eqref{4.34}, \eqref{5.14} and \eqref{5.43} we have $F_\EcL D^{n-1}_\Eulc=D^{n-1}_\Lagr$, so formula \eqref{5.46} follows by \eqref{4.48} and Theorem~\ref{Theorem 5.2}--II).
\end{proof}

To complete the study of the relation between $(\Eulc)$ and $(\Lagr)$ we also point out that we can  simply combine Corollary~\ref{Corollary 1.2} with\eqref{4.75} and Theorem~\ref{Theorem 5.4} as follows.
\begin{cor}[\bf Regularity for $(\Eulc) \rightleftarrows (\Lagr)$]\label{Corollary 5.2}
For all $n\in\widetilde{\N}$, $2\le n\le r$, the bijective isomorphisms $\Psi_\LEc$ and $\Psi_\EcL$ in Corollary~\ref{Corollary 1.2} respectively restrict to bijective isomorphisms
\begin{equation}\label{5.47}
  \dot{\Psi}_\EcL^n\in\mathcal{L}(\mathcal{S}^n_\Lagr;\mathcal{S}^n_\Eulc),
  \qquad   \dot{\Psi}_\EcL^n\in\mathcal{L}(\mathcal{S}^n_\Eulc;\mathcal{S}^n_\Lagr).
\end{equation}
\end{cor}

By combining Theorem~\ref{Theorem 5.4} with Corollary~\ref{Corollary 5.2} one gets that the commutative diagrams in ~\eqref{1.41}  extend to classes of more regular solutions.
\def\cprime{$'$}
\providecommand{\bysame}{\leavevmode\hbox to3em{\hrulefill}\thinspace}
\providecommand{\MR}{\relax\ifhmode\unskip\space\fi MR }
\providecommand{\MRhref}[2]{%
  \href{http://www.ams.org/mathscinet-getitem?mr=#1}{#2}
}
\providecommand{\href}[2]{#2}

\end{document}